\documentclass[12pt]{article}

\usepackage[latin1]{inputenc}
\usepackage{amsmath,amsthm,amsfonts,amscd,amssymb,xspace,color,rotating}
\usepackage[french]{babel}

\newcommand{\GL}[2]{{\rm GL}_{#1}(#2)}

\newcommand{\td}[2]{\xrightarrow[#1\rightarrow #2]{}}

\newcommand{\N}[1]{\left\|#1\right\|}

\newcommand{\abs}[1]{\left|#1\right|}
\newcommand{\de}{{\rm d}}

\def\Aut{\qopname\relax o{Aut}}

\newtheorem{Thm}{Th\'eor\`eme}[section]
\newtheorem{Prop}[Thm]{Proposition}
\newtheorem{Lem}[Thm]{Lemme}
\newtheorem{Cor}[Thm]{Corollaire}

\newtheorem*{Thm*}{Th\'eor\`eme}
\newtheorem*{Cor*}{Corollaire}
\newtheorem*{Prop*}{Proposition}

\newcounter{ploum}

\newcounter{ex}[section]
\newcounter{rem}[section]

\numberwithin{equation}{section}

\begin{document}

\newenvironment{rem}
{\addtocounter{rem}{1}\setlength{\topsep}{1em}\par\trivlist\item{\bf Remarque
\arabic{section}.\arabic{rem}.} }{\endtrivlist}

\newenvironment{ex}
{\addtocounter{ex}{1}\setlength{\topsep}{1em}\par\trivlist\item{\bf Exemple
\arabic{section}.\arabic{ex}.} }{\endtrivlist}

\newenvironment{demo}[1][D\'emonstration]
{\setlength{\topsep}{1em}\par\trivlist\item{\em #1.} }
{\openbox\endtrivlist}

\newenvironment{cond}{\begin{list}{{\it (\roman{ploum})}}{\usecounter{ploum}}}
{\end{list}}

\newenvironment{souscond}{\begin{list}{{\hspace{5ex}\it (\alph{ploum})}}{\usecounter{ploum}}}
{\end{list}}

\begin{center}

{\Huge\bf Analyse harmonique sur le graphe de Pascal}

\ \\

{\Large\bf J.-F. Quint}

\ \\

\end{center}

\section{Introduction}

\begin{figure}\begin{center}\input{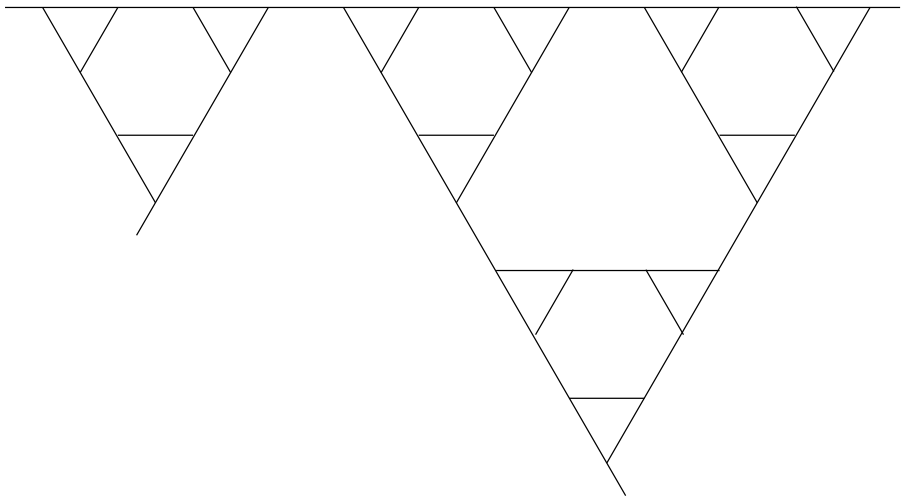}\caption{Le graphe de Pascal}
\label{graphep}\end{center}\end{figure}

Dans tout cet article, nous appellerons graphe de Pascal, et nous
noterons $\Gamma$, le graphe infini, connexe et régulier de valence
$3$ représenté par la figure \ref{graphep}. Ce graphe peut se
construire de la fa\c con suivante. On écrit le triangle de Pascal
et on en efface les valeurs paires des coefficients du binôme. Dans
ce dessin, on joint chaque point à ceux de ses voisins les plus
proches qui n'ont pas été effacés. On obtient ainsi un graphe dans
lequel tout point a trois voisins, sauf le sommet du triangle, qui
en a deux. On prend alors deux copies de ce graphe, qu'on joint par
leurs sommets : on obtient bien un graphe régulier de valence $3$.
C'est le graphe $\Gamma$.

Soit $\varphi$ une fonction de $\Gamma$ dans $\mathbb C$. Pour $p$
dans $\Gamma$, on pose  $\Delta\varphi(p)=\sum_{q\sim p}\varphi(q)$.
Alors, comme  tout point de $\Gamma$ a exactement trois voisins,
l'o\-pé\-ra\-teur linéaire $\Delta$ est auto-adjoint pour la mesure
de comptage sur $\Gamma$, c'est-à-dire que, pour toutes fonctions
$\varphi$ et $\psi$ à support fini, on a
$\sum_{p\in\Gamma}\varphi(p)(\Delta\psi(p))=\sum_{p\in\Gamma}(\Delta\varphi(p))\psi(p)$.
Dans cet article, nous allons complètement déterminer les invariants
spectraux de l'opérateur $\Delta$ dans l'espace $\ell^2(\Gamma)$ des
fonctions de carré intégrable sur $\Gamma$.

Pour énoncer nos résultats, notons $f:\mathbb R\rightarrow\mathbb
R,x\mapsto x^2-x-3$. Soit $\Lambda$ l'ensemble de Julia de $f$,
c'est-à-dire, dans ce cas, l'ensemble des $x$ dans $\mathbb R$ pour
lesquels la suite $(f^n(x))_{n\in\mathbb N}$ reste bornée. C'est un
ensemble de Cantor contenu dans l'intervalle $[-2,3]$. Plus
précisément, si on pose
$I_{-2}=\left[-2,\frac{1-\sqrt{5}}{2}\right]$ et
$I_3=\left[\frac{1+\sqrt{5}}{2},3\right]$, pour tout
$\varepsilon=(\varepsilon_n)_{n\in\mathbb N}$ dans
$\{-2,3\}^{\mathbb N}$, il existe un unique $x$ dans $\Lambda$ tel
que, pour tout $n$ dans $\mathbb N$, on ait $f^n(x)\in
I_{\varepsilon_n}$ et l'application $\{-2,3\}^{\mathbb N}\rightarrow
\Lambda$ ainsi définie est un homéomorphisme bi-höldérien qui
conjugue $f$ et l'application de décalage dans $\{-2,3\}^{\mathbb
N}$. Pour $x$ dans $\Lambda$, posons $\rho(x)=\frac{x}{2x-1}$ et, si
$\varphi$ est une fonction continue sur $\Lambda$,
$L_\rho\varphi(x)=\sum_{f(y)=x}\rho(y)\varphi(y)$. On vérifie
aisément qu'on a $L_\rho(1)=1$. Alors, d'après le théorème de
Ruelle-Perron-Frobenius (voir \cite[§ 2.2]{PP}), il existe une unique
mesure borélienne de probabilité $\nu_\rho$ sur $\Lambda$ telle que
$L_\rho^*\nu_\rho=\nu_\rho$. La mesure $\nu_\rho$ est diffuse et
$f$-invariante. Enfin, on remarque que, si $h$ désigne la fonction
$\Lambda\rightarrow\mathbb R,x\mapsto 3-x$, on a $L_\rho h=2$ et,
donc, $\int_{\Lambda}h\de \nu_\rho=2$.

Notons $p_0$ et $p_0^\vee$ les deux sommets des triangles infinis
qu'on a recollés pour construire le graphe $\Gamma$. Soit
$\varphi_0$ la fonction sur $\Gamma$ qui vaut $1$ en $p_0$, $-1$ en
$p_0^\vee$ et $0$ partout ailleurs. Nous avons le

\begin{Thm} \label{spectrePascal}
Le spectre de $\Delta$ dans $\ell^2(\Gamma)$ est constitué de la
réunion de $\Lambda$ et de l'ensemble $\bigcup_{n\in\mathbb
N}f^{-n}(0)$. La mesure spectrale de $\varphi_0$ pour $\Delta$ dans
$\ell^2(\Gamma)$ est la mesure $h\nu_\rho$, les valeurs propres de
$\Delta$ dans $\ell^2(\Gamma)$ sont les éléments de
$\bigcup_{n\in\mathbb N}f^{-n}(0)$ et de $\bigcup_{n\in\mathbb
N}f^{-n}(-2)$ et les sous-espaces propres associés sont engendrés
par des fonctions à support fini. Enfin, l'orthogonal de la somme
des sous-espace propres de $\Delta$ dans $\ell^2(\Gamma)$ est le
sous-espace cyclique engendré par $\varphi_0$.\end{Thm}

\begin{figure}\begin{center}\input{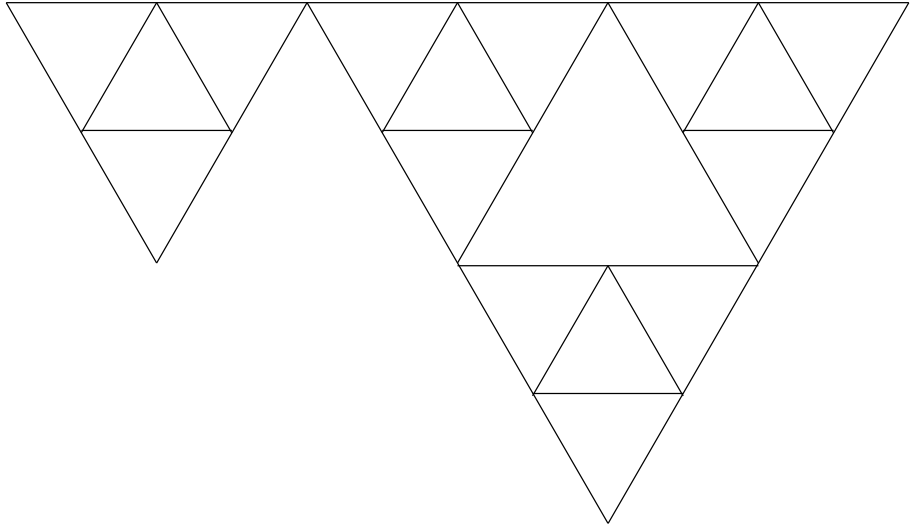}\caption{Le graphe de Sierpi{\'n}ski}
\label{graphes}\end{center}\end{figure}

Un problème semblable a été abordé par A. Teplyaev dans \cite{Tep},
qui a étudié le graphe de Sierpi{\'n}ski, représenté par la figure
\ref{graphes}.
Le graphe de Sierpi{\' n}ski peut être vu comme le graphe des arêtes
du graphe de Pascal, où l'on joint deux arêtes quand elles ont un
point commun. En particulier, notre description du spectre de
$\Delta$ découle des travaux de Teplyaev. En revanche, la
description exacte des composantes cycliques de $\Delta$ et de son
spectre continu sont nouvelles et répondent à la question posée par
Teplyaev dans \cite[§ 6.6]{Tep}. \`A la section \ref{PascalSierp},
nous expliquerons précisément comment faire le lien entre l'étude
spectrale du graphe de Sierpi{\' n}ski et celle du graphe de Pascal

\begin{figure}\begin{center}\input{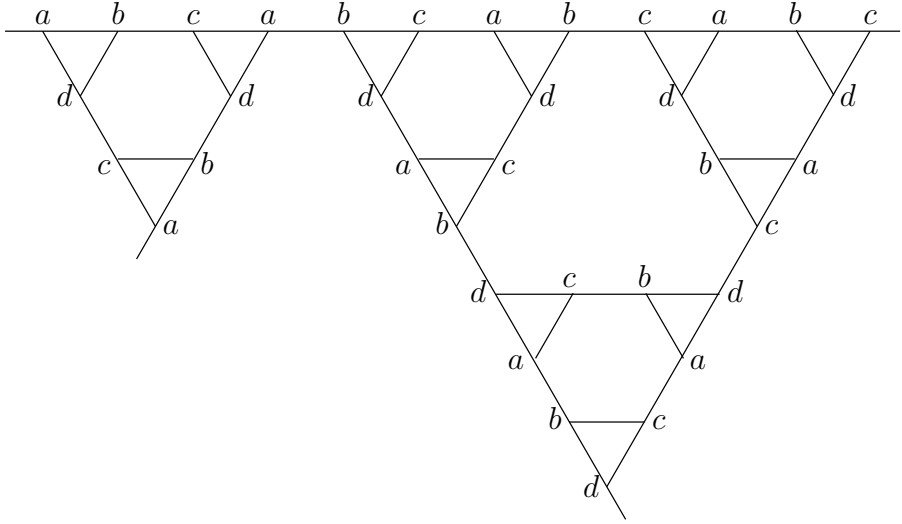}\caption{Un revêtement $\Gamma\rightarrow\Gamma_0$}
\label{revetement}\end{center}\end{figure}

Les méthodes développées dans cet article permettent de décrire la
théorie spectrale d'autres opérateurs, liés au graphe $\Gamma$.
Notons $\Gamma_0$ le graphe complet à quatre sommets $a$, $b$, $c$
et $d$. Le graphe $\Gamma$ est un revêtement du graphe $\Gamma_0$,
ainsi que le montre la figure \ref{revetement}.
Construisons, pour tout entier $n$ un graphe fini de la fa\c con
suivante : si le graphe $\Gamma_n$ a été construit, le graphe
$\Gamma_{n+1}$ est le graphe obtenu en rempla\c cant chaque point de
$\Gamma_n$ par un triangle (ce procédé est détaillé plus
formellement à la section \ref{decritgraphes}). On note toujours
$\Delta$ l'opérateur de somme sur les voisins, agissant sur les
fonctions définies sur $\Gamma_n$. Posons, pour $x$ dans $\mathbb
R$, $k(x)=x+2$, $l(x)=x$ et $m(x)=x-2$. Nous montrerons le

\begin{Thm} \label{spectrePascalfini} Pour tous entiers naturels
$m\geq n$, il existe des revêtements $\Gamma_m\rightarrow\Gamma_n$
et $\Gamma\rightarrow\Gamma_n$. Le polynôme caractéristique de
$\Delta$ dans $\ell^2(\Gamma_n)$ est
$$(X-3)(X+1)^3\prod_{p=0}^{n-1} (m\circ f^p(X))^3(l\circ f^p(X))^{2.3^{n-1-p}}(k\circ f^p(X))^{1+2.3^{n-1-p}}.$$
\end{Thm}

Venons-en à la motivation initiale de cet article, qui était l'étude
d'un phénomène de systèmes dynamiques. Soit $X\subset\left(\mathbb
Z/2\mathbb Z\right)^{\mathbb Z^2}$ le système à trois points,
c'est-à-dire l'ensemble des familles $(p_{k,l})_{(k,l)\in\mathbb
Z^2}$ d'éléments de $\mathbb Z/2\mathbb Z$ telles que, pour tous
entiers relatifs $k$ et $l$, on ait $p_{k,l}+p_{k+1,l}+p_{k,l+1}=0$
(dans $\mathbb Z/2\mathbb Z$). On munit $X$ de l'action naturelle de
$\mathbb Z^2$, engendrée par les applications $T:(x_{k,l})\mapsto
(x_{k+1,l})$ et $S:(x_{k,l})\mapsto (x_{k,l+1})$. Ce système est un
analogue de l'extension naturelle de l'action du doublement et du
triplement de l'angle sur le cercle et, comme dans la conjecture de
Fürstenberg, le problème de la classification des probabilités
boréliennes sur $X$ qui sont $\mathbb Z^2$-invariantes est ouvert.
On note $Y$ l'ensemble des $p$ dans $X$ tels que $p_{0,0}=1$. Alors,
si $p$ est un point de $Y$, il existe exactement trois éléments
$(k,l)$ de l'ensemble $\{(1,0),(0,1),(-1,1),(-1,0),(0,-1),(1,-1)\}$
tels que $T^kS^l x$ appartienne à $Y$. Cette relation munit l'espace
$Y$ d'une structure de graphe régulier de valence $3$ (avec des
arêtes multiples). Si $p$ est un point de $Y$, sa composante connexe
$Y_p$ dans cette structure de graphe est exactement l'ensemble des
points de l'orbite de $p$ sous l'action de $\mathbb Z^2$ qui
appartiennent à $Y$, c'est-à-dire la classe d'équivalence de $p$
dans la relation d'équivalence induite sur $Y$ par l'action de
$\mathbb Z^2$. Pour toute fonction continue $\varphi$ sur $Y$, on
pose, pour tout $p$ dans $Y$,
$$\bar{\Delta}\varphi(p)=\sum_{\substack{(k,l)\in\{(1,0),(0,1),(-1,1),(-1,0),(0,-1),(1,-1)\}
\\ T^kS^lp\in Y}}\varphi(T^kS^lp).$$
Si $\lambda$ est une probabilité borélienne invariante par l'action
de $\mathbb Z^2$ sur $X$ telle que $\lambda(Y)>0$ (c'est-à-dire
telle que $\lambda$ ne soit pas la masse de Dirac en la famille
nulle), la restriction $\mu$ de $\lambda$ à $Y$ vérifie
$\bar{\Delta}^*\mu=3\mu$ et $\bar{\Delta}$ est un opérateur
auto-adjoint de ${\rm L}^2(Y,\mu)$.

\`A l'origine de ce travail, nous souhaitions nous intéresser aux
phénomènes d'intersection homocline de $X$. Rappelons que, si $\phi$
est un difféomorphisme d'une variété compacte $M$ et si $p$ est un
point fixe hyperbolique de $f$, une intersection homocline est un
point d'intersection $q$ de la feuille stable de $p$ et de sa
feuille instable. Pour un tel point, on a en particulier $\phi^n
(q)\td{n}{\pm\infty}p$. Cette notion possède un analogue en
dynamique symbolique. Notons $M$ l'espace $\left(\mathbb Z/2\mathbb
Z\right)^{\mathbb Z}$, $\phi$ l'application de décalage et $p$ le
point de $M$ dont toutes les composantes sont nulles. Si $q$ est un
élément de $M$ dont toutes les coordonnées sauf un nombre fini sont
nulles, on a $\phi^n (q)\td{n}{\pm\infty}p$. En particulier, le
point $q$ tel que $q_0=1$ et dont toutes les autres composantes sont
nulles possède cette propriété. Dans notre situation, on vérifie
qu'il existe un unique élément $q$ de $X$ tel qu'on ait $q_{0,0}=1$,
que, pour tous $k$ dans $\mathbb Z$ et $l\geq 1$, on ait $q_{k,l}=0$
et que, pour tous $k\geq 1$ et $l\leq -1-k$, on ait $q_{k,l}=0$.
Alors, il existe un isomorphisme du graphe $\Gamma$ dans $Y_q$
envoyant $p_0$ sur $q$ : c'est l'origine de la représentation plane
de $\Gamma$ donnée à la figure \ref{graphep}. Dorénavant, on
identifiera $q$ à $p_0$ et $Y_q$ à $\Gamma$. Notons que, si $p$
désigne l'élément de $X$ dont toutes les composantes sont nulles,
pour tous entiers $l>k>0$, on a $(T^{-k}S^l)^n
(p_0)\td{n}{\pm\infty}p$. Soit $\bar{\Gamma}$ l'adhérence de
$\Gamma$ dans $Y$ : on peut voir $\bar{\Gamma}$ comme un ensemble de
graphes pointés plans. Notre objectif est de déterminer les
structures induites sur $\bar{\Gamma}$ par l'action de $\mathbb Z^2$
sur $X$. Pour tout $p$ dans $\bar{\Gamma}$, on notera $\Gamma_p$
pour $Y_p$. Notons qu'alors, si $\Theta_p$ désigne le graphe des
arêtes de $\Gamma_p$, les graphes $\Theta_p$ sont exactement les
graphes étudiés par Teplyaev dans \cite{Tep}. En particulier,
d'après \cite[§ 5.4]{Tep}, si $\Gamma_p$ ne contient pas $p_0$ où
l'une de ses six images par l'action naturelle du groupe diédral
d'ordre $6$ sur l'espace $X$, le spectre de l'opérateur $\Delta$
dans $\ell^2(\Gamma_p)$ est discret.

Tout point $p$ de $\bar{\Gamma}$ appartient à un unique triangle
dans $\Gamma_p$. Notons $a$ l'ensemble des éléments $p$ de
$\bar{\Gamma}$ pour lesquels ce triangle est
$\{p,T^{-1}p,S^{-1}p\}$, $b$ celui des points $p$ pour lesquels il
est de la forme $\{p,Tp,TS^{-1}p\}$ et $c$ l'ensemble des points $p$
pour lesquels ce triangle est $\{p,Sp,T^{-1}Sp\}$. L'ensemble
$\bar{\Gamma}$ est la réunion disjointe de $a$, $b$ et $c$. Notons
$\theta_1:\bar{\Gamma}\rightarrow\{a,b,c\}$ l'application naturelle
associée à cette partition. Nous dirons qu'une fonction $\varphi$ de
$\bar{\Gamma}$ dans $\mathbb C$ est $1$-triangulaire si elle
factorise à travers $\theta_1$. Nous noterons $E_1$ l'espace des
fonctions $1$-triangulaires $\varphi$ telles que
$\varphi(a)+\varphi(b)+\varphi(c)=0$ : il s'identifie naturellement
à $\mathbb C^3_0=\{(s,t,u)\in\mathbb C^3|s+t+u=0\}$. On munira
$\mathbb C^3_0$ du produit scalaire égal à un tiers du produit
canonique.

Notons $\zeta:\Lambda\rightarrow\mathbb R_+^*,x\mapsto
\frac{1}{3}\frac{(x+3)(x-1)}{2x-1}$ et, comme ci-dessus, désignons
par $L_\zeta$ l'opérateur de transfert associé à $\zeta$ pour la
dynamique du polynôme $f$. Comme on a $L_\zeta(1)=1$, il existe une
unique mesure borélienne de probabilité $\nu_\zeta$ sur $\Lambda$
telle que $L_\zeta^*\nu_\zeta=\nu_\zeta$. Alors, si $j$ désigne la
fonction $\Lambda\rightarrow \mathbb
R,x\mapsto\frac{1}{3}\frac{3-x}{x+3}$, on a $L_\zeta(j)=1$ et, donc,
$\int_{\Lambda}j\de \nu_\zeta=1$.

\begin{Thm} \label{spectrePascalcomp}
Pour tout $p$ dans $\bar{\Gamma}$, l'ensemble $\Gamma_p$ est dense
dans $\bar{\Gamma}$. Il existe une unique probabilité borélienne
$\mu$ sur $\bar{\Gamma}$ telle que $\bar{\Delta}^*\mu=3\mu$ et
l'opérateur $\bar{\Delta}$ est auto-adjoint dans ${\rm
L}^2\left(\bar{\Gamma},\mu\right)$. Le spectre de l'opérateur
$\bar{\Delta}$ dans ${\rm L}^2\left(\bar{\Gamma},\mu\right)$ est le
même que celui de $\Delta$ dans $\ell^2(\Gamma)$. Pour tout
$\varphi$ dans $E_1$, la mesure spectrale de $\varphi$ pour
$\bar{\Delta}$ dans ${\rm L}^2\left(\bar{\Gamma},\mu\right)$ est
$\N{\varphi}^2_2j\nu_\zeta$ et la somme des sous-espaces cycliques
engendrés par les éléments de $E_1$ est isométrique à ${\rm
L}^2\left(j\nu_\zeta,\mathbb C^3_0\right)$. Le spectre de
$\bar{\Delta}$ dans l'orthogonal de ce sous-espace est discret et
ses valeurs propres sont $3$, qui est simple, et les éléments de
$\bigcup_{n\in\mathbb N}f^{-n}(0)\cup\bigcup_{n\in\mathbb
N}f^{-n}(-2)$.\end{Thm}

Le plan de l'article est le suivant.

Les sections \ref{decritgraphes}, \ref{secspectre},
\ref{mesurespectrale}, \ref{valeurp} et \ref{decompo} sont
consacrées à l'étude du graphe $\Gamma$. Dans la section
\ref{decritgraphes}, nous construisons précisément $\Gamma$ et
établissons des propriétés élémentaires de sa géométrie. Dans la
section \ref{secspectre}, nous déterminons le spectre de $\Delta$
dans $\ell^2(\Gamma)$ et, dans la section \ref{mesurespectrale},
nous démontrons un résultat essentiel en vue du calcul des mesures
spectrales des éléments de cet espace. Dans la section
\ref{valeurp}, nous décrivons la structure des espaces propres de
$\Delta$ dans $\ell^2(\Gamma)$. Enfin, à la section \ref{decompo},
nous appliquons l'ensembles de ces résultats préliminaires à la
démonstration du théorème \ref{spectrePascal}.

\`A la section \ref{graphesfinis}, nous utilisons les techniques
mises aux points précédem\-ment pour démontrer le théorème
\ref{spectrePascalfini}.

Dans les sections \ref{compactifie}, \ref{fonctiontriangulaire},
\ref{deltacomp}, \ref{secspectrecomp}, \ref{valeurpcomp} et
\ref{decompocomp}, nous étudions l'espace $\bar{\Gamma}$. \`A la
section \ref{compactifie}, nous décrivons en détail la géométrie de
l'espace $\bar{\Gamma}$ et, à la section \ref{fonctiontriangulaire},
nous introduisons des espaces de fonctions localement constantes
remarquables sur cet espace. La section \ref{deltacomp} est
consacrée à la définition de l'opérateur $\bar{\Delta}$ et à la
démonstration de l'unicité de sa mesure harmonique. La section
\ref{secspectrecomp} étend à $\bar{\Gamma}$ les propriétés
démontrées pour $\Gamma$ dans les sections \ref{secspectre} et
\ref{mesurespectrale}. Dans la section \ref{valeurpcomp}, nous
étudions les espaces propres de $\bar{\Delta}$ dans ${\rm
L}^2\left(\bar{\Gamma},\mu\right)$. Enfin, dans la section
\ref{decompocomp}, nous achevons la démonstration du théorème
\ref{spectrePascalcomp}.

Dans la section \ref{PascalSierp}, nous expliquons rapidement
comment transférer nos résultats du graphe de Pascal au graphe de
Sierpi\'nski.

\section{Préliminaires géométriques}
\label{decritgraphes}

Dans toute la suite, nous appellerons graphe un ensemble  $\Phi$
muni d'une relation symétrique $\sim$ telle que, pour tout $p$ dans
$\Phi$, on n'ait pas $p\sim p$. Pour $p$ dans $\Phi$, nous
appellerons voisins de $p$ l'ensemble des éléments $q$ de $\Phi$
tels que $p\sim q$. Nous dirons que $\Phi$ est régulier si tous ses
éléments ont le même nombre (fini) de voisins. Dans ce cas, nous
appellerons ce nombre la valence de $\Phi$. Nous dirons que $\Phi$
est connexe si, pour tous $p$ et $q$ dans $\Phi$, il existe une
suite de points $r_0=p,r_1,\ldots,r_n=q$ de $\Phi$ telle que, pour
tout $1\leq i\leq n$, on ait $r_{i-1}\sim r_{i}$. Nous appellerons
une telle suite un chemin de $p$ à $q$ et l'entier $n$ la longueur
de ce chemin. Si $\Phi$ est connexe et si $\varphi$ est une fonction
sur $\Phi$ telle que, pour tous points $p$ et $q$ de $\Phi$ avec
$p\sim q$ on ait $\varphi(p)=\varphi(q)$, $\varphi$ est constante.

Nous dirons qu'une partie $\mathcal T$ d'un graphe $\Phi$ est un
triangle si $\mathcal T$ est constituée d'exactement trois points
$p$, $q$ et $r$ et qu'on a $p\sim q$, $q\sim r$ et $r\sim p$.

\begin{figure}\begin{center}\input{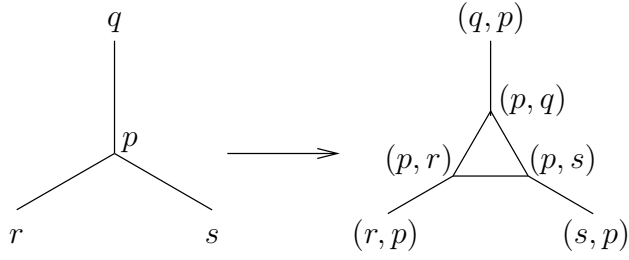}
\caption{Construction du graphe $\hat{\Phi}$}\label{nouveaugraphe}\end{center}\end{figure}

Soit $\Phi$ un graphe régulier de valence $3$. On note $\hat{\Phi}$
l'ensemble des couples $(p,q)$ d'éléments de $\Phi$ avec $p\sim q$
et on le munit de la structure de graphe pour laquelle, si $p$ est
un point de $\Phi$, de voisins $q$, $r$ et $s$, les voisins de
$(p,q)$ sont $(q,p)$, $(p,r)$ et $(p,s)$. Si $\Phi$ est connexe,
$\hat{\Phi}$ est connexe. Géométriquement, $\hat{\Phi}$ est le
graphe obtenu en rempla\c cant chaque point de $\Phi$ par un
triangle. Ce procédé est représenté à la figure \ref{nouveaugraphe}.
On note $\Pi$ l'application $\hat{\Phi}\rightarrow\Phi,(p,q)\mapsto
p$.

Notons $\ell^2(\Phi)$ l'espace des fonctions
$\varphi:\Phi\rightarrow\mathbb C$ telles que
$\sum_{p\in\Phi}\abs{\varphi(p)}^2<\infty$, muni de sa structure
naturelle d'espace de Hilbert $\langle .,.\rangle$. Si $\Phi$ est
régulier de valence $3$, l'application $\Pi$ induit une application
linéaire bornée de norme $3$,
$\Pi^*:\varphi\mapsto\varphi\circ\Pi,\ell^2(\Phi)\rightarrow
\ell^2\left(\hat{\Phi}\right)$. Par abus de langage, on note encore
$\Pi$ l'adjoint de $\Pi^*$ : c'est l'opérateur linéaire borné
$\ell^2\left(\hat{\Phi}\right)\rightarrow\ell^2(\Phi)$ qui, à une
fonction $\varphi$ dans $\ell^2\left(\hat{\Phi}\right)$, associe la
fonction dont la valeur en un point $p$ de $\Phi$ est $\sum_{q\sim
p}\varphi(p,q)$. On a $\Pi\Pi^*=3$.

\'Etendons la définition du graphe des triangles à des graphes plus
gé\-né\-raux. Nous dirons qu'un graphe $\Phi$ est régulier de
valence $3$ à bord si tous les points de $\Phi$ ont deux ou trois
voisins. Dans ce cas, nous appellerons bord de $\Phi$ et nous
noterons $\partial \Phi$ l'ensemble des points de $\Phi$ qui ont
deux voisins. Si $\Phi$ est un graphe régulier de valence $3$ à
bord, nous noterons $\hat{\Phi}$ l'ensemble formé de la réunion de
$\partial \Phi $ et de l'ensemble des couples $(p,q)$ d'éléments de
$\Phi$ avec $p\sim q$. Nous munirons $\hat{\Phi}$ de la structure de
graphe pour laquelle, si $p$ est un point de $\Phi-\partial \Phi$,
de voisins $q$, $r$ et $s$, les voisins de $(p,q)$ sont $(q,p)$,
$(p,r)$ et $(p,s)$ et, si $p$ est un point de $\partial\Phi$, de
voisins $q$ et $r$, les voisins de $p$ dans $\hat{\Phi}$ sont
$(p,q)$ et $(p,r)$ et les voisins de $(p,q)$ sont $(q,p)$, $p$ et
$(p,r)$. Ainsi, $\hat{\Phi}$ est à nouveau un graphe régulier de
valence $3$ à bord et le bord de $\hat{\Phi}$ est en bijection
naturelle avec celui de $\Phi$. \`A nouveau, si $\Phi$ est connexe,
$\hat{\Phi}$ est connexe. On note encore $\Pi$ l'application
naturelle  $\hat{\Phi}\rightarrow\Phi$ et $\Pi^*$ et $\Pi$ les
opérateurs bornés associés $\ell^2(\Phi)\rightarrow
\ell^2\left(\hat{\Phi}\right)$ et
$\ell^2\left(\hat{\Phi}\right)\rightarrow\ell^2(\Phi)$.

\begin{Lem} \label{triangles}
Soit $\Phi$ un graphe régulier de valence $3$ à bord. Alors, les
triangles de  $\hat{\Phi}$ sont exactement les parties de la forme
$\Pi^{-1}(p)$ où $p$ est un point de  $\Phi$. En particulier, tout
point de  $\hat{\Phi}$ appartient à un unique triangle.\end{Lem}

\begin{demo} Si $p$ est un point de $\Phi$, l'ensemble $\Pi^{-1}(p)$ constitue clairement un triangle.
Réciproquement, donnons-nous un point $p$ de $\Phi-\partial\Phi$, de
voisins $q$, $r$ et $s$. Alors, les voisins de $(p,q)$ sont $(q,p)$,
$(p,r)$ et $(p,s)$. Par définition, comme $p\neq q$, le point
$(q,p)$ ne peut pas être voisin de $(p,r)$ ou de $(p,s)$. Par
conséquent, le seul triangle contenant $(q,p)$ est $\Pi^{-1}(p)$. De
même, si $p$ appartient à $\partial\Phi$, et si les voisins de $p$
sont $q$ et $r$, comme $p$ n'a que deux voisins dans  $\hat{\Phi}$,
$p$ n'appartient qu'à un seul triangle et, comme $p\neq q$, aucun
des voisins de $(q,p)$ n'est voisin de $(p,q)$, donc $(p,q)$
n'appartient qu'à un seul triangle.\end{demo}

Si $\Phi$ est un graphe, nous dirons qu'une bijection
$\sigma:\Phi\rightarrow\Phi$ est un automorphisme du graphe $\Phi$
si, pour tous $p$ et $q$ dans $\Phi$ avec $p\sim q$, on a
$\sigma(p)\sim\sigma(q)$. L'ensemble des automorphismes de $\Phi$
constitue un sous-groupe du groupe des bijections de $\Phi$, qu'on
note $\Aut\Phi$. Si $\Phi$ est régulier de valence $3$ à bord et si
$\sigma$ est un automorphisme de $\Phi$, on a
$\sigma(\partial\Phi)=\partial\Phi$ et il existe un unique
automorphisme $\hat{\sigma}$ de $\hat{\Phi}$ tel que
$\Pi\hat{\sigma}=\sigma\Pi$.

\begin{Lem} \label{automorphismes}
Soit $\Phi$ un graphe régulier de valence $3$ à bord.
L'application $\sigma\mapsto\hat{\sigma},\Aut\Phi\rightarrow\Aut\hat{\Phi}$ est un isomorphisme de groupes.\end{Lem}

\begin{demo} Cette application étant clairement un morphisme injectif, il suffit
de démontrer qu'elle est surjective. Soit donc $\tau$ un
automorphisme de $\hat{\Phi}$. Comme $\tau$ permute les triangles de
$\hat{\Phi}$, d'après le lemme \ref{triangles}, il existe une unique
bijection $\sigma:\Phi\rightarrow\Phi$ telle que
$\Pi\tau=\sigma\Pi$. Soient $p$ et $q$ des points de $\Phi$ tels que
$p\sim q$. Alors, on a $(p,q)\sim (q,p)$, donc $\tau(p,q)\sim
\tau(q,p)$ et, comme ces deux points de $\hat{\Phi}$ n'appartiennent
pas à un même triangle, $\sigma(p)=\Pi\tau(p,q)\sim
\Pi\tau(q,p)=\sigma(q)$. Donc $\sigma$ est un automorphisme de
$\Phi$ et $\tau=\hat{\sigma}$, ce qu'il fallait démontrer.\end{demo}

Nous allons à présent définir une famille importante de graphes
réguliers de valence $3$ à bord. Si $a$, $b$ et $c$ sont trois
éléments distincts, on note $\mathcal T(a,b,c)=\mathcal T_1(a,b,c)$
l'ensemble $\{a,b,c\}$ muni de la structure de graphe pour laquelle
on a $a\sim b$, $b\sim c$ et $c\sim a$ et on dit que $\mathcal
T(a,b,c)$ est le triangle ou le $1$-triangle de sommets $a$, $b$ et
$c$. On le considère comme un graphe régulier de valence $3$ à bord.
On définit alors par récurrence une famille de graphes réguliers à
bord en posant, pour tout $n\geq 1$, $\mathcal
T_{n+1}(a,b,c)=\widehat{\mathcal T_{n}(a,b,c)}$. Pour tout $n\geq
1$, on appelle $\mathcal T_n(a,b,c)$ le $n$-triangle de sommets $a$,
$b$ et $c$.

On note $\mathfrak S(a,b,c)$ le groupe des permutations de
l'ensemble $\{a,b,c\}$. Par définition et d'après le lemme
\ref{automorphismes}, on a le

\begin{Lem} \label{grostriangles}
Soit $a$, $b$ et $c$, trois éléments distincts. Alors, pour tout
$n\geq 1$, $\mathcal T_n(a,b,c)$ est un graphe connexe régulier de
valence $3$ à bord et $\partial \mathcal T_n(a,b,c)=\{a,b,c\}$.
L'application qui, à un automorphisme $\sigma$ de  $\mathcal
T_n(a,b,c)$, associe sa restriction à $\{a,b,c\}$ induit un
isomorphisme de groupes de $\Aut \mathcal T_n(a,b,c)$ dans
$\mathfrak S(a,b,c)$.\end{Lem}

Si $\Phi$ est un graphe et $n$ un entier $\geq 1$, nous dirons
qu'une partie $\mathcal T$ de $\Phi$ est un $n$-triangle s'il existe
des points $p$, $q$ et $r$ de $\mathcal T$ tels que la partie
$\mathcal T$, munie de la restriction de la relation $\sim$, soit
isomorphe au graphe $\mathcal T_n(p,q,r)$. Par abus de langage, nous
appellerons $0$-triangles les points de $\Phi$.

Soit $\Phi$ un graphe régulier de valence $3$ à bord. Notons
$\hat\Phi^{(0)}=\Phi$, $\hat{\Phi}^{(1)}=\hat{\Phi}$ et, pour tout
entier $n$, $\hat{\Phi}^{(n+1)}=\widehat{\hat{\Phi}^{(n)}}$. Par
récurrence, pour tout entier $n$, tout automorphisme $\sigma$ de
$\Phi$ induit un unique automorphisme $\hat{\sigma}^{(n)}$ de
$\hat{\Phi}^{(n)}$ tel que $\Pi^n\hat{\sigma}^{(n)}=\sigma \Pi^n$.

Des lemmes \ref{triangles} et \ref{automorphismes}, on déduit immédiatement, par récurrence, le

\begin{Lem} \label{structuregrostriangles}
Soient $\Phi$ un graphe régulier de valence $3$ à bord et $n$ un
entier naturel. Les $n$-triangles de  $\hat{\Phi}^{(n)}$ sont
exactement les parties de la forme $\Pi^{-n}(p)$ où $p$ est un point
de $\Phi$. En particulier, tout point de $\hat{\Phi}^{(n)}$
appartient à un unique $n$-triangle. L'application
$\sigma\mapsto\hat{\sigma}^{(n)},\Aut\Phi\rightarrow\Aut\hat{\Phi}^{(n)}$
est un isomorphisme de groupes.\end{Lem}

\begin{Cor}\label{sommettriangles}
Soient $\Phi$ un graphe régulier de valence $3$ à bord, $n\geq m\geq
1$ des entiers naturels, $p$ un point de $\hat{\Phi}^{(n)}$,
$\mathcal T$ le $n$-triangle contenant $p$ et $\mathcal S$ le
$m$-triangle contenant $p$. On a $\mathcal S\subset\mathcal T$. Si
$p$ est un sommet de $\mathcal T$ et si $p$ n'appartient pas à
$\partial\hat{\Phi}^{(n)}$, l'unique voisin de $p$ dans
$\hat{\Phi}^{(n)}$ qui n'appartient pas au $1$-triangle contenant
$p$ est lui-même le sommet d'un $n$-triangle de $\hat{\Phi}^{(n)}$.
\end{Cor}

\begin{demo} D'après le lemme \ref{structuregrostriangles}, on a
$\mathcal T=\Pi^{-n}(\Pi^np)$ et $\mathcal S=\Pi^{-m}(\Pi^mp)$, donc
$\mathcal S\subset\mathcal T$. Supposons que $p$ appartient à
$\partial\mathcal T-\partial\hat{\Phi}^{(n)}$. Alors, $p$ possède un
unique voisin dans $\hat{\Phi}^{(n)}$ qui n'appartient pas à
$\mathcal T$. Soient $q$ un voisin de $p$ et $\mathcal R$ le
$n$-triangle contenant $q$. Si $q$ n'est pas un sommet de $\mathcal
R$, tous les voisins de $\mathcal Q$ appartiennent à $\mathcal R$
et, donc, on a $p\in\mathcal R$. Comme, d'après le lemme
\ref{structuregrostriangles}, $p$ appartient à un unique
$n$-triangle de $\hat{\Phi}^{(n)}$, on a $\mathcal T=\mathcal R$,
et, donc, $q$ appartient à $\mathcal T$, si bien que le voisin de
$p$ qui n'appartient pas à $\mathcal T$ est un sommet du
$n$-triangle qui le contient.
\end{demo}

\begin{Cor} \label{decoupetriangles}
Soient $n$ un entier $\geq 2$ et $a$, $b$ et $c$ des éléments
distincts. Alors, il existe des éléments uniques $ab$, $ba$, $ac$,
$ca$, $bc$ et $cb$ de $\mathcal T_n(a,b,c)$ tels que $\mathcal
T_n(a,b,c)$ soit la réunion des trois $(n-1)$-triangles $\mathcal
T_{n-1}(a,ab,ac)$, $\mathcal T_{n-1}(b,ba,bc)$ et $\mathcal
T_{n-1}(c,ca,cb)$ et que l'on ait $ab\sim ba$, $ac\sim ca$ et
$bc\sim cb$.
\end{Cor}

\begin{demo} Pour $n=2$, le corollaire se démontre directement. On en déduit
le cas général en appliquant le lemme \ref{structuregrostriangles}.\end{demo}

Donnons-nous à présent un élément $a$ et deux suites d'éléments
distincts $(b_n)_{n\geq 1}$ et $(c_n)_{n\geq 1}$ telles que, pour
tous $n\geq 1$, on ait $b_n\neq a$, $c_n\neq a$ et $b_n\neq c_n$.
D'après le corollaire \ref{decoupetriangles}, pour tout $n\geq 1$,
on peut identifier $\mathcal T_n(a,b_n,c_n)$ à une partie de
$\mathcal T_{n+1}(a,b_{n+1},c_{n+1})$ grâce à l'unique isomorphisme
de graphes envoyant $a$ sur $a$, $b_n$ sur $ab_{n+1}$ et  $c_n$ sur
$ac_{n+1}$. On appelle alors triangle infini issu de $a$ et on note
$\mathcal T_{\infty}(a)$ l'ensemble $\bigcup_{n\geq 1}\mathcal
T_n(a,b_n,c_n)$ muni de la structure de graphe qui induit sur chacun
des $\mathcal T_n(a,b_n,c_n)$, $n\geq 1$, sa structure de
$n$-triangle. Des études précédentes, on déduit le

\begin{Lem}\label{triangleinfini}
Soit $a$ un élément. Le graphe $\mathcal T_{\infty}(a)$ est connexe,
régulier de valence $3$ à bord et $\partial \mathcal
T_{\infty}(a)=\{a\}$. Si $b$ et $c$ sont les deux voisins de $a$
dans $\mathcal T_{\infty}(a)$, il existe un unique isomorphisme de
$\mathcal T_{\infty}(a)$ dans $\widehat{\mathcal T_{\infty}(a)}$ qui
envoie $a$ sur $a$, $b$ sur $(a,b)$ et $c$ sur $(a,c)$. Pour tout
entier naturel $n$, cet isomorphisme induit une bijection naturelle
entre les points de $\mathcal T_{\infty}(a)$ et les $n$-triangles de
$\mathcal T_{\infty}(a)$ et tout point de  $\mathcal T_{\infty}(a)$
appartient à un unique $n$-triangle.
Enfin, $\mathcal T_{\infty}(a)$ possède un
unique automorphisme non-trivial ; cet automorphisme est une
involution qui fixe $a$ et qui, pour tout $n\geq 1$, échange les
deux sommets différents  de $a$ du $n$-triangle contenant
$a$.\end{Lem}

Dans tout cet article, on fixe deux éléments distincts $p_0$ et
$p_0^\vee$. On appelle graphe de Pascal et on note $\Gamma$
l'ensemble $\mathcal T_{\infty}(p_0)\cup \mathcal
T_{\infty}(p_0^\vee)$ muni de la structure de graphe qui induit la
structure de triangle infini sur $\mathcal T_{\infty}(p_0)$ et sur
$\mathcal T_{\infty}(p_0^\vee)$ et pour laquelle on a $p_0\sim
p_0^\vee$. Du lemme \ref{triangleinfini}, on déduit la

\begin{Prop}\label{structurePascal} Le graphe de Pascal est un
graphe infini, connexe et régulier de valence $3$. Si $q_0$ et $r_0$
sont les deux voisins de $p_0$ dans $\mathcal T_{\infty}(p_0)$ et
$q_0^\vee$ et $r_0^\vee$ les deux voisins de $p_0^\vee$ dans
$\mathcal T_{\infty}(p_0^\vee)$, il existe un unique isomorphisme de
$\Gamma$ dans $\hat{\Gamma}$ qui envoie $p_0$ sur $(p_0,p_0^\vee)$,
$p_0^\vee$ sur $(p_0^\vee,p_0)$, $q_0$ sur $(p_0,q_0)$, $r_0$ sur
$(p_0,r_0)$, $q_0^\vee$ sur $(p_0^\vee,q_0^\vee)$ et $r_0^\vee$ sur
$(p_0^\vee,r_0^\vee)$. Pour tout entier naturel $n$, cet
isomorphisme induit une bijection naturelle entre les points de
$\Gamma$ et les $n$-triangles de $\Gamma$ et tout point de $\Gamma$
appartient à un unique $n$-triangle.
\end{Prop}

Une représentation plane du graphe de Pascal est donnée à la figure
\ref{graphep}.

On identifiera dorénavant $\Gamma$ et $\hat{\Gamma}$ par
l'isomorphisme décrit dans la proposition \ref{structurePascal}. En
particulier, on considérera donc désormais que $\Pi^*$ et $\Pi$ sont
des endomorphismes bornés de $\ell^2(\Gamma)$.

On notera $\Theta$ le graphe des arêtes de $\Gamma$. Plus
précisément, $\Theta$ sera constitué de l'ensemble des paires
$\{p,q\}$ d'éléments de $\Gamma$ avec $p\sim q$, muni de la relation
pour laquelle, si $p$ et $q$ sont deux points voisins de $\Gamma$,
si $r$ et $s$ sont les deux autres voisins de $p$ et $t$ et $u$ les
deux autres voisins de $q$, les voisins de $\{p,q\}$ sont $\{p,r\}$,
$\{p,s\}$, $\{q,t\}$ et $\{q,u\}$. On appellera $\Theta$ le graphe
de Sierpi{\'n}ski. C'est un graphe infini, connexe et régulier de
valence $4$. Une représentation plane en est donnée à la figure
\ref{graphes}.

Si $\Phi$ est un graphe régulier de valence $k$, pour toute fonction
$\varphi$ de $\Phi$ dans $\mathbb C$, on notera $\Delta\varphi$ la
fonction $p\mapsto\sum_{q\sim p}\varphi(q)$. Alors, $\Delta$ induit
un opérateur auto-adjoint de norme $\leq k$ de l'espace
$\ell^2(\Phi)$. On appellera spectre de $\Phi$ le spectre de cet
opérateur.

\section{Le spectre de $\Gamma$}
\label{secspectre}

Soit $\Phi$ un graphe régulier de valence $3$. Dans cette section,
nous allons étudier le lien entre les propriétés spectrales de
$\Phi$ et celle de $\hat{\Phi}$. Notre étude est fondée sur le

\begin{Lem} \label{relation1}
On a $(\Delta^2-\Delta-3)\Pi^*=\Pi^*\Delta$ et $\Pi(\Delta^2-\Delta-3)=\Delta\Pi$.\end{Lem}

\begin{demo}
Soient $\varphi$ une fonction sur $\Phi$, $p$ un point de $\Phi$ et
$q,r,s$ les trois voisins de $p$. Supposons que $\varphi(p)=a$,
$\varphi(q)=b$,  $\varphi(r)=c$ et   $\varphi(s)=d$. Alors, on a
$\Pi^*\varphi(p,q)=a$, $\Delta\Pi^*\varphi(p,q)=2a+b$ et
$\Delta^2\Pi^*\varphi(p,q)=(2b+a)+(2a+c)+(2a+d)=5a+2b+c+d$. Il vient
bien
$(\Delta^2-\Delta-3)\Pi^*\varphi(p,q)=b+c+d=\Pi^*\Delta\varphi(p,q)$.
La seconde relation s'obtient en passant aux opérateurs adjoints
dans la première.
\end{demo}

Nous allons utiliser le lemme \ref{relation1} pour déterminer le
spectre de $\Delta$ dans  $\ell^2\left(\hat{\Phi}\right)$. Nous
aurons recours à des résultats élémentaires d'analyse fonctionnelle.

\begin{Lem}\label{spectral1} Soient $E$ un espace de Banach et $T$ un opérateur borné de $E$.
On suppose que tous les éléments du spectre de $T$ ont une partie
réelle strictement positive. Alors, si $F\subset E$ est un
sous-espace fermé stable par $T^2$, $F$ est stable par $T$.\end{Lem}

\begin{demo} Soient $0<\alpha<\beta$ et
$\gamma>0$ tels que le spectre $S$ de $T$ soit contenu dans
l'intérieur du rectangle $R=[\alpha,\beta]+[-\gamma,\gamma]i$ et
$U\supset R$ et $V$ des ouverts de $\mathbb C$ tels que
l'application $\lambda\mapsto\lambda^2$ induise un bi-holomorphisme
de $U$ dans $V$. Il existe une fonction holomorphe $r$ sur $V$ telle
que, pour tout $\lambda$ dans $U$, on ait $r(\lambda^2)=\lambda$.
Comme $R$ est simplement connexe, d'après le théorème de Runge, il
existe une suite $(r_n)_{n\in\mathbb N}$ de polynômes dans $\mathbb
C[X]$ qui converge uniformément vers $r$ sur $R^2$. Alors, comme le
spectre de $T^2$ est $S^2$, qui est contenu dans l'intérieur de
$R^2$, la suite $r_n(T^2)$ converge vers $T$ dans l'espace des
endomorphismes de $E$. Pour tout entier $n$, $r_n(T^2)$ laisse
stable $F$, donc $T$ laisse stable $F$.
\end{demo}

\begin{Lem}\label{spectral2} Soient $H$ un espace de Hilbert, $T$ un endomorphisme
auto-adjoint borné de $H$ et $\pi$ un polynôme du second degré à
coefficients réels. On suppose qu'il existe un sous-espace fermé $K$
de $H$ tel que $\pi(T)K\subset K$ et que $K$ et $TK$ engendrent $H$.
Alors, l'image par $\pi$ du spectre de $T$ dans $H$ est égale au
spectre de $\pi(T)$ dans $K$ et, si on a en outre $T^{-1}K\cap
K=\{0\}$, le spectre de $T$ dans $H$ est exactement l'ensemble des
$\lambda$ dans $\mathbb R$ tels que $\pi(\lambda)$ appartienne au
spectre de $\pi(T)$ dans $K$.\end{Lem}

\begin{demo} Après avoir écrit $\pi$ sous sa forme canonique,
on peut supposer qu'on a $\pi(X)=X^2$. Notons $E$ la résolution
spectrale de $T$ : pour tout borélien $B$ de $\mathbb R$, $E(B)$ est
un projecteur de $H$ qui commute à $T$. Soit $B$ un borélien de
$\mathbb R$ tel que $B=-B$. Alors, pour toute mesure de Radon $\mu$
sur $\mathbb R$, dans ${\rm L}^2(\mu)$, la fonction caractéristique
de $B$ est la limite d'une suite de polynômes pairs. On a donc
$E(B)K\subset K$.

Le spectre de $T^2$ dans $H$ est exactement l'ensemble des carrés
des éléments du spectre de $T$ dans $H$. Comme $T^2$ est
auto-adjoint et que $K$ est stable par $T^2$, le spectre de $T^2$
dans $K$ est contenu dans son spectre dans $H$, et donc, dans
l'ensemble des carrés des éléments du spectre de $T$.
Réciproquement, supposons qu'il existe des éléments du spectre de
$T$ dont le carré n'appartienne pas au spectre de $T^2$ dans $K$.
Alors, il existe un ouvert $V$ de $\mathbb R$, symétrique et tel que
$V$ contienne des éléments du spectre de $T$, mais que $V^2$ ne
contienne pas d'éléments de celui de $T^2$ dans $K$. On a
$E(V)K\subset K$, mais, comme $V^2$ ne contient pas d'éléments du
spectre de $T^2$ dans $K$, $E(V)K=0$. Or, puisque $K$ et $TK$
engendrent $H$, $E(V)K$ et $TE(V)K=E(V)TK$ engendrent $E(V)H$. Il
vient $E(V)H=0$, ce qui contredit le fait que $V$ contient des
valeurs spectrales de $T$. Le spectre de $T^2$ dans $K$ est donc
exactement l'ensemble des carrés des éléments du spectre de $T$ dans
$H$.

Supposons à présent qu'on a  $T^{-1}K\cap K=\{0\}$. Pour conclure,
il nous reste à prouver que le spectre de $T$ est symétrique.
Supposons que ce ne soit pas le cas. Alors, on peut, quitte à
remplacer $T$ par $-T$, trouver des réels $0<\alpha<\beta$ tel que
$U=]\alpha,\beta[$ contienne des éléments du spectre de $T$ et que
$-U$ n'en contienne pas. Mais alors, on a $E(U)=E(U\cup(-U))$ et,
donc, $E(U)K\subset K$. Si $L$ désigne l'image de $E(U)$ dans $H$,
on a donc $T^2(K\cap L)\subset K\cap L$. Comme le spectre de la
restriction de $T$ à $L$ est contenu dans $\mathbb R_+^*$, on a,
d'après le lemme \ref{spectral1}, $T(K\cap L)\subset K\cap L$ et,
donc, par hypothèse, $K\cap L=0$. Comme $E(U)K\subset K$, il vient
$E(U)=0$ sur $K$. Comme $K$ et $TK$ engendrent $H$, on a $E(U)=0$,
ce qui contredit le fait que $U$ contient des éléments du spectre de
$T$. Le spectre de $T$ est donc symétrique. Le lemme en découle.
\end{demo}

Pour appliquer ces résultats dans des espaces de fonctions de carré
inté\-grable sur des graphes, nous aurons besoin de résultats de
géométrie des graphes. Soient $\Phi$ un graphe connexe et $P$ et $Q$
deux sous-ensembles disjoints de $\Phi$ tels que $\Phi=P\cup Q$.
Nous dirons que $\Phi$ est partagé par la partition $\{P,Q\}$ si
tout voisin d'un élément de $P$ appartient à $Q$ et tout voisin d'un
élément de $Q$ appartient à $P$. Nous dirons que $\Phi$ est
partageable s'il existe une partition de $\Phi$ en deux
sous-ensembles qui le partage. On vérifie aisément que $\Phi$ est
partageable si et seulement si, pour tous $p$ et $q$ dans $\Phi$,
les chemins joignant $p$ à $q$ sont soit tous de longueur paire,
soit tous de longueur impaire. En particulier, si $\Phi$ est
partageable, la partition $\{P,Q\}$ qui le partage est unique, deux
points $p$ et $q$ appartenant au même atome si et seulement s'ils
peuvent être joints par un chemin de longueur paire.

\begin{Lem}\label{graphepartageable} Soit $\Phi$ un graphe connexe et soit
$L$ l'espace des fonctions $\varphi$ sur $\Phi$ telles que, pour
tout $p$ dans $\Phi$, $\varphi$ est constante sur les voisins de
$p$. Alors, si $\Phi$ n'est pas partageable, $L$ est égal à l'espace
des fonctions constantes. Si $\Phi$ est partagé par la partition
$\{P,Q\}$, $L$ est engendré par les fonctions constantes et par la
fonction $1_P-1_Q$.\end{Lem}

\begin{demo} Soit $\varphi$ dans $L$, $p$ et $q$ des points de $\Phi$ et
$r_0=p,r_1,\ldots,r_n=q$ un chemin de $p$ à $q$. Pour tout $1\leq
i\leq n-1$, on a $r_{i-1}\sim r_i\sim r_{i+1}$, donc
$\varphi(r_{i-1})=\varphi(r_{i+1})$ et, si $n$ est pair,
$\varphi(p)=\varphi(q)$. Par conséquent, si
$\varphi(p)\neq\varphi(q)$ et si
$P=\{r\in\Phi|\varphi(r)=\varphi(p)\}$ et
$Q=\{r\in\Phi|\varphi(r)=\varphi(q)\}$, la partition $\{P,Q\}$
partage $\Phi$. Le lemme s'en déduit facilement.
\end{demo}

Nous utiliserons le principe du maximum sous la forme du

\begin{Lem}\label{maximum}
Soit $\Phi$ un graphe connexe régulier de valence $3$. Soit
$\varphi$ dans $\ell^2(\Phi)$ telle que $\Delta\varphi=3\varphi$. Si
$\Phi$ est infini, on a $\varphi=0$ et si $\Phi$ est fini, $\varphi$
est constante. Soit $\psi$ dans $\ell^2(\Phi)$ telle que
$\Delta\psi=-3\psi$. Si $\Phi$ est infini ou non partageable, on a
$\psi=0$ et, si $\Phi$ est fini et partagé par la partition
$\{P,Q\}$, $\psi$ est proportionnelle à $1_P-1_Q$.\end{Lem}

\begin{demo} Comme $\varphi$ est dans $\ell^2(\Phi)$, l'ensemble
$M=\{p\in\Phi|\varphi(p)=\max_\Phi\varphi\}$ n'est pas vide. Comme
$\Delta\varphi=3\varphi$, pour tout $p$ dans $M$, les voisins de $p$
appartiennent tous à $M$ et, donc, comme $\Phi$ est connexe,
$M=\Phi$ et $\varphi$ est constante. Si $\Phi$ est infini, comme
$\varphi$ est dans $\ell^2(\Phi)$, elle est nulle. De même,
supposons $\psi\neq 0$ et posons
$P=\{p\in\Phi|\psi(p)=\max_\Phi\psi\}$ et
$Q=\{q\in\Phi|\psi(q)=\min_\Phi\psi\}$. Comme $\Delta\psi=-3\psi$,
on a $\min_\Phi\psi=-\max_\Phi\psi$ et les voisins des points de $P$
appartiennent à $Q$ tandis que les voisins des points de $Q$
appartiennent à $P$. Par connexité, on a $P\cup Q=\Phi$, le graphe
$\Phi$ est partageable et $\psi$ est proportionnelle à $1_P-1_Q$.
Enfin, comme $\psi$ est dans $\ell^2(\Phi)$, le graphe $\Phi$ est
fini.
\end{demo}

Rappelons qu'on a noté $f$ le polynôme $x^2-x-3$. Du lemme
\ref{relation1}, on déduit le

\begin{Cor} \label{transformespectre1}
Soient $\Phi$ un graphe connexe régulier de valence $3$ et $H$ le
sous-espace fermé de $\ell^2\left(\hat{\Phi}\right)$ engendré par
l'image de $\Pi^*$ et par l'image de $\Delta\Pi^*$. Alors $H$ est
stable par $\Delta$ et le spectre de la restriction de $\Delta$ à
$H$ est,
\begin{cond}
\item si $\Phi$ est infini, l'image inverse par $f$ du spectre de $\Delta$ dans $\ell^2(\Phi)$.
\item si $\Phi$ est fini, mais non partageable, l'image inverse par $f$ du
spectre de $\Delta$ dans $\ell^2(\Phi)$ privée de $-2$.
\item si $\Phi$ est fini et partageable, l'image inverse par $f$ du
spectre de $\Delta$ dans $\ell^2(\Phi)$ privée de $-2$ et de $0$.\end{cond}
\end{Cor}

\begin{demo} Notons $K$ l'image de $\Pi^*$. Comme $\frac{1}{\sqrt{3}}\Pi^*$
induit une isométrie de $\ell^2(\Phi)$ dans $K$,
d'après le lemme \ref{relation1}, le spectre de $f(\Delta)$ dans $K$
est égal au spectre de $\Delta$ dans $\ell^2(\Phi)$. Nous allons
appliquer le lemme \ref{spectral2} à l'espace $H$ et à l'opérateur
$\Delta$. Pour cela, étudions l'espace $\Delta^{-1}K\cap K$. Soit
$L$ l'espace des $\varphi$ dans
 $\ell^2(\Phi)$ tels que $\Delta\Pi^*\Phi$ appartienne à $K$
et soit $\varphi$ dans $L$. Si $p$ est un point de $\Phi$, dont les
voisins sont $q$, $r$ et $s$, posons $\varphi(p)=a$, $\varphi(q)=b$,
$\varphi(r)=c$ et $\varphi(s)=d$. Alors, on a
$\Delta\Pi^*\varphi(p,q)=2a+b$, $\Delta\Pi^*\varphi(p,r)=2a+c$ et
$\Delta\Pi^*\varphi(p,s)=2a+d$. Comme $\Delta\Pi^*\varphi$
appartient à $K$, il vient $b=c=d$. Réciproquement, si $\varphi$ est
un élément de $\ell^2(\Phi)$ qui, pour tout point $p$ de $\Phi$, est
constant sur l'ensemble des voisins de $p$, alors $\varphi$
appartient à $L$.

Si $\Phi$ est infini, d'après le lemme \ref{graphepartageable}, on a
$L=\{0\}$, et l'on peut appliquer le lemme \ref{spectral2} à $H$. Le
spectre de $\Delta$ dans $H$ est alors bien l'image réciproque par
$f$ de celui de $\Delta$ dans $\ell^2(\Phi)$. Si $\Phi$ est fini,
mais non partageable, d'après le lemme \ref{graphepartageable}, $L$
est la droite des fonctions constantes et l'on peut appliquer le
lemme \ref{spectral2} à l'orthogonal des fonctions constantes dans
$H$. On obtient le résultat puisque $f(-2)=f(3)=3$ et que, d'après
le lemme \ref{maximum}, les fonctions constantes sont les seules
fonctions propres de valeur propre $3$ pour $\Delta$ dans
$\ell^2(\Phi)$. Enfin, si $\Phi$ est fini et partagé par la
partition $\{P,Q\}$, d'après le lemme \ref{graphepartageable}, $L$
est engendré par les fonctions constantes et par la fonction
$1_P-1_Q$. Alors, $\Pi^*(1_P-1_Q)$ est un vecteur propre de valeur
propre $1$ dans $H$. On applique le lemme \ref{spectral2} à
l'orthogonal du sous-espace de $H$ engendré par les fonctions
constantes et par $\Pi^*(1_P-1_Q)$. Le résultat en découle puisque
$f(0)=f(1)=-3$ et que, toujours d'après le lemme \ref{maximum}, les
fonctions propres de valeur propre $-3$ dans $\ell^2(\Phi)$ sont les
multiples de $1_P-1_Q$.
\end{demo}

Il nous reste à déterminer le spectre de $\Delta$ dans l'orthogonal de $H$.
C'est l'objet du

\begin{Lem} \label{spectreresiduel} Soient $\Phi$ un graphe connexe régulier
de valence $3$ et $H$ le sous-espace fermé de
$\ell^2\left(\hat{\Phi}\right)$ engendré par l'image de $\Pi^*$ et
par l'image de $\Delta\Pi^*$. Le spectre de $\Delta$ dans
l'orthogonal de $H$ est égal à $\{0,-2\}$. L'espace propre associé à
la valeur propre $0$ dans $\ell^2\left(\hat{\Phi}\right)$ est
l'espace des fonctions $\varphi$ dans
$\ell^2\left(\hat{\Phi}\right)$ telles que $\Pi\varphi=0$ et que,
pour tous $p$ et $q$ voisins dans $\Phi$, on ait
$\varphi(p,q)=\varphi(q,p)$. L'espace propre associé à la valeur
propre $-2$ dans $\ell^2\left(\hat{\Phi}\right)$ est l'espace des
fonctions $\varphi$ dans $\ell^2\left(\hat{\Phi}\right)$ telles que
$\Pi\varphi=0$ et que, pour tous $p$ et $q$ voisins dans $\Phi$, on
ait $\varphi(p,q)=-\varphi(q,p)$.
\end{Lem}

\begin{demo} Soit $\varphi$ dans l'orthogonal de $H$ et soit $p$ un point de $\Phi$, de voisins $q,r,s$.
Posons $a=\varphi(p,q)$, $b=\varphi(p,r)$, $c=\varphi(q,p)$ et
$d=\varphi(r,p)$. Enfin, notons $\psi$ la fonction indicatrice de
$\{p\}$ sur $\Phi$. Comme $\varphi$ est orthogonale à $\Pi^*\psi$ et
à $\Delta\Pi^*\psi$, on a $\varphi(p,s)=-a-b$ et
$\varphi(s,p)=-c-d$. Il vient $\Delta\varphi(p,q)=c-a$ et
$\Delta^2\varphi(p,q)=(a-c)+(d-b)+(-c-d+a+b)=2a-2c$. On a donc, dans
l'orthogonal de $H$, $\Delta^2+2\Delta=0$, et, pour $\varphi$ dans
cet espace, on a $\Delta\varphi=0$ si et seulement si, pour tous $p$
et $q$ voisins dans $\Phi$,  $\varphi(p,q)=\varphi(q,p)$ et
$\Delta\varphi=-2\varphi$ si et seulement si, pour tous $p$ et $q$
voisins dans $\Phi$,  $\varphi(p,q)=-\varphi(q,p)$.

Pour terminer la démonstration du lemme, il nous reste à montrer que
$\Pi$ n'a pas de vecteur propre de valeur propre $0$ ou $-2$ dans
$H$. Pour cela, donnons-nous donc $\varphi$ dans $H$ tel que
$\Delta\varphi=-2\varphi$. On a alors, d'après le lemme
\ref{relation1},
$\Delta\Pi\varphi=\Pi(\Delta^2-\Delta-3)\varphi=3\varphi$. Si $\Phi$
est infini, d'après le lemme \ref{maximum}, $\Pi\varphi$ est nulle.
On a donc $\Pi\varphi=0$ et $\Pi\Delta\varphi=-2\Pi\varphi=0$. Donc,
comme $\varphi$ appartient à $H$ qui est engendré par l'image de
$\Pi^*$ et par l'image de $\Delta\Pi^*$, on a $\varphi=0$. Si $\Phi$
est fini, toujours d'après le lemme \ref{maximum}, $\Pi\varphi$ est
constante. Comme $\varphi$ est orthogonale aux fonctions constantes,
on a donc, à nouveau, $\Pi\varphi=0$ et $\Pi\Delta\varphi=0$, ce qui
implique $\varphi=0$.

Si, à présent, $\varphi$ est un élément de $H$ tel que
$\Delta\varphi=0$, on a $\Delta\Pi\varphi=-3\Pi\varphi$. \`A
nouveau, d'après le lemme \ref{maximum}, si $\Phi$ est infini ou non
partageable, on a $\Pi\varphi=0$ et, donc, $\varphi=0$, tandis que,
si $\Phi$ est partagé par la partition $\{P,Q\}$, $\Pi\varphi$ est
proportionnelle à $1_P-1_Q$. Or, $\Pi^*(1_P-1_Q)$  est un vecteur
propre de valeur propre $1$ pour $\Delta$, donc
$\langle\Pi\varphi,1_P-1_Q\rangle=\langle\varphi,\Pi^*(1_P-1_Q)\rangle=0$
et $\Pi\varphi=0$, si bien que $\varphi=0$.
\end{demo}

\begin{figure}\begin{center}\input{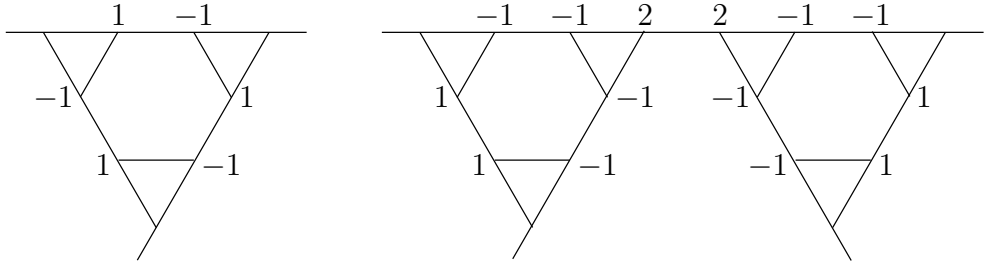}\caption{Fonctions propres sur $\hat{\Phi}^{(2)}$}
\label{fonctionspropres}\end{center}\end{figure}

Rappelons que, pour tout entier $n$, on a noté $\hat{\Phi}^{(n)}$ le
graphe obtenu en rempla\c cant chaque point de $\Phi$ par un
$n$-triangle. L'espace $\ell^2\left(\hat{\Phi}^{(2)}\right)$
contient des fonctions propres à support fini de valeur propre $-2$
et $0$, comme on peut le voir sur la figure \ref{fonctionspropres},
où l'on n'a représenté que les valeurs non nulles des fonctions.

Nous avons donc le

\begin{Lem} \label{existevecteurpropre}
Pour tout $n\geq 2$, l'espace $\ell^2\left(\hat{\Phi}^{(n)}\right)$
contient des fonctions propres de valeur propre $-2$ et
$0$.\end{Lem}

Rappelons qu'on a noté $\Lambda$ l'ensemble de Julia de $f$. En
appliquant le corollaire \ref{transformespectre1} et le lemme
\ref{spectreresiduel} à $\Gamma$, on obtient le

\begin{Cor} \label{spectre}
Le spectre de $\Gamma$ est constitué de la réunion de $\Lambda$
et de l'ensemble $\bigcup_{n\in\mathbb N}f^{-n}(0)$.\end{Cor}

\begin{demo} D'après la proposition \ref{structurePascal}, $\hat{\Gamma}$ est isomorphe à $\Gamma$.
Par consé\-quent, d'après le lemme \ref{existevecteurpropre} et le
corollaire \ref{transformespectre1}, le spectre de $\Gamma$ contient
$\bigcup_{n\in\mathbb N}f^{-n}(0)$ et, donc, l'adhérence de cet
ensemble, qui est précisément la réunion de $\Lambda$ et de
$\bigcup_{n\in\mathbb N}f^{-n}(0)$. Or, comme $-2$ appartient à
$\Lambda$, si $x$ est un point du spectre de $\Gamma$ qui
n'appartient pas à
 $\bigcup_{n\in\mathbb N}f^{-n}(0)$, d'après le corollaire
\ref{transformespectre1}, pour tout $n$ dans $\mathbb N$, le réel
$f^n(x)$ appartient au spectre de $\Gamma$ et, donc, la suite
$(f^n(x))_{n\in\mathbb N}$ est bornée, si bien que $x$ appartient à
$\Lambda$.
\end{demo}

\section{Les mesures spectrales de $\Gamma$}
\label{mesurespectrale}

Soit toujours $\Phi$ un graphe connexe régulier de valence $3$. Dans
cette section, nous allons expliquer comment calculer les mesures
spectrales de certains éléments de $\ell^2\left(\hat{\Phi}\right)$.
Pour cela nous utiliserons le

\begin{Lem} \label{relation2} On a $\Pi\Delta\Pi^*=6+\Delta$ et,
donc, pour tous $\varphi$ et $\psi$ dans $\ell^2(\Phi)$,
$$\langle\Delta\Pi^*\varphi,\Pi^*\psi\rangle=6\langle\varphi,\psi\rangle+\langle\Delta\varphi,\psi\rangle
=2\langle\Pi^*\varphi,\Pi^*\psi\rangle+\frac{1}{3}\langle(\Delta^2-\Delta-3)\Pi^*\varphi,\Pi^*\psi\rangle.$$
\end{Lem}

\begin{demo} Soient $\varphi$ dans $\ell^2(\Phi)$ et $p$ un point de
$\Phi$, de voisins $q,r,s$. Posons $a=\varphi(p)$, $b=\varphi(q)$,
$c=\varphi(r)$ et $d=\varphi(s)$. Alors, on a
$\Delta\Pi^*\varphi(p,q)=2a+b$ et, donc,
$\Pi\Delta\Pi^*\varphi(p)=(2a+b)+(2a+c)+(2a+d)=6a+(b+c+d)$, d'où la
première identité. La deuxième en découle, en appliquant le lemme
\ref{relation1} et la relation $\Pi\Pi^*=3$.
\end{demo}

\'Etudions à présent les conséquences abstraites de ce type
d'identité.

\begin{Lem} \label{mesurespectrale1}
Soient $H$ un espace de Hilbert, $T$ un endomorphisme auto-adjoint
borné de $H$, $K$ un sous-espace fermé de $H$ et $\pi(x)=(x-u)^2+m$
un polynôme unitaire à coefficients réels du second degré. On
suppose qu'on a $\pi(T)K\subset K$, que $K$ et $TK$ engendrent $H$
et qu'il existe des nombres réels $a$ et $b$ tels que, pour tous $v$
et $w$ dans $K$, on ait $\langle Tv,w\rangle=a\langle
v,w\rangle+b\langle \pi(T)v,w\rangle$. Alors, pour tout $x\neq u$
dans le spectre de $T$, on a
$$1+\frac{a-u+b\pi(x)}{x-u}\geq 0.$$
\end{Lem}

\begin{demo} Après avoir mis $\pi$ sous forme canonique, on peut supposer
qu'on a $\pi(x)=x^2$. Soit $v$ un vecteur unitaire de $K$. Alors,
pour tout nombre réel $s$, on a
\begin{align*}0\leq \langle Tv+sv,Tv+sv\rangle
&=\langle Tv,Tv\rangle+2s\langle Tv,v\rangle+s^2\\
&=\langle T^2v,v\rangle+2s(a+b\langle
T^2v,v\rangle)+s^2.\end{align*} D'après le lemme \ref{spectral2},
les carrés des éléments du spectre de $T$ dans $H$ appartiennent au
spectre de $T^2$ dans $K$. Si $x$ est un élément du spectre de $T$,
il existe donc des vecteurs unitaires $v$ de $K$ tels que $\langle
T^2v,v\rangle$ soit aussi proche que l'on veut de $x^2$. Alors, par
la remarque ci-dessus, pour tout nombre réel $s$, on a
$x^2+2s(a+bx^2)+s^2\geq 0$. Le discriminant de ce polynôme du second
degré est donc négatif, c'est-à-dire qu'on a $x^2-(a+bx^2)^2\geq 0$.
Le lemme en découle.
\end{demo}

Soient  $\pi(x)=(x-u)^2+m$ un polynôme unitaire du second degré à
coefficients réels. Donnons-nous une fonction borélienne $\theta$
sur $\mathbb R-\{u\}$. Alors, si $\alpha$ est une fonction
borélienne sur $\mathbb R-\{u\}$, on note, pour tout $y$ dans
$]m,\infty[$,
$L_{\pi,\theta}\alpha(y)=\sum_{\pi(x)=y}\theta(x)\alpha(x)$. Soit
$\mu$ une mesure borélienne positive sur $]m,\infty[$. Si, pour
$\mu$-presque tout $y$ dans $]m,\infty[$, $\theta$ est positive sur
les deux images inverses de $y$ par $\pi$, on note
$L_{\pi,\theta}^*\mu$ la mesure borélienne $\nu$ sur $\mathbb
R-\{u\}$ telle que, pour toute fonction borélienne positive $\alpha$
sur $\mathbb R-\{u\}$, on ait $\int_{\mathbb
R-\{u\}}\alpha\de\nu=\int_{]m,\infty[}L_{\pi,\theta}\alpha\de\mu$.

Nous avons le

\begin{Lem} \label{mesurespectrale2}
Soient $H$ un espace de Hilbert, $T$ un endomorphisme auto-adjoint
borné de $H$, $K$ un sous-espace fermé de $H$ et $\pi(x)=(x-u)^2+m$
un polynôme unitaire à coefficients réels du second degré. On
suppose qu'on a $\pi(T)K\subset K$, que $K$ et $TK$ engendrent $H$
et qu'il existe des nombres réels $a$ et $b$ tels que, pour tous $v$
et $w$ dans $K$, on ait $\langle Tv,w\rangle=a\langle
v,w\rangle+b\langle \pi(T)v,w\rangle$. Alors, pour tout $v$ dans
$K$, si $\mu$ est la mesure spectrale de $v$ pour $\pi(T)$ et $\nu$
sa mesure spectrale pour $T$, si $\mu(m)=0$, on a $\nu(u)=0$ et
$\nu=L_{\pi,\theta}^*\mu$ où, pour tout $x\neq u$, on a
$$\theta(x)=\frac{1}{2}\left(1+\frac{a-u+b\pi(x)}{x-u}\right).$$\end{Lem}

\begin{demo} Notons que, comme $\mu(m)=0$, d'après le lemme
\ref{spectral2}, la mesure $\mu$ est concentrée sur $]m,\infty[$. De
plus, si $w$ est un vecteur de $H$ tel que $Tw=uw$, on a
$\pi(T)w=mw$ et, par hypothèse, $\langle v,w\rangle=0$. On a donc
$\nu(u)=0$.

D'après le lemme \ref{mesurespectrale1}, la fonction $\theta$ est
positive sur le spectre de $T$ privé de $u$. Soit $n$ dans $\mathbb
N$. D'une part, on a $$\int_{\mathbb
R-\{u\}}\pi(x)^n\de\nu(x)=\langle\pi(T)^nv,v\rangle=
\int_{]m,\infty[}y^n\de\mu(y).$$ D'autre part, pour tout $x\neq u$,
on a $\theta(x)+\theta(2u-x)=1$ et, donc, pour tout $y$ dans
$]m,\infty[$, $L_{\pi,\theta}\pi^n(y)=y^n$. Il vient bien
$\int_{\mathbb
R-\{u\}}\pi^n\de\nu=\int_{]m,\infty[}L_{\pi,\theta}\pi^n\de\mu$. De
même, pour tout tout $x$ dans $\mathbb R$, posons
$\alpha(x)=x\pi(x)^n$. D'une part, on a alors
\begin{multline*}\int_{\mathbb R-\{u\}}\alpha(x)\de\nu(x)
=\langle T\pi(T)^nv,v\rangle\\=a\langle \pi(T)^nv,v\rangle+b\langle
\pi(T)^{n+1}v,w\rangle
=a\int_{]m,\infty[}y^n\de\mu(y)+b\int_{]m,\infty[}y^{n+1}\de\mu(y).\end{multline*}
D'autre part, pour tout $x\neq u$, comme $(2u-x)-u=u-x$, on a
\begin{align*}&\theta(x)\alpha(x)+\theta(2u-x)\alpha(2u-x)\\
&\quad=\left(\frac{1}{2}\left(x+\left(2u-x\right)\right)+\frac{1}{2}\left(x-\left(2u-x\right)\right)\left(\frac{a-u+b\pi(x)}{x-u}\right)\right)\pi(x)^n\\
&\quad=a\pi(x)^n+b\pi(x)^{n+1}\end{align*} et, donc, pour tout $y$
dans $]m,\infty[$, $L_{\pi,\theta}\alpha(y)=ay^n+by^{n+1}$. Il vient
à nouveau $\int_{\mathbb
R-\{u\}}\alpha\de\nu=\int_{]m,\infty[}L_{\pi,\theta}\alpha\de\mu$.
Par conséquent, pour tout polynôme $\alpha$, on a $\int_{\mathbb
R-\{u\}}\alpha\de\nu=\int_{]m,\infty[}L_{\pi,\theta}\alpha\de\mu$.
En particulier, la mesure positive $L_{\pi,\theta}^*\mu$ est finie
et, donc, pour toute fonction continue $\alpha$ à support compact
dans $\mathbb R$, on a encore $\int_{\mathbb
R-\{u\}}\alpha\de\nu=\int_{]m,\infty[}L_{\pi,\theta}\alpha\de\mu$,
si bien que $\nu=L_{\pi,\theta}^*\mu$.
\end{demo}

Des lemmes \ref{relation2} et \ref{mesurespectrale2}, on déduit le

\begin{Cor}\label{transformemesure}
Soient $\varphi$ dans $\ell^2(\Phi)$, $\mu$ la mesure spectrale de
$\varphi$ pour $\Delta$ dans $\ell^2(\Phi)$ et $\nu$ la mesure
spectrale de $\Pi^*\varphi$ pour $\Delta$ dans
$\ell^2\left(\hat{\Phi}\right)$. Alors, on a
$\nu\left(\frac{1}{2}\right)=0$ et, si, pour tout
$x\neq\frac{1}{2}$, on pose $\theta(x)=\frac{x(x+2)}{2x-1}$, on a
$\nu=L_{f,\theta}^*\mu$.
\end{Cor}

\begin{demo} La valeur minimale du polynôme $f$ sur $\mathbb R$ est
$f\left(\frac{1}{2}\right)=-\frac{13}{4}<-3\leq-\N{\Delta}_2$. On a
donc $\mu\left(-\frac{13}{4}\right)=0$ et le corollaire découle des
lemmes \ref{relation2} et \ref{mesurespectrale2} par un calcul
élémentaire.\end{demo}

\section{Fonctions propres dans $\ell^2(\Gamma)$}
\label{valeurp}

Dans cette section, nous allons compléter les informations données
par le lemme \ref{spectreresiduel} en décrivant plus précisément les
espaces propres de $\Delta$ dans $\ell^2(\Gamma)$ pour les valeurs
propres $-2$ et $0$. Nous étendrons ces résultats aux valeurs
propres dans $\bigcup_{n\in\mathbb N}f^{-n}(-2)$ et dans
$\bigcup_{n\in\mathbb N}f^{-n}(0)$ grâce au

\begin{Lem} \label{transformevalp}
Soient $\Phi$ un graphe régulier de valence $3$ et $H$ le
sous-espace fermé de $\ell^2\left(\hat{\Phi}\right)$ engendré par
l'image de $\Pi^*$ et par celle de $\Delta\Pi^*$. Alors, pour tout
$x$ dans $\mathbb R-\{0,-2\}$, $x$ est valeur propre de $\Delta$
dans $H$ si et seulement si $y=f(x)$ est valeur propre de $\Delta$
dans $\ell^2(\Phi)$. Dans ce cas, l'application $R_x$ qui, à une
fonction propre $\varphi$ de valeur propre $y$ dans $\ell^2(\Phi)$,
associe $(x-1)\Pi^*\varphi+\Delta\Pi^*\varphi$ induit un
isomorphisme entre l'espace propre de valeur propre $y$ dans
$\ell^2(\Phi)$ et l'espace propre de valeur propre $x$ dans $H$ et,
pour tout $\varphi$, on a
$\N{R_x\varphi}_2^2=x(x+2)(2x-1)\N{\varphi}_2^2$.
\end{Lem}

\begin{demo} Soit $\psi\neq 0$ dans $H$ tel que $\Delta\psi=x\psi$.
Puisque $\psi$ est dans $H$, on a $\Pi\psi\neq 0$ ou
$\Pi\Delta\psi\neq 0$. Comme $\Pi\Delta\psi=x\Pi\psi$, on a
$\Pi\psi\neq 0$. D'après le lemme \ref{relation1}, on a
$\Delta\Pi\psi=y\Pi\psi$, donc $y$ est valeur propre de $\Delta$
dans $\ell^2(\Phi)$. En particulier, comme
$f\left(\frac{1}{2}\right)=-\frac{13}{4}<-3\leq-\N{\Delta}_2$, on a
$x\neq\frac{1}{2}$. Réciproquement, si $\varphi$ est un élément de
$\ell^2(\Phi)$ tel que $\Delta\varphi=y\varphi$, on a, d'après le
lemme \ref{relation1},
\begin{align*}\Delta R_x\varphi&=\Delta((x-1)\Pi^*\varphi+\Delta\Pi^*\varphi)\\
&=(x-1)\Delta\Pi^*\varphi+\Pi^*\Delta\varphi+(\Delta+3)\Pi^*\varphi\\
&=x\Delta\Pi^*\varphi+(x^2-x)\Pi^*\varphi=xR_x\varphi.\end{align*}
Or, d'après le lemme \ref{relation2}, on a
$\Pi\Delta\Pi^*=6+\Delta$, donc, si $\Delta\varphi=y\varphi$,  on a,
par un calcul immédiat, $\Pi R_x\varphi=x(x+2)\varphi$ et, comme on
a supposé $x(x+2)\neq 0$, $R_x$ est injectif et fermé. Il nous reste
à montrer que $R_x$ est surjectif. Pour cela, considérons $\psi$
dans $H$ tel que $\Delta\psi=x\psi$ mais que $\psi$ soit orthogonal
à l'image de $R_x$. Pour tout $\varphi$ dans $\ell^2(\Phi)$ tel que
$\Delta\varphi=y\varphi$, on a $\langle\psi,R_x\varphi\rangle=
(x-1)\langle\Pi\psi,\varphi\rangle+\langle\Pi\Delta\psi,\varphi\rangle
=(2x-1)\langle\Pi\psi,\varphi\rangle$ et donc, comme
$x\neq\frac{1}{2}$, $\langle\Pi\psi,\varphi\rangle=0$. Comme
$\Delta\Pi\psi=y\Pi\psi$, on a $\Pi\psi=0$. Comme $\psi$ est dans
$H$, on a $\psi=0$. L'opérateur $R_x$ est donc un isomorphisme. Le
calcul de norme est alors direct, en utilisant les lemmes
\ref{relation1} et \ref{relation2}.
\end{demo}

Commen\c cons par nous intéresser aux valeurs propres dans
$\bigcup_{n\in\mathbb N}f^{-n}(0)$. Soit $n$ un entier $\geq 1$.
Rappelons que, d'après le corollaire \ref{sommettriangles}, si
$\mathcal T$ est un $n$-triangle de $\Gamma$ et si $p$ est un sommet
de $\mathcal T$, le voisin de $p$ qui n'appartient pas à $\mathcal
T$ est le sommet d'un $n$-triangle. On appellera arêtes extérieures
aux $n$-triangles les arêtes joignant deux points qui sont des
sommets d'un $n$-triangle. On notera $\Theta_n$ l'ensemble des
arêtes extérieures aux $n$-triangles et on le munira de la structure
de graphe pour laquelle deux arêtes différentes sont voisines si
deux de leurs extrémités sont sommets d'un même $n$-triangle. On
vérifie aisément qu'alors le graphe $\Theta_n$ est naturellement
isomorphe au graphe de Sierpi\'nski $\Theta$ introduit à la fin de
la section \ref{decritgraphes}. On identifiera désormais $\Theta$ et
$\Theta_n$. Si $\varphi$ est une fonction sur $\Gamma$ qui est
constante sur les arêtes extérieures aux $n$-triangles, on notera
$P_n\varphi$ la fonction sur $\Theta$ dont la valeur en un point de
$\Theta$ est la valeur de $\varphi$ sur l'arête extérieure aux
$n$-triangles associée. Enfin, rappelons que, comme $\Theta$ est
régulier de valence $4$, la norme de $\Delta$ dans $\ell^2(\Theta)$
est $\leq 4$.

D'après le lemme \ref{spectreresiduel}, les fonctions propres de
valeur propres $0$ sont constantes sur les arêtes extérieures aux
$1$-triangles. On a le

\begin{Lem} \label{espacep01}
L'application $P_2$ induit un isomorphisme d'espaces de Banach de
l'espace propre de $\ell^2(\Gamma)$ associé à la valeur propre $0$
dans $\ell^2(\Theta)$. On note $Q_{0}$ sa réciproque. Pour tout
$\psi$ dans $\ell^2(\Theta)$, on a $\N{Q_
{0}\psi}^2_{\ell^2(\Gamma)}=3\N{\psi}^2_{\ell^2(\Theta)}-\frac{1}{2}\langle\Delta\psi,\psi\rangle_{\ell^2(\Theta)}$.
\end{Lem}

\begin{demo} En utilisant la caractérisation du lemme \ref{spectreresiduel},
on vérifie ai\-sé\-ment qu'étant données trois valeurs $a$, $b$ et
$c$ aux sommets d'un $2$-triangle, une fonction propre de valeur
propre $0$ prenant ces valeurs aux trois sommets doit avoir à
l'intérieur du triangle les valeurs représentées sur la figure
\ref{valeurp0}.

\begin{figure}\begin{center}\input{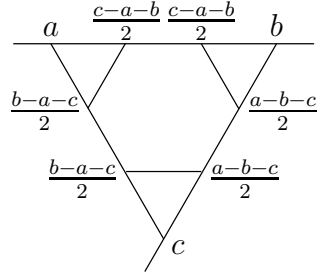}\caption{Les fonctions propres de valeur propre $0$}
\label{valeurp0}\end{center}\end{figure}

Rappelons qu'on a identifié les graphes $\Theta$ et $\Theta_2$. Pour
toute fonction $\psi$ sur $\Theta$, notons $Q_{0}\psi$ la fonction
sur $\Gamma$ qui, sur une arête extérieure à un $2$-triangle de
$\Gamma$, est constante, de valeur la valeur de $\psi$ en le point
de $\Theta$ associé à cette arête, et dont les valeurs à l'intérieur
des $2$-triangles sont celles décrites par la figure \ref{valeurp0}.
Par un calcul élémentaire, pour tous réels $a$, $b$ et $c$, la somme
des carrés des valeurs représentés sur la figure \ref{valeurp0} est
$$\frac{3}{2}(a^2+b^2+c^2)-ab-ac-bc=\frac{1}{2}(3a^2-ab-ac)+\frac{1}{2}(3b^2-ab-bc)+\frac{1}{2}(3c^2-ac-bc),$$
si bien que, pour toute fonction $\psi$ sur $\Theta$, on a $\N{Q_
{0}\psi}^2_{\ell^2(\Gamma)}=3\N{\psi}^2_{\ell^2(\Theta)}-\frac{1}{2}\langle\Delta\psi,\psi\rangle_{\ell^2(\Theta)}$.
Comme
$-4\N{\psi}^2_{\ell^2(\Theta)}\leq\langle\Delta\psi,\psi\rangle_{\ell^2(\Theta)}\leq
4\N{\psi}^2_{\ell^2(\Theta)}$, la fonction $\psi$ appartient à
$\ell^2(\Theta)$ si et seulement si $Q_ {0}\psi$ appartient à
$\ell^2(\Gamma)$. Le lemme en découle.
\end{demo}

Des lemmes \ref{transformevalp} et \ref{espacep01}, nous allons déduire une description
des espaces propres associés aux éléments de $\bigcup_{n\in\mathbb N}f^{-n}(0)$.
Pour $x$ dans $\bigcup_{n\in\mathbb N}f^{-n}(0)$, notons $n(x)$ l'entier $n$ tel que $f^n(x)=0$ et
$$\kappa(x)=\prod_{k=0}^{n(x)-1}\frac{f^k(x)(2f^k(x)-1)}{f^k(x)+2}.$$
On a la

\begin{Prop} \label{espacep02}
Soit $x$ dans $\bigcup_{n\in\mathbb N}f^{-n}(0)$. Les fonctions
propres de valeur propre $x$ dans $\ell^2(\Gamma)$ sont constantes
sur les arêtes extérieures aux $(n(x)+1)$-triangles dans $\Gamma$.
L'application $P_{n(x)+2}$ induit un isomorphisme d'espaces de
Banach de l'espace propre de $\ell^2(\Gamma)$ associé à la valeur
propre $x$ dans $\ell^2(\Theta)$. On note $Q_{x}$ sa réciproque.
Alors, pour tout $\psi$ dans $\ell^2(\Theta)$, on a
$\N{Q_{x}\psi}^2_{\ell^2(\Gamma)}=
\kappa(x)\left(3\N{\psi}^2_{\ell^2(\Theta)}-\frac{1}{2}\langle\Delta\psi,\psi\rangle_{\ell^2(\Theta)}\right)$.
\end{Prop}

\begin{demo} On démontre ce résultat par récurrence sur $n(x)$. Le cas $n(x)=0$ a fait
l'object du lemme \ref{espacep01}. Supposons le lemme démontré pour $n(y)$ avec $y=f(x)$.
Et choisissons $\varphi$ dans $\ell^2(\Gamma)$ tel que
$\Delta\varphi=x\varphi$. Alors, comme $n=n(x)\geq 1$, on a
$x\notin\{-2,0\}$ et, donc, d'après le lemme \ref{spectreresiduel},
$\varphi$ appartient à $H$. D'après le lemme \ref{transformevalp},
on a donc $\varphi=R_x\psi$, pour une fonction $\psi$ telle que
$\Delta\psi=y\psi$. Par récurrence, $\psi$ est constante sur les
arêtes extérieures aux $n$-triangles. Soient donc $p$ un sommet d'un
$(n+1)$-triangle dans $\Gamma$ et $q$ son voisin extérieur. Les
points $\Pi p$ et $\Pi q$ sont des sommets voisins de $n$-triangles
de $\Gamma$. On a donc
 $\Pi^*\psi(p)=\psi(\Pi p)=\psi(\Pi q)=\Pi^*\psi(q)$ et
$\Delta\Pi^*\psi(p)=2\psi(\Pi p))+\psi(\Pi q))=3\Pi^*\psi(p)$. Il
vient $\varphi(p)=R_x\psi(p)=(x+2)\psi(p)=\varphi(q)$ : la fonction
$\varphi$ est constante sur les arêtes extérieures aux
$(n+1)$-triangles et $P_{n+2}\varphi=(x+2)P_{n+1}\psi$. Comme, par
récurrence, $P_{n+1}$ induit un isomorphisme de l'espace propre de
valeur propre $y$ sur $\ell^2(\Theta)$, d'après le lemme
\ref{transformevalp}, $P_{n+2}$ induit un isomorphisme de l'espace
propre de valeur propre $x$ sur $\ell^2(\Theta)$. Le calcul de norme
est alors une conséquence de la récurrence et de la formule
$P_{n+2}R_x=(x+2)P_{n+1}$.
\end{demo}

\begin{Cor} \label{espacep03}
Pour tout $x$ dans $\bigcup_{n\in\mathbb N}f^{-n}(0)$, l'espace
propre associé à $x$ dans $\ell^2(\Gamma)$ est de dimension infinie
et engendré par des fonctions à support fini.
\end{Cor}

Pour les éléments de $\bigcup_{n\in\mathbb N}f^{-n}(-2)$, on ne
dispose pas d'un analogue à la proposition \ref{espacep02}.
Cependant, nous allons étendre le corollaire \ref{espacep03}.
Commen\c cons par traiter le cas de la valeur propre $-2$. Rappelons
qu'on a noté $p_0$ et $p_0^\vee$ les sommets des deux triangles
infinis de $\Gamma$.

\begin{Lem} \label{espacep-21}
Soit $\varphi$ une fonction propre de valeur propre $-2$ dans
$\ell^2(\Gamma)$. Alors, pour tout $n\geq 1$ la somme des valeurs de
$\varphi$ sur les sommets de chaque $n$-triangle de $\Gamma$ est
nulle et $\varphi(p_0)=\varphi(p_0^\vee)=0$. L'espace propre associé
à la valeur propre $-2$ est de dimension infinie et engendré par des
fonctions à support fini.\end{Lem}

\begin{demo} Un calcul immédiat utilisant le lemme \ref{spectreresiduel}
montre que les valeurs de $\varphi$ dans un $2$-triangle vérifient les règles décrites par la figure \ref{valeurp-2}.
\begin{figure}\begin{center}\input{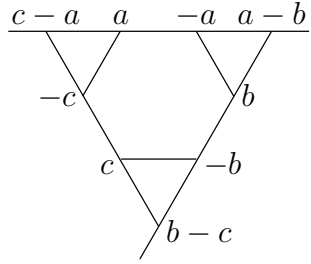}\caption{Les fonctions propres de valeur propre $-2$}
\label{valeurp-2}\end{center}\end{figure} En particulier, la somme
de leurs valeurs aux sommets de chaque $2$-triangle est nulle. Par
récurrence, en employant le corollaire \ref{decoupetriangles}, on en
déduit que, pour tout $n\geq 1$, la somme de leurs valeurs aux
sommets de chaque $n$-triangle est nulle. Notons alors, pour tout
$n\geq 1$ , $p_n$ et $q_n$ les deux autres sommets du $n$-triangle
issu de $p_0$, de fa\c con à ce qu'on ait, avec les notations du
corollaire \ref{decoupetriangles}, $p_0p_{n+1}=p_n$ et
$p_0q_{n+1}=q_n$. On a, pour tout $n\geq 1$,
$\varphi(p_0)=-\varphi(p_n)-\varphi(q_n)$. Comme $\varphi$ est de
carré intégrable, on a $\varphi(p_n)\td{n}{\infty}0$ et
$\varphi(q_n)\td{n}{\infty}0$. Il vient $\varphi(p_0)=0$ et, de
même, $\varphi(p_0^\vee)=0$. En particulier, pour tout $n\geq 1$,
$\varphi(p_n)+\varphi(q_n)=0$.

Montrons que $\varphi$ est limite d'un suite de fonctions à support
fini. Notons toujours, comme dans la section \ref{decritgraphes},
$\mathcal T_\infty(p_0)$ le triangle infini de $\Gamma$ issu de
$p_0$. Alors, comme $\varphi(p_0)=0$, $\varphi 1_{\mathcal
T_\infty(p_0)}$ est encore un vecteur propre de valeur propre $-2$
et on peut supposer qu'on a $\varphi=0$ sur $\mathcal
T_\infty(p_0^\vee)$. Pour tout entier $n$, notons $\varphi_n$ la
fonction sur $\Gamma$ qui est nulle en dehors du $(n+1)$-triangle
$\mathcal T_{n+1}(p_0,p_{n+1},q_{n+1})$, qui est égale à $\varphi_n$
sur le $n$-triangle $\mathcal T_{n}(p_0,p_{n},q_{n})$ et qui est
invariante par l'action des éléments de signature $1$ du groupe
$\mathfrak S(p_0,p_{n+1},q_{n+1})$ sur le $(n+1)$-triangle $\mathcal
T_{n+1}(p_0,p_{n+1},q_{n+1})$. Vu le corollaire
\ref{decoupetriangles}, les valeurs de $\varphi_n$ aux sommets des
$n$-triangles sont celles représentées par la figure
\ref{fonctionapprochee}.
\begin{figure}\begin{center}\input{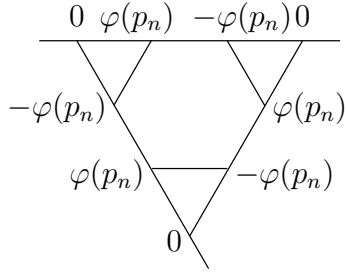}\caption{La fonction $\varphi_n$}
\label{fonctionapprochee}\end{center}\end{figure} Alors, d'après le
lemme \ref{spectreresiduel}, pour tout $n\geq 1$, $\varphi_n$ est un
vecteur propre de valeur propre $-2$ et on a $\N{\varphi_n}_2^2\leq
3\N{\varphi}_2^2$. La suite $(\varphi_n)$ tend faiblement vers
$\varphi$ dans $\ell^2(\Gamma)$. La fonction $\varphi$ appartient à
l'adhérence faible du sous-espace engendré par les fonctions propres
de valeur propre $-2$ à support fini, et, donc, à son adhérence
forte.

Enfin, l'espace des fonctions propres de valeur propre $-2$ est de
dimension infinie puisque, d'après la figure \ref{fonctionspropres},
tout $2$-triangle contient le support d'une fonction propre de
valeur propre $-2$.
\end{demo}

Pour $x$ dans  $\bigcup_{n\in\mathbb N}f^{-n}(-2)$, notons $n(x)$
l'entier $n$ tel que $f^{n(x)}(x)=-2$. Un raisonnement par
récurrence fondé sur le lemme \ref{transformevalp} permet de déduire
du lemme \ref{espacep-21} le

\begin{Cor} \label{espacep-22}
Soient $x$ dans $\bigcup_{n\in\mathbb N}f^{-n}(-2)$ et $\varphi$ une fonction propre
de valeur propre $x$ dans $\ell^2(\Gamma)$.
Alors, les valeurs de $\varphi$ sur les arêtes extérieures aux
$(n(x)+1)$-triangles sont opposées, pour tout $n\geq n(x)+1$ la
somme des valeurs de $\varphi$ sur les sommets de chaque
$n$-triangle est nulle et $\varphi(p_0)=\varphi(p_0^\vee)=0$.
L'espace propre associé à la valeur propre $x$ est de dimension
infinie et engendré par des fonctions à support fini.\end{Cor}

\section{Décomposition spectrale de $\ell^2(\Gamma)$}
\label{decompo}

Notons $\varphi_0$ la fonction sur $\Gamma$ qui vaut $1$ en $p_0$,
$-1$ en $p_0^\vee$ et $0$ partout ailleurs. Dans ce paragraphe, nous
allons montrer que $\ell^2(\Gamma)$ est la somme directe orthogonale
des espaces propres associés aux éléments de $\bigcup_{n\in\mathbb
N}f^{-n}(-2)\cup\bigcup_{n\in\mathbb N}f^{-n}(0)$ et du sous-espace
cyclique engendré par $\varphi_0$. Commen\c cons par décrire ce
dernier. Un calcul immédiat démontre le

\begin{Lem}\label{relation3} On a $\Pi^*\varphi_0=(\Delta+2)\varphi_0$.\end{Lem}

Cette relation et le corollaire \ref{transformemesure} vont nous
permettre de déterminer la mesure spectrale de $\varphi_0$.
Rappelons pour cela les propriétés des opérateurs de transfert que
nous serons amenés à utiliser : elles découlent de la version du
théorème de Ruelle-Perron-Frobenius donnée dans \cite[§ 2.2]{PP}. Si
$\kappa$ est une fonction borélienne sur $\Lambda$, on notera
$L_\kappa$ pour $L_{f,\kappa}$.

\begin{Lem}\label{transfert} Soit $\kappa:\Lambda\rightarrow\mathbb R_+^*$ une fonction höldérienne.
Munissons l'espace $\mathcal C^0(\Lambda)$ de la topologie de la
convergence uniforme. Alors, si $\lambda_\kappa>0$ est le rayon
spectral de l'opérateur $L_\kappa$ dans $\mathcal C^0(\Lambda)$, il
existe une unique probabilité borélienne $\nu_\kappa$ sur $\Lambda$
et une unique fonction continue strictement positive $l_\kappa$ sur
$\Lambda$ telles qu'on ait $L_\kappa l_\kappa=\lambda_\kappa
l_\kappa$, $L_\kappa^* \nu_\kappa=\lambda_\kappa\nu_\kappa$ et
$\int_\Lambda l_\kappa\de\nu_\kappa=1$. Le rayon spectral de
$L_\kappa$ dans l'espace des fonctions d'intégrale nulle pour
$\nu_\kappa$ est $<\lambda_\kappa$ et, en particulier, pour tout $g$
dans $\mathcal C^0(\Lambda)$, la suite
$\left(\frac{1}{\lambda_\kappa^n}L_\kappa^n(g)\right)_{n\in\mathbb
N}$ converge uniformément vers $\int_\Lambda g\de\nu_\kappa$. La
mesure $\nu_\kappa$ est diffuse est son support est $\Lambda$.
\end{Lem}

Posons, pour tout $x$ dans $\mathbb R$, $h(x)=3-x$, $k(x)=x+2$ et, pour $x\neq\frac{1}{2}$, $\rho(x)=\frac{x}{2x-1}$.
On a $h\circ f=hk$.
Du lemme \ref{relation3}, nous déduisons, grâce au corollaire \ref{transformemesure}, le

\begin{Cor} \label{spectrediffus1}
Soit $\nu_\rho$ l'unique probabilité borélienne sur $\Lambda$ telle que $L_\rho^*\nu_\rho=\nu_\rho$.
La mesure spectrale de $\varphi_0$ est $h\nu_\rho$.
\end{Cor}

\begin{demo} Soit $\mu$ la mesure spectrale de $\varphi_0$. Pour $x\neq\frac{1}{2}$,
posons $\theta(x)=\frac{x(x+2)}{2x-1}=k(x)\rho(x)$. La mesure
spectrale de $(\Delta+2)\varphi_0$ est $k^2\mu$. Par conséquent,
d'après le corollaire \ref{transformemesure} et le lemme
\ref{relation3}, on a $\mu(\frac{1}{2})=0$ et
$k^2\mu=L^*_\theta\mu$. Or, d'après le lemme \ref{espacep-21}, si
$\varphi$ est une fonction propre de valeur propre $-2$ dans
$\ell^2(\Gamma)$, on a $\varphi(p_0)=\varphi(p_0^\vee)=0$ et, donc,
$\langle\varphi,\varphi_0\rangle=0$, si bien que $\mu(-2)=0$. Par
conséquent, on a $L_{\frac{1}{k}\rho}^*\mu=\mu$.

Par ailleurs, d'après le lemme \ref{maximum}, on a $\mu(3)=0$. Par
conséquent, comme $h\circ f=hk$, on a
$L^*_\rho(\frac{1}{h}\mu)=\frac{1}{h}L_{\frac{1}{k}\rho}^*(\mu)=\frac{1}{h}\mu$.

La mesure borélienne $\frac{1}{h}\mu$ sur $\mathbb R$ est concentrée
sur le spectre de $\Delta$. Or, d'après la proposition
\ref{espacep02}, pour tous $x$ dans $\bigcup_{n\in\mathbb
N}f^{-n}(0)$ et $\varphi$ dans $\ell^2(\Gamma)$ tel que
$\Delta\varphi=x\varphi$, on a $\varphi(p_0)=\varphi(p_0^\vee)$ et,
donc, $\langle\varphi,\varphi_0\rangle=0$. Par conséquent, on a
$\mu\left(\bigcup_{n\in\mathbb N}f^{-n}(0)\right)=0$ et, d'après le
corollaire \ref{spectre}, $\mu$ est concentrée sur $\Lambda$.

La fonction $\rho$ est höldérienne et strictement positive sur
$\Lambda$ et on a $L_\rho(1)=1$ sur $\Lambda$. D'après le lemme
\ref{transfert}, il existe une unique mesure borélienne de
probabilité $\nu_\rho$ sur $\Lambda$ telle que
$L_\rho^*(\nu_\rho)=\nu_\rho$ et, pour toute fonction continue $g$
sur $\Lambda$, la suite $(L_\rho^n(g))_{n\in\mathbb N}$ converge
uniformément vers la fonction constante de valeur
$\int_{\Lambda}g\de\nu_\rho$. Montrons que la mesure borélienne
positive $\frac{1}{h}\mu$ est finie et, donc, proportionnelle à
$\nu_\rho$. Donnons-nous une fonction continue positive $g$ sur
$\Lambda$, nulle au voisinage de $3$, telle qu'on ait
$0<\int_\Lambda\frac{1}{h}g\de\mu<\infty$. Il existe un entier $n$
et un réel $\varepsilon>0$ tels que, pour tout $x$ dans $\Lambda$,
on ait $L^n_\rho(g)(x)\geq\varepsilon$. Comme
$\int_\Lambda\frac{1}{h}g\de\mu=\int_\Lambda\frac{1}{h}L^n_\rho(g)\de\mu$,
il vient  $\int_\Lambda\frac{1}{h}\de\mu<\infty$. Par conséquent, la
mesure $\frac{1}{h}\mu$ est un multiple de $\nu_\rho$. Or, on a
$\mu(\Lambda)=\N{\varphi_0}_2^2=2$ et, par un calcul immédiat,
$L_\rho(h)=2$, si bien que $\int_{\Lambda}h\de\nu_\rho=2$. Il vient
donc bien $\mu=h\nu_\rho$.
\end{demo}

Pour tout polynôme $p$ dans $\mathbb C[X]$, notons $\hat{p}$ la
fonction $p(\Delta)\varphi_0$ sur $\Gamma$. Par définition,
l'application $g\mapsto\hat{g}$ se prolonge en une isométrie de
${\rm L}^2(h\nu_\rho)$ dans le sous-espace cyclique $\Phi$ de
$\ell^2(\Gamma)$ engendré par $\varphi_0$. Notons $l$ la fonction
$x\mapsto x$ sur $\Lambda$. On a la

\begin{Prop} \label{spectrediffus2} Le sous-espace $\Phi$ est stable par les opérateurs $\Delta$, $\Pi$ et $\Pi^*$.
Pour tout $g$ dans ${\rm L}^2(h\nu_\rho)$, on a
\begin{align*}\Delta\hat{g}&=\widehat{lg}\\
\Pi\hat{g}&=\widehat{L_\rho g}\\
\Pi^*\hat{g}&=\widehat{k(g\circ f)}.\end{align*}
\end{Prop}

\begin{demo} Par définition, $\Phi$ est stable par $\Delta$ et on a la formule concernant $\Delta$.

Par un calcul direct, on montre qu'on a $L_\rho(l)=1$. Soit $n$ dans
$\mathbb N$. On a $L_\rho(f^n)=l^nL_\rho(1)=l^n$ et
$L_\rho(f^nl)=l^nL_\rho(l)=l^n$. Or, d'après le lemme
\ref{relation1}, on a
$\Pi(f(\Delta)^n\varphi_0)=\Delta^n\Pi\varphi_0$ et
$\Pi(f(\Delta)^n\Delta\varphi_0)=\Delta^n\Pi\Delta\varphi_0$ et,
donc, comme $\Pi\varphi_0=\Pi\Delta\varphi_0=\varphi_0$, l'espace
$\Phi$ est stable par $\Pi$ et, pour tout $p$ dans $\mathbb C[X]$,
$\Pi\hat{p}=\widehat{L_\rho p}$. Enfin, par convexité, pour toute
fonction mesurable $g$ sur $\Lambda$, on a $\abs{L_\rho(g)}^2\leq
L_\rho\left(\abs{g}^2\right)$, si bien que, pour $g$ dans ${\rm
L}^2(h\nu_\rho)$, on a
\begin{multline*}\int_\Lambda \abs{L_\rho(g)}^2h\de\nu_\rho
\leq \int_\Lambda  L_\rho\left(\abs{g}^2\right)h\de\nu_\rho
=\int_\Lambda \abs{g}^2 (h\circ f)\de\nu_\rho\\
=\int_\Lambda \abs{g}^2 kh\de\nu_\rho \leq 5 \int_\Lambda \abs{g}^2
h\de\nu_\rho,\end{multline*} donc l'opérateur $L_\rho$ est continu
dans  ${\rm L}^2(h\nu_\rho)$ et, par densité, pour tout $g$ dans
${\rm L}^2(h\nu_\rho)$, $\Pi\hat{g}=\widehat{L_\rho g}$.

Enfin, d'après les lemmes \ref{relation1} et \ref{relation3}, pour
tout polynôme $p$ dans $\mathbb C[X]$, on a
$\Pi^*(p(\Delta)\varphi_0)=p(f(\Delta))\Pi^*\varphi_0=p(f(\Delta))(\Delta+2)\varphi_0$.
Par conséquent, l'espace $\Phi$ est stable par $\Pi^*$ et, pour tout
$p$ dans $\mathbb C[X]$, $\Pi^*\hat{p}=\widehat{k (p\circ f)}$. Or,
pour tout $g$ dans ${\rm L}^2(h\nu_\rho)$, on a
\begin{multline*}\int_\Lambda \abs{k (g\circ f)}^2h\de\nu_\rho =\int_\Lambda
k\abs{g\circ f}^2(h\circ f)\de\nu_\rho\\ =\int_\Lambda
L_\rho(k)\abs{g}^2h\de\nu_\rho
=3\int_\Lambda\abs{g}^2h\de\nu_\rho\end{multline*} et, donc, par
densité, pour tout $g$ dans  ${\rm L}^2(h\nu_\rho)$,
$\Pi^*\hat{g}=\widehat{k(g\circ f)}$.
\end{demo}

Pour déterminer la structure spectrale complète de $\Delta$, nous
allons nous intéresser à d'autres éléments remarquables de
$\ell^2(\Gamma)$. Commen\c cons par noter $\psi_0$ la fonction sur
$\Gamma$ qui vaut $1$ en $p_0$ et en $p_0^\vee$ et $0$ partout
ailleurs. Nous avons le

\begin{Lem}\label{relation4} On a $\Pi^*\psi_0=\Delta\psi_0$.\end{Lem}

Nous en déduisons le

\begin{Cor} \label{spectrediscret1}
La mesure spectrale de $\psi_0$ est discrète. Plus précisément, la
fonction $\psi_0$ est contenue dans la somme directe des espaces
propres de $\Delta$ associés aux éléments de $\bigcup_{n\in\mathbb
N}f^{-n}(0)$.\end{Cor}

\begin{demo} Soit $\mu$ la restriction de la mesure spectrale de
$\psi_0$ à $\Lambda$. D'après le corollaire \ref{spectre}, il s'agit
de montrer qu'on a $\mu=0$.

Pour $x\notin\{0,\frac{1}{2}\}$, posons
$\tau(x)=\frac{(x+2)}{x(2x-1)}$ et $\sigma(x)=\frac{1}{x(2x-1)}$. La
fonction $\sigma$ est höldérienne et strictement positive sur
$\Lambda$. En raisonnant comme dans la démonstration du corollaire
\ref{spectrediffus1}, on voit que, comme $0\notin\Lambda$, on a,
d'après le lemme \ref{relation4}, $L_{\tau}^*\mu=\mu$. Or, comme
$h\circ f=hk$, pour tout $x$ dans $\mathbb R$, pour tout entier $n$,
on a $L_{\tau}^n(h)=hL^n_\sigma(1)$. Notons $\lambda_\sigma$ le
rayon spectral de $L_\sigma$ et $\nu_\sigma$ sa mesure d'équilibre,
comme dans le lemme \ref{transfert}. Par un calcul direct, on montre
que, pour tout $x$ dans $\Lambda$, on a
$L_\sigma(1)(x)=\frac{1}{x+3}$. En particulier, pour $x\neq -2$, on
a $L_\sigma(1)(x)<1$, si bien que $\lambda_\sigma=\int_\Lambda
L_\sigma(1)\de\nu_\sigma<1$ et, donc, la suite
$(L_\sigma^n(1))_{n\in\mathbb N}$ converge uniformément vers $0$ sur
$\Lambda$. Par conséquent, la suite $(L_{\tau}^n(h))_{n\in\mathbb
N}$ converge uniformément vers $0$ sur $\Lambda$ et on a
$\int_\Lambda h\de\mu=0$, si bien que $\mu(\Lambda-\{3\})=0$.
D'après le lemme \ref{maximum}, on a $\mu(3)=0$, et, donc, $\mu=0$.
\end{demo}

Notons que, en utilisant le lemme \ref{relation4}, on pourrait établir
une formule donnant, pour tout $x$ dans  $\bigcup_{n\in\mathbb N}f^{-n}(0)$,
la valeur de la norme de la projection de $\psi_0$ sur l'espace des fonctions propres de valeur propre $x$.

\'Etudions les invariants spectraux d'un dernier élément de
$\ell^2(\Gamma)$. Pour cela, notons $q_0$ et $r_0$ les deux voisins
de $p_0$ différents de $p_0^\vee$ et $\chi_0$ la fonction sur
$\Gamma$ qui vaut $1$ en $q_0$, $-1$ en $r_0$ et $0$ partout
ailleurs. De même, on note $q_0^\vee$ et $r_0^\vee$ les deux voisins
de $p_0^\vee$ différents de $p_0$ et $\chi_0^\vee$ la fonction sur
$\Gamma$ qui vaut $1$ en $q_0^\vee$, $-1$ en $r_0^\vee$ et $0$
partout ailleurs. On pourrait à nouveau remarquer qu'on a
$\Pi^*\chi_0=(\Delta^2+2\Delta)\chi_0$ et étudier la mesure
spectrale de $\chi_0$ en utilisant les mêmes méthodes que dans les
corollaires \ref{spectrediffus1} et \ref{spectrediscret1}. Nous
allons utiliser une autre approche, plus géométrique, analogue à
celle de la démonstration du lemme \ref{espacep-21}.

D'après le lemme \ref{triangleinfini}, il existe un unique
automorphisme $\iota$ du graphe $\Gamma$ tel que $\iota(q_0)=r_0$ et
que $\iota(q_0^\vee)=r_0^\vee$ et $\iota$ est une involution. Notons
$H$ l'espace des éléments $\varphi$  de $\ell^2(\Gamma)$ tels que
$\iota(\varphi)=-\varphi$ et $K$ (resp. $K^\vee$) le sous-espace de
$H$ formé des éléments de $H$ qui sont nuls sur le triangle infini
issu de $p_0^\vee$ (resp. de $p_0$). On a $H=K\oplus K^\vee$,
$\chi_0\in K$, $\chi_0^\vee\in K^\vee$ et les sous-espaces $K$ et
$K^\vee$ sont stables par les endomorphismes $\Delta$, $\Pi$ et
$\Pi^*$. Pour tout $n\geq 1$, notons $\mathcal T_n$ le $n$-triangle
issu de $p_0$ et $\mathcal T_n^\vee$ le $n$-triangle issu de
$p_0^\vee$. Les groupe de permutations $\mathfrak S(\partial
\mathcal T_n)$  et $\mathfrak S(\partial \mathcal T_n^\vee)$
agissent sur les triangles $\mathcal T_n$ et $\mathcal T_n^\vee$. On
note $K_n$ (resp. $K_n^\vee$) l'espace des fonctions $\varphi$ sur
$\mathcal T_n$ (resp. $\mathcal T_n^\vee$) telles que, pour tout $s$
dans $\mathfrak S(\partial \mathcal T_n)$ (resp. dans $\mathfrak
S(\partial \mathcal T_n^\vee)$), on ait
$s\varphi=\varepsilon(s)\varphi$, où $\varepsilon$ est le morphisme
de signature. On identifie $K_n$ et $K_n^\vee$ à des sous-espaces de
dimension finie de $K$ et $K^\vee$. On a alors $\Delta K_n\subset
K_n$, $\Pi^* K_n\subset K_{n+1}$ et, si $n\geq 2$,  $\Pi K_n\subset
K_{n-1}$, et les identités analogues dans $K^\vee$.

Nous avons le

\begin{Lem} \label{finidense}
Les espaces $K$ et $K^\vee$ sont topologiquement engendrés par les
ensembles $\bigcup_{n\geq 1} K_n$ et $\bigcup_{n\geq 1} K_n^\vee$.
\end{Lem}

\begin{demo} Soit $\varphi$ une fonction dans $K$. Pour tout entier $n\geq 2$,
on note $\varphi_n$ l'unique élément de $K_n$ qui est égal à $\varphi$ sur $\mathcal T_{n-1}$.
On a  $\N{\varphi_n}_2\leq\sqrt{3}\N{\varphi}_2$. Alors, pour tout
$\varphi$, la suite $(\varphi_n)$ tend faiblement vers $\varphi$
dans $\ell^2(\Gamma)$. Donc l'ensemble  $\bigcup_{n\geq 1} K_n$ est
faiblement dense dans $H$ et l'espace vectoriel qu'il engendre est,
par conséquent, fortement dense. Le résultat pour $K^\vee$ s'en
déduit par symétrie.
\end{demo}

\begin{Cor} \label{spectrediscret2}
Le spectre de $\Delta$ dans $H$ est discret. Ses valeurs propres sont exactement les éléments de l'ensemble
$\bigcup_{n\in\mathbb N}f^{-n}(-2)\cup\bigcup_{n\in\mathbb N}f^{-n}(0)$.\end{Cor}

\begin{demo} Comme, pour tout $n$, les sous-espaces $K_n$ et $K_n^\vee$ sont stables par $\Delta$ et de dimension finie,
le caractère discret du spectre de $\Delta$ dans $H$ est une
conséquence immédiate du lemme \ref{finidense}. La détermination
exacte des valeurs propres s'obtient en raisonnant comme dans la
section \ref{secspectre}. Une formule pour le polynôme
caractéristique de $\Delta$ dans $K_n$ est donnée à la proposition
\ref{spectrediscret1comp}.\end{demo}

Nous pouvons à présent terminer la démonstration du théorème \ref{spectrePascal} grâce à la

\begin{Prop} \label{spectrePascalfinal}
Soit $\Phi^\perp$ l'orthogonal dans $\ell^2(\Gamma)$ du sous-espace
cyclique $\Phi$ engendré par $\varphi_0$. Alors le spectre de
$\Delta$ dans $\Phi^\perp$ est discret et l'ensemble de ses valeurs
propres est exactement $\bigcup_{n\in\mathbb
N}f^{-n}(-2)\cup\bigcup_{n\in\mathbb N}f^{-n}(0)$.
\end{Prop}

Avant de démontrer cette proposition, établissons un résultat
prélimi\-naire. Pour $\varphi$ et $\psi$ dans $\ell^2(\Gamma)$,
notons $\mu_{\varphi,\psi}$ l'unique mesure complexe boré\-lienne
sur $\mathbb R$ telle que, pour tout polynôme $p$ dans $\mathbb
C[X]$, on ait $\int_{\mathbb R}p\de\mu_{\varphi,\psi}=\langle
p(\Delta)\varphi,\psi\rangle$. On a le

\begin{Lem} \label{poussemesure}
Pour tous $\varphi$ et $\psi$ dans  $\ell^2(\Gamma)$, on a $\mu_{\Pi\varphi,\psi}=f_*\mu_{\varphi,\Pi^*\psi}$.\end{Lem}

\begin{demo} Pour $p$ dans $\mathbb C[X]$, on a, d'après le lemme \ref{relation1},
$$\int_{\mathbb R}p\de\mu_{\Pi\varphi,\psi}=\langle p(\Delta)\Pi\varphi,\psi\rangle
=\langle p(f(\Delta))\varphi,\Pi^*\psi\rangle=\int_{\mathbb
R}(p\circ f)\de\mu_{\varphi,\Pi^*\psi}.$$
\end{demo}

\begin{demo}[Démonstration de la proposition \ref{spectrePascalfinal}]
D'après les corollaires \ref{espacep03} et \ref{espacep-22},
les espaces propres associés aux éléments de $\bigcup_{n\in\mathbb
N}f^{-n}(-2)\cup\bigcup_{n\in\mathbb N}f^{-n}(0)$ sont non-triviaux.
Notons $P$ le projecteur orthogonal sur $\Phi^\perp$ dans
$\ell^2(\Gamma)$. D'après la proposition \ref{spectrediffus2},
l'opérateur $P$ commute à $\Delta$, $\Pi$ et $\Pi^*$. Pour démontrer
la proposition, il suffit d'établir que, pour tout $\varphi$  à
support fini, pour tout $\psi$ dans $\ell^2(\Gamma)$, la mesure
$\mu_{P\varphi,\psi}$ est atomique et concentrée sur l'ensemble
$\bigcup_{n\in\mathbb N}f^{-n}(-2)\cup\bigcup_{n\in\mathbb
N}f^{-n}(0)$.

Soient toujours $q_0$, $r_0$, $q_0^\vee$ et $r_0^\vee$ les voisins
de $p_0$ et $p_0^\vee$ et, pour tout entier $n$, $\mathcal T_n$ le
$n$-triangle issu de $p_0$ et $\mathcal T_n^\vee$ le $n$-triangle
issu de $p_0^\vee$. On note $L_n$ l'espace des fonctions sur
$\Gamma$ dont le support est contenu dans la réunion de $\mathcal
T_n$, de $\mathcal T_n^\vee$, et des voisins des sommets de
$\mathcal T_n$ et de $\mathcal T_n^\vee$. On a, pour $n\geq 1$, $\Pi
L_n\subset L_{n-1}$ et $\Pi\Delta L_n\subset L_{n-1}$. Montrons, par
récurrence sur $n$, que, pour toute fonction $\varphi$ dans $L_n$,
pour tout $\psi$ dans $\ell^2(\Gamma)$, la mesure
$\mu_{P\varphi,\psi}$ est atomique et concentrée sur l'ensemble
$\bigcup_{n\in\mathbb N}f^{-n}(-2)\cup\bigcup_{n\in\mathbb
N}f^{-n}(0)$.

Pour $n=0$, $L_0$ est l'espace des fonctions qui sont nulles en
dehors de l'ensemble $\{p_0,q_0,r_0,p_0^\vee,q_0^\vee,r_0^\vee\}$.
On vérifie aisément que cet espace est engendré par les fonctions
$\varphi_0$, $\Delta\varphi_0$, $\psi_0$, $\Delta\psi_0$, $\chi_0$
et $\chi_0^\vee$. Dans ce cas, la description des mesures spectrales
découle immédiatement des corollaires \ref{spectrediffus1},
\ref{spectrediscret1} et \ref{spectrediscret2}.

Si le résultat est vrai pour un entier $n$, donnons-nous $\varphi$
dans $L_{n+1}$. Alors, les fonctions $\Pi\varphi$ et
$\Pi\Delta\varphi$ sont dans $L_n$ et, par récurrence, pour tout
$\psi$ dans  $\ell^2(\Gamma)$, les mesures $\mu_{\Pi
P\varphi,\psi}=\mu_{P\Pi \varphi,\psi}$ et $\mu_{\Pi\Delta
P\varphi,\psi}=\mu_{P\Pi\Delta\varphi,\psi}$ sont atomiques et
concentrées sur l'ensemble $\bigcup_{n\in\mathbb
N}f^{-n}(-2)\cup\bigcup_{n\in\mathbb N}f^{-n}(0)$. D'après le lemme
\ref{poussemesure}, les mesures $\mu_{P\varphi,\Pi^*\psi}$ et
$\mu_{\Delta P\varphi,\Pi^*\psi}=\mu_{P\varphi,\Delta\Pi^*\psi}$
sont donc atomiques et concentrées sur l'ensemble $\bigcup_{n\geq
1}f^{-n}(-2)\cup\bigcup_{n\geq 1}f^{-n}(0)$. Or, d'après le lemme
\ref{spectreresiduel}, le spectre de $\Delta$ dans l'orthogonal du
sous-espace de
 $\ell^2(\Gamma)$ engendré par l'image de $\Pi^*$ et par celle de $\Delta\Pi^*$ est
égal à $\{-2,0\}$. Par conséquent, pour tout $\psi$ dans
$\ell^2(\Gamma)$, la mesure $\mu_{P\varphi,\psi}$ est atomique et
concentrée sur l'ensemble $\bigcup_{n\in\mathbb
N}f^{-n}(-2)\cup\bigcup_{n\in\mathbb N}f^{-n}(0)$. Le résultat en
découle.
\end{demo}

\section{Quotients finis de $\Gamma$}
\label{graphesfinis}

Dans ce paragraphe, nous appliquons les méthodes développées
précédem\-ment à la description du spectre de certains graphes finis
fortement reliés à $\Gamma$.

Soient $\Phi$ et $\Psi$ des graphes. Nous dirons qu'une application
$\varpi:\Phi\rightarrow\Psi$ est un revêtement si, pour tout $p$
dans $\Phi$, l'application $\varpi$ induit une bijection de
l'ensemble des voisins de $p$ dans celui des voisins de $\varpi(p)$.
La composition de deux revêtements est un revêtement. Si $\Phi$ et
$\Psi$ sont des graphes réguliers de valence $3$ et si
$\varpi:\Phi\rightarrow\Psi$ est un revêtement, il existe un unique
revêtement $\hat{\varpi}:\hat{\Phi}\rightarrow\hat{\Psi}$ tel que
$\Pi\hat{\varpi}=\varpi\Pi$. Réciproquement, en raisonnant comme
dans le lemme \ref{automorphismes}, on montre que tout revêtement
$\hat{\Phi}\rightarrow\hat{\Psi}$ est de cette forme.

Fixons quatre éléments distincts, $a$, $b$, $c$ et $d$. On note
$\Gamma_{0}$ le graphe obtenu en munissant l'ensemble $\{a,b,c,d\}$
de la relation qui lie tous les points distincts : c'est un graphe
régulier de valence $3$. Le groupe de ses automorphismes est égal au
groupe $\mathfrak S(a,b,c,d)$ des permutations de l'ensemble
$\{a,b,c,d\}$.

\begin{Lem} \label{revetements}
Soient $\Phi$ un graphe régulier de valence $3$ et
$\varpi:\Phi\rightarrow \Gamma_0$ un revêtement. Alors,
l'application
$\tilde{\varpi}:\hat{\Phi}\rightarrow\Gamma_0,(p,q)\mapsto\varpi(q)$
est un revêtement. L'application $\varpi\mapsto\tilde{\varpi}$ est
une bijection $\mathfrak S(a,b,c,d)$-équivariante de l'ensemble des
revêtements $\Phi\rightarrow \Gamma_0$ dans l'ensemble des
revêtements $\hat{\Phi}\rightarrow \Gamma_0$.
\end{Lem}

La construction du revêtement $\tilde{\varpi}$ est représentée par la figure \ref{nouveaurevetement}.

\begin{figure}\begin{center}\input{nouveaurevetement.pstex_t}\caption{Construction du revêtement $\tilde{\varpi}$}
\label{nouveaurevetement}\end{center}\end{figure}

\begin{demo} Soit $p$ un point de $\Phi$ et soient $q$, $r$ et $s$ les voisins de $p$.
Quitte à permuter les éléments de $\{a,b,c,d\}$,
supposons que l'on a $\varpi(p)=a$, $\varpi(q)=b$,  $\varpi(r)=c$ et
$\varpi(s)=d$. Alors, on a $\tilde{\varpi}(p,q)=b$,
$\tilde{\varpi}(q,p)=a$, $\tilde{\varpi}(p,r)=c$ et
$\tilde{\varpi}(p,s)=d$ et, donc, $\tilde{\varpi}$ est bien un
revêtement.

Réciproquement, soit $\omega:\hat{\Phi}\rightarrow \Gamma_0$ un revêtement.
Soit toujours $p$ un point de $\Phi$, de voisins $q$, $r$ et $s$.
\`A nouveau, quitte à permuter, supposons qu'on a $\omega(p,q)=b$,  $\omega(p,r)=c$
et $\omega(p,s)=d$. Alors, comme $\omega$ est un revêtement, on a nécessairement
$\omega(q,p)=\omega(r,p)=\omega(s,p)=a$. Il existe donc une application
$\varpi:\Phi\rightarrow\Gamma_0$ telle que, pour tous $p$ et $q$ dans $\Phi$ avec $p\sim q$, on ait
$\omega(q,p)=\varpi(p)$. Par construction, $\varpi$ est un revêtement et on a $\tilde{\varpi}=\omega$.
\end{demo}

Pour tout entier naturel $n$, on note
$\Gamma_n=\hat{\Gamma}_0^{(n)}$ le graphe obtenu en rempla\c cant
les points de $\Gamma_0$ par des $n$-triangles. D'après le lemme
\ref{automorphismes}, le groupe des automorphismes de $\Gamma_n$
s'identifie naturellement à $\mathfrak S(a,b,c,d)$. Du lemme
\ref{revetements}, on déduit le

\begin{Cor} Pour tous entiers naturels $n\leq m$, il existe des revêtements $\Gamma_m\rightarrow\Gamma_n$.
Le groupe $\mathfrak S(a,b,c,d)$ agit simplement transitivement
sur l'ensemble de ces revêtements.
\end{Cor}

\begin{demo} Les revêtements $\Gamma_0\rightarrow\Gamma_0$ sont
exactement les bijections de $\Gamma_0$ dans lui-même et le
corollaire est donc vrai pour $m=n=0$. Par récurrence, d'après le
lemme \ref{revetements}, le corollaire est vrai pour tout entier $m$
et pour $n=0$. Enfin, comme, si $\Phi$ et $\Psi$ sont des graphes
réguliers de valence $3$, les revêtements $\Phi\rightarrow\Psi$ et
$\hat{\Phi}\rightarrow\hat{\Psi}$ sont en bijection naturelle, à
nouveau par récurrence, le corollaire est vrai pour tous entiers
$m\geq n$.\end{demo}

Revenons à présent à $\Gamma$.
Notons $q_0$ et $r_0$ les deux voisins de $p_0$ différents de $p_0^\vee$
et $q_0^\vee$ et $r_0^\vee$ les deux voisins de $p_0^\vee$ différents de $p_0$.
Nous avons le

\begin{Lem}\label{quotientPascal} Il existe un unique revêtement $\varpi:\Gamma\rightarrow\Gamma_0$
tel que $\varpi(p_0)=a$, $\varpi(q_0)=c$,
$\varpi(r_0)=d$, $\varpi(p_0^\vee)=b$, $\varpi(q_0^\vee)=c$ et $\varpi(r_0^\vee)=d$.
\end{Lem}

Ce revêtement est représenté par la figure \ref{revetement}.

\begin{demo} Soient $n$ un entier $\geq 1$ ou l'infini et $\mathcal T$ un $n$-triangle.
Soit $\varpi:\mathcal T\rightarrow\Gamma_0$. Disons que $\varpi$ est
un quasi-revêtement si, pour tout point $p$ de $\mathcal T-\partial
\mathcal T$, $\varpi$ induit une bijection de l'ensemble des voisins
de $\mathcal T$ dans $\Gamma_0-\{\varpi(p)\}$ et si, pour tout point
$p$ de $\partial \mathcal T$, les valeurs de $\varpi$ sur les
voisins de $p$ sont des éléments distincts de
$\Gamma_0-\{\varpi(p)\}$. Dans ce cas, on note encore
$\tilde{\varpi}$ l'application $\hat{\mathcal T}\rightarrow\Gamma_0$
telle que, pour tous $p$ et $q$ dans $\mathcal T$ avec $p\sim q$, on
ait $\tilde{\varpi}(p,q)=\varpi(q)$ et que, pour tout $p$ dans
$\partial \mathcal T=\partial\hat{\mathcal T}$, si les voisins de
$p$ dans $\mathcal T$ sont $q$ et $r$, $\tilde{\varpi}(p)$ soit
l'unique élément de $\Gamma_0-\{\varpi(p),\varpi(q),\varpi(r)\}$. En
raisonnant comme dans la démonstration du lemme \ref{revetements},
on vérifie aisément que l'application $\varpi\mapsto\tilde{\varpi}$
est une bijection $\mathfrak S(a,b,c,d)$-équivariante de l'ensemble
des quasi-revêtements $\mathcal T\rightarrow \Gamma_0$ dans
l'ensemble des quasi-revêtements $\hat{\mathcal T}\rightarrow
\Gamma_0$.

Par conséquent, pour tout $n\geq 1$, si $\mathcal T_n$ est le
$n$-triangle de $\Gamma$ contenant $p_0$, il existe un unique
quasi-revêtement $\varpi_n$ de $\mathcal T_n$ dans $\Gamma_0$ tel
que $\varpi_n(p_0)=a$, $\varpi_n(q_0)=c$ et  $\varpi_n(r_0)=d$. Par
unicité, $\varpi_n$ et $\varpi_{n+1}$ sont égaux sur $\mathcal T_n$.
Il existe donc un unique quasi-revêtement $\varpi_\infty$ du
triangle infini $\mathcal T_\infty$ issu de  $p_0$ dans $\Gamma_0$
tel que $\varpi_\infty(p_0)=a$, $\varpi_\infty(q_0)=c$ et
$\varpi_\infty(r_0)=d$. De même, si $\mathcal T^\vee_\infty$ désigne
le triangle infini issu de  $p_0^\vee$, il existe un unique
quasi-revêtement $\varpi_\infty^\vee$ de $\mathcal T_\infty^\vee$
dans $\Gamma_0$ tel que $\varpi_\infty(p_0^\vee)=b$,
$\varpi_\infty(q_0^\vee)=c$ et $\varpi_\infty(r_0^\vee)=d$.
L'application $\varpi:\Gamma\rightarrow\Gamma_0$ dont la restriction
à $\mathcal T_\infty$ est $\varpi_\infty$ et la restriction à
$\mathcal T_\infty^\vee$ est $\varpi_\infty^\vee$ est donc bien
l'unique revêtement de $\Gamma$ dans $\Gamma_0$ ayant les propriétés
demandées.
\end{demo}

\`A nouveau, des lemmes \ref{revetements} et \ref{quotientPascal},
on déduit le

\begin{Cor} Pour tout entier naturel $n$, il existe des revêtements $\Gamma\rightarrow\Gamma_n$.
Le groupe $\mathfrak S(a,b,c,d)$ agit simplement sur l'ensemble de
ces revêtements. Cette action possède deux orbites : d'une part,
l'ensemble de revêtements $\varpi$ tels que
$\varpi(q_0)=\varpi(q_0^\vee)$ et, d'autre part, l'ensemble de
revêtements $\varpi$ tels que $\varpi(q_0)=\varpi(r_0^\vee)$.
\end{Cor}

Nous allons à présent décrire, pour tout entier $n$, la théorie
spectrale du graphe $\Gamma_n$. Notons $k$ le polynôme $X+2$, $l$ le
polynôme $X$, $m$ le polynôme $X-2$, et, comme toujours, $f$ le
polynôme $X^2-X-3$. Les méthodes développées dans la section
\ref{secspectre} et la section \ref{valeurp} nous permettent de
démontrer la

\begin{Prop} \label{spectrePascalfinifinal}
Pour tout entier naturel $n$, le polynôme caractéristique de $\Delta$ dans $\ell^2(\Gamma_n)$ est
$$(X-3)(X+1)^3\prod_{p=0}^{n-1} (m\circ f^p(X))^3(l\circ f^p(X))^{2.3^{n-1-p}}(k\circ f^p(X))^{1+2.3^{n-1-p}}.$$
\end{Prop}

Rappelons que, à la section \ref{secspectre}, nous avons introduit
la notion de graphe partageable. La démonstration utilise le

\begin{Lem} \label{evite-3}
Soit $\Phi$ un graphe connexe régulier de valence $3$. Le graphe
$\hat{\Phi}$ n'est pas partageable. En particulier, pour tout entier
naturel $n$, le graphe $\Gamma_n$ n'est pas partageable.\end{Lem}

\begin{demo} Comme tout point de $\hat{\Phi}$ est contenu dans un
$1$-triangle, tout point peut être joint à lui-même par un chemin de
longueur impaire et, donc, $\hat{\Phi}$ n'est pas partageable. De
même, tout point de $\Gamma_0$ peut être joint à lui-même par un
chemin de longueur impaire. Le lemme en découle.
\end{demo}

\begin{demo}[Démonstration de la proposition \ref{spectrePascalfinifinal}]
Nous allons démontrer ce résultat par récurrence sur $n$. Pour
$n=0$, l'espace $\ell^2(\Gamma_0)$ est de dimension $4$ et, pour
l'action naturelle du groupe  $\mathfrak S(a,b,c,d)$, il est la
somme de deux sous-espaces irréductibles non isomorphes, l'espace
des fonctions constantes et l'espace des fonctions $\varphi$ telles
que $\varphi(a)+\varphi(b)+\varphi(c)+\varphi(d)=0$. L'opérateur
$\Delta$ commute à l'action de  $\mathfrak S(a,b,c,d)$ et laisse
donc stables ces deux sous-espaces. Dans le premier, il agit par
multiplication par $3$ et, dans le deuxième, par multiplication par
$-1$. Son polynôme caractéristique est donc $(X-3)(X+1)^3$.

Supposons le résultat démontré pour $n$. D'après le lemme
\ref{evite-3}, $\Gamma_{n}$ n'est pas partageable. Par conséquent,
si $H$ est le sous-espace vectoriel de $\ell^2(\Gamma_{n+1})$
engendré par l'image de $\Pi^*$ et par celle de $\Delta\Pi^*$,
d'après le corollaire \ref{transformespectre1} et le lemme
\ref{transformevalp}, le polynôme caractéristique de $\Delta$ dans
l'orthogonal dans $H$ des fonctions constantes est
$$(f(X)+1)^3\prod_{p=0}^{n-1} (m\circ f^{p+1}(X))^3(l\circ f^{p+1}(X))^{2.3^{n-1-p}}(k\circ f^{p+1}(X))^{1+2.3^{n-1-p}}$$
et, donc, comme $f(X)+1=(X+1)(X-2)$, le polynôme caractéristique de $\Delta$ dans $H$ s'écrit
$$(X-3)(X+1)^3\prod_{p=0}^{n} (m\circ f^{p}(X))^3\prod_{p=1}^{n} (l\circ f^{p}(X))^{2.3^{n-p}}(k\circ f^{p}(X))^{1+2.3^{n-p}}.$$

Il nous reste à déterminer les dimensions des espaces propres de
valeur propre $0$ et $-2$ dans l'orthogonal de $H$ dans
$\ell^2(\Gamma_{n+1})$. Ceux-ci sont décrits par le lemme
\ref{spectreresiduel}. Or, si $n\geq 1$, les $2$-triangles de
$\Gamma_{n+1}$ sont les images inverses par $\Pi^2$ des points de
$\Gamma_{n-1}$ et tout point de $\Gamma_{n+1}$ appartient à un
unique $2$-triangle. En raisonnant comme dans le lemme
\ref{espacep01}, on voit alors que l'espace propre de valeur propre
$0$ de $\Delta$ dans $\ell^2(\Gamma_{n+1})$ est isomorphe à l'espace
des fonctions sur les arêtes de $\Gamma_{n-1}$. Comme $\Gamma_{n-1}$
est un graphe régulier de valence $3$ et qu'il contient $4.3^{n-1}$
points, il contient $2.3^n$ arêtes et l'espace propre de valeur
propre $0$ de $\Delta$ dans $\ell^2(\Gamma_{n+1})$ est de dimension
$2.3^n$. Si $n=0$, en utilisant la caractérisation du lemme
\ref{spectreresiduel}, on vérifie par un calcul direct que la
dimension de l'espace propre de valeur propre $0$ dans
$\ell^2(\Gamma_1)$ est $2$. Alors, comme, d'après le corollaire
\ref{transformespectre1} et le lemme \ref{transformevalp}, $H$ est
de dimension $2\dim\ell^2(\Gamma_{n})-1=8.3^{n}-1$, l'orthogonal de
$H$ et de l'espace propre de valeur propre $0$ est de dimension
$4.3^{n+1}-(8.3^{n}-1)-2.3^n=2.3^n+1$. D'après le lemme
\ref{spectreresiduel}, cet espace est l'espace propre de valeur
propre $-2$ de $\Delta$ et le polynôme caractéristique de $\Delta$
dans $\ell^2(\Gamma_{n+1})$ a bien la forme donnée par l'énoncé.
\end{demo}

\section{Le compactifié plan de $\Gamma$}
\label{compactifie}

Nous considérons désormais l'ensemble $X$ des éléments
$(p_{k,l})_{(k,l)\in\mathbb Z^2}$ de $(\mathbb Z/2\mathbb
Z)^{\mathbb Z^2}$ tels que, pour tous entiers $k$ et $l$, on ait
$p_{k,l}+p_{k+1,l}+p_{k,l+1}=0$ dans $\mathbb Z/2\mathbb Z$. C'est
un espace topologique compact pour la topologie induite par la
topologie produit. On note $T$ et $S$ les deux applications de $X$
dans $X$ telles que, pour tout $p$ dans $X$, on ait
$Tp=(p_{k+1,l})_{(k,l)\in\mathbb Z^2}$ et
$Sp=(p_{k,l+1})_{(k,l)\in\mathbb Z^2}$. Les homéomorphismes $T$ et
$S$ engendrent l'action naturelle de $\mathbb Z^2$ sur $X$.

Pour $p$ dans $X$ et $k$ et $l$ dans $\mathbb Z$, on a
$p_{k,l}+p_{k-1,l+1}+p_{k-1,l}=0$ et
$p_{k,l}+p_{k,l-1}+p_{k+1,l-1}=0$. Or, le sous-groupe fini
$\mathfrak S$ de $\GL{2}{\mathbb Z}$, engendré par les matrices
$\begin{pmatrix}-1&-1\\0&1\end{pmatrix}$ et
$\begin{pmatrix}0&1\\-1&-1\end{pmatrix}$ permute les trois paires de
vecteurs $\{(1,0),(0,1)\}$, $\{(-1,1),(-1,0)\}$ et
$\{(0,-1),(1,-1)\}$ de $\mathbb Z^2$. En particulier, le groupe
$\mathfrak S$ agit naturellement sur $X$ : pour tous $p$ dans $X$ et
$s$ dans $\mathfrak S$, pour tous $k$ et $l$ dans $\mathbb Z^2$, on
a $(sp)_{k,l}=p_{s^{-1}(k,l)}$. L'action de $\mathfrak S$ sur les
trois paires de vecteurs $\{(1,0),(0,1)\}$, $\{(-1,1),(-1,0)\}$ et
$\{(0,-1),(1,-1)\}$ identifie $\mathfrak S$ et le groupe des
permutations de cet ensemble à trois éléments.

Soit $Y$ l'ensemble des points $p$ de $X$ tels que $p_{0,0}=1$.
L'ensemble $Y$ est stable par l'action de $\mathfrak S$. Pour tout
$p$ dans $Y$, on note $Y_p$ l'ensemble des points de l'orbite de $p$
sous l'action de $\mathbb Z^2$ qui appartiennent à $Y$. Si $p$ est
un point de $Y$, on a $p_{1,0}+p_{0,1}=1$ dans $\mathbb Z/2\mathbb
Z$ et, donc, un et un seul des points $Tp$ et $Sp$ appartient à $Y$.
De même, un et un seul des points $T^{-1}p$ et $T^{-1}Sp$ appartient
à $Y$ et un et un seul des points $S^{-1}p$ et $TS^{-1}p$ appartient
à $Y$.  Pour $p$ et $q$ dans $Y$, notons $p\sim q$ si $q$ appartient
à l'ensemble $\{Tp,Sp,T^{-1}Sp,T^{-1}p,S^{-1}p,TS^{-1}p\}$. Cette
relation est symétrique et $\mathfrak S$-invariante.

De même, si $p$ appartient à $Y$, on note $\tilde{Y}_p$ l'ensemble
des $(k,l)$ dans $\mathbb Z^2$ tels que $p_{k,l}=1$ et, pour $(i,j)$
et $(k,l)$ dans $\tilde{Y}_p$ on note $(i,j)\sim (k,l)$ si
$(i-k,j-l)$ appartient à l'ensemble
$\{(1,0),(0,1),(-1,1),(-1,0),(0,-1),(1,-1)\}$. Alors, $\tilde{Y}_p$
est un graphe régulier de valence $3$. Si le stabilisateur dans
$\mathbb Z^2$ de $p$ est trivial, $Y_p$ est un graphe régulier de
valence $3$ et l'application naturelle $\tilde{Y_p}\rightarrow Y_p$
est un isomorphisme de graphes.

Notons $u$ l'unique élément de $(\mathbb Z/2\mathbb Z)^{\mathbb
Z^2}$ tel que $u_{k,l}=0$ si et seulement si $k-l$ est congru à $0$
modulo $3$. L'élément $u$ est périodique sous l'action de $\mathbb
Z^2$ et son stabilisateur est l'ensemble des $(k,l)$ dans $\mathbb
Z^2$ tels que $k-l$ soit congru à $0$ modulo $3$. On vérifie que $u$
appartient à $X$. Son orbite sous l'action de $\mathbb Z^2$ est
égale à $\{u,Tu,Su\}$ et elle est stable par l'action de $\mathfrak
S$ : les éléments de signature $1$ de $\mathfrak S$ fixent $u$, $Tu$
et $Su$ et les éléments de signature $-1$ fixent $u$ et échangent
$Tu$ et $Su$. Nous avons le

\begin{Lem}\label{grapheinduit}
Soit $p$ dans $Y$. Le graphe $\tilde{Y}_p$ est connexe. Si $p$ est
différent de $Tu$ et de $Su$, l'ensemble $Y_p$, muni de la relation
$\sim$, est un graphe connexe régulier de valence $3$ et
l'application naturelle $\tilde{Y}_p\rightarrow Y_p$ est un
revêtement.\end{Lem}

\begin{demo} Montrons que $\tilde{Y}_p$ est connexe. Soit $(k,l)$
dans $Y_p$. Quitte à faire agir le groupe $\mathfrak S$ et à
échanger les rôles de $(0,0)$ et de $(k,l)$, on peut supposer $k$ et
$l$ positifs. Montrons dans ce cas, par récurrence sur $k+l$, que
$(k,l)$ appartient à la même composante connexe de $\tilde{Y}_p$ que
$(0,0)$. Si $k+l=0$, c'est immédiat. Supposons donc $k+l>0$ et
considérons les $p_{h,k+l-h}$ avec $0\leq h\leq k+l$. Si tous ces
éléments de $\mathbb Z/2\mathbb Z$ étaient égaux à $1$, on aurait,
pour tous entiers positifs $i$ et $j$ avec $i+j\leq k+l-1$,
$p_{i,j}=0$, ce qui est impossible, vu que $p_{0,0}=1$. Quitte à
encore faire agir $\mathfrak S$, on peut donc supposer qu'il existe
un entier naturel $0\leq i\leq l-1$ tel que, pour tout $0\leq j\leq
i$, on ait $p_{k+j,l-j}=1$, mais $p_{k+i+1,l-i-1}=0$. Cette situation est représentée à la figure \ref{chemin}. 
\begin{figure}\begin{center}\input{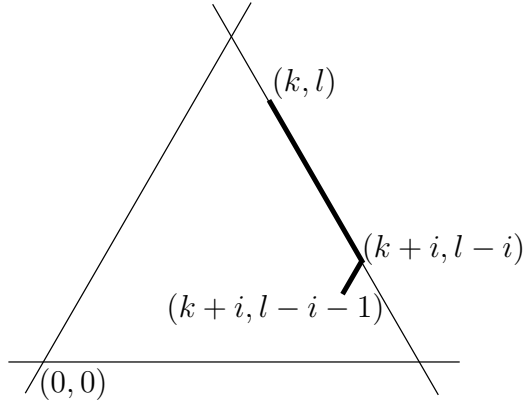}\caption{Connexité de  $\tilde{Y}_p$}
\label{chemin}\end{center}\end{figure}
Alors, les
points $(k+i,l-i)$ et $(k,l)$ appartiennent à la même composante
connexe de $\tilde{Y}_p$ et, comme $p_{k+i+1,l-i-1}=0$, on a
$p_{k+i,l-i-1}=1$ et $(k+i,l-i-1)$ appartient à $\tilde{Y}_p$. Comme
$(k+i)+(l-i-1)=k+l-1$, le résultat en découle par récurrence.

Comme l'application naturelle $\tilde{Y}_p\rightarrow Y_p$ est
surjective, pour conclure, il nous reste à montrer que, pour
$p\notin\{Tu,Su\}$, les points de l'ensemble
$\{p,Tp,Sp,T^{-1}Sp,T^{-1}p,S^{-1}p,TS^{-1}p\}$ sont deux à deux
distincts. Posons $V_p=\{Tp,Sp,T^{-1}Sp,T^{-1}p,S^{-1}p,TS^{-1}p\}$
et commen\-\c cons par supposer que $p$ appartient à $V_p$. Alors,
quitte à faire agir le groupe $\mathfrak S$, on peut supposer qu'on
a $p=Tp$ et $p_{0,-1}=1$ et, donc, $p_{1,-1}=p_{0,-1}+p_{0,0}=0$, ce
qui contredit le fait que $Tp=p$. On a donc $p\notin V_p$. Supposons
à présent que deux des éléments de l'ensemble $V_p$ sont égaux.
Quitte à encore faire agir $\mathfrak S$, on peut supposer qu'on a
$Tp=Sp$, $Tp=T^{-1}p$ ou $Tp=T^{-1}Sp$. Si $Tp=Sp$, on a
$S^{-1}Tp=p$ et on vient de montrer que c'est impossible. Si
$Tp=T^{-1}p$, on a $T^2p=p$ et la famille
$q=(p_{2k,2l})_{(k,l)\in\mathbb Z^2}$ appartient à $Y$ et vérifie
$Tq=q$ : à nouveau, on vient de montrer que c'est impossible. Enfin,
si $T^{-2}Sp=p$, supposons toujours, quitte à permuter, qu'on a
$p_{0,-1}=1$. Alors, on a $p_{1,-1}=0$, et, donc, comme
$T^{-2}Sp=p$, $p_{-1,0}=0$. De même, on a $p_{-2,0}=p_{0,-1}=1$ et
$p_{-1,-1}=p_{0,-1}+p_{-1,0}=1$. \`A nouveau, ceci implique
$p_{-3,0}=p_{-1,-1}=1$, $p_{-2,-1}=p_{-1,-1}+p_{-2,0}=0$ et, enfin,
$p_{-3,-1}=p_{-2,-1}+p_{-3,0}=1$, si bien que le point $q=T^{-3}p$
vérifie à nouveau $T^{-2}Sq=q$ et $q_{0,0}=q_{0,-1}=1$. Par
récurrence, on en déduit que, pour tout entier $k\leq 0$, on a
$p_{k,0}=1$ si $k$ est congru à $0$ ou $1$ modulo $3$ et que
$p_{k,0}=0$ si $k$ est congru à $2$ modulo $3$. En raisonnant de la
même fa\c con, on voit que $p_{-1,1}=p_{-1,0}+p_{0,0}=1$ et que,
comme  $T^{2}S^{-1}p=p$, $p_{1,0}=1$. Il vient
$p_{2,0}=p_{0,1}=p_{0,0}+p_{1,0}=0$, puis
$p_{3,0}=p_{1,1}=p_{1,0}+p_{2,0}=1$ et $p_{3,-1}=p_{1,0}=1$. Le
point $r=T^3p$ vérifie donc aussi $T^{-2}Sr=r$ et
$r_{0,0}=r_{0,-1}=1$, si bien que, pour tout $k$ dans $\mathbb Z$,
on a $p_{k,0}=0$ si et seulement si $k$ est congru à $2$ modulo $3$.
En particulier, la suite  $(p_{k,0})_{k\in\mathbb Z}$ est
$3$-périodique. Comme $T^{-2}Sp=p$, pour tout $l$ dans $\mathbb Z$,
la suite $(p_{k,l})_{k\in\mathbb Z}$ est $3$-périodique et, donc,
$T^3p=p$. Par conséquent, pour tous $k$ et $l$ dans $\mathbb Z$, si
$k-l$ est congru à $0$ modulo $3$, on a $T^kS^lp=p$. Comme on a
$p_{0,0}=u_{1,0}$, $p_{-1,0}=u_{-1,1}$ et $p_{-1,1}=u_{-1,2}$, il
vient $p=Tu$. Par conséquent, si $p$ n'appartient pas à $\{Tu,Su\}$,
la relation $\sim$ munit bien l'ensemble $Y_p$ d'une structure de
graphe régulier de valence $3$. Par définition, l'application
naturelle $\tilde{Y}_p\rightarrow Y_p$ est alors un revêtement. En
particulier, $Y_p$ est connexe.
\end{demo}

Soient $\varepsilon$ et $\eta$ dans $\{0,1\}$. Notons
$X^{(\varepsilon,\eta)}$ l'ensemble des éléments $p$ de $X$ tels
que, pour tous $k$ et $l$ dans $\mathbb Z$, si $(k,l)$ est congru à
$(\varepsilon,\eta)$ dans $(\mathbb Z/2\mathbb Z)^2$, on a
$p_{k,l}=0$. Si $p$ appartient à $X^{(\varepsilon,\eta)}$, pour tous
$k$ et $l$ dans $\mathbb Z$, on a
$p_{2k+1+\varepsilon,2l+1+\eta}=p_{2k+\varepsilon,2l+1+\eta}=p_{2k+1+\varepsilon,2l+\eta}$.
En particulier, pour $(\varepsilon',\eta')\neq(\varepsilon,\eta)$,
on a $X^{(\varepsilon',\eta')}\cap X^{(\varepsilon,\eta)}=\{0\}$ et,
si $p$ est un point de $Y^{(\varepsilon,\eta)}=Y\cap
X^{(\varepsilon,\eta)}$ (on a alors $(\varepsilon,\eta)\neq (0,0)$),
le point $p$ appartient à un triangle dans le graphe $Y_p$.

\begin{figure}\begin{center}\input{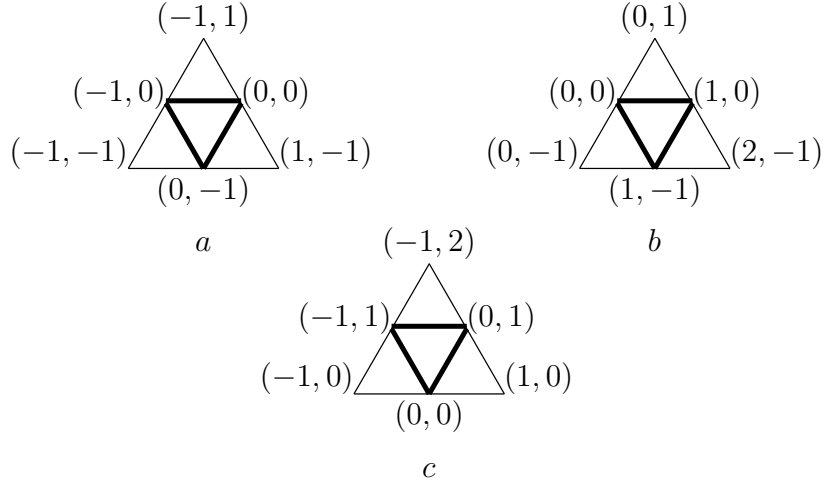}\caption{Les ensembles $Y^a$, $Y^b$ et $Y^c$}
\label{partition}\end{center}\end{figure}

Le groupe $\mathfrak S$ agit naturellement sur $(\mathbb Z/2\mathbb
Z)^2$ et, pour tout $s$ dans $\mathfrak S$, pour tout
$(\varepsilon,\eta)$ dans $(\mathbb Z/2\mathbb Z)^2$, on a
$X^{s(\varepsilon,\eta)}=sX^{(\varepsilon,\eta)}$. Dorénavant, on
notera $a=(1,1)$, $b=(0,1)$, $c=(1,0)$ et $\mathcal T_1=\{a,b,c\}$.
On considérera $\mathcal T_1$ comme un $1$-triangle. Le groupe
$\mathfrak S$ s'identifie au groupe de permutations $\mathfrak
S(a,b,c)$. On pose $\hat{Y}=Y^{a}\cup Y^{b}\cup Y^{c}$ : cette
réunion est disjointe et l'ensemble $\hat{Y}$ est invariant par
l'action de $\mathfrak S$. Les éléments de $Y^a$, $Y^b$ et $Y^c$ sont décrits par la figure \ref{partition}.

Soit $p$ un point de $Y$. Nous
noterons $\hat{p}^a$, $\hat{p}^b$ et $\hat{p}^c$ les points de
$(\mathbb Z/2\mathbb Z)^{\mathbb Z^2}$ tels que, pour tous $k$ et
$l$ dans $\mathbb Z$, on ait,\begin{cond}
\item $\hat{p}^{a}_{2k,2l}=\hat{p}^{a}_{2k-1,2l}=\hat{p}^{a}_{2k,2l-1}=p_{k,l}$
et $\hat{p}^{a}_{2k-1,2l-1}=0$.
\item $\hat{p}^{b}_{2k,2l}=\hat{p}^{b}_{2k+1,2l}=\hat{p}^{b}_{2k+1,2l-1}=p_{k,l}$ et $\hat{p}^{b}_{2k+2,2l-1}=0$.
\item $\hat{p}^{c}_{2k,2l}=\hat{p}^{c}_{2k,2l+1}=\hat{p}^{c}_{2k-1,2l+1}=p_{k,l}$
et $\hat{p}^{c}_{2k-1,2l+2}=0$.
\end{cond}
On vérifie que, par construction, on a
$\hat{p}^a=T\hat{p}^b=S\hat{p}^c$ et que, pour tout $s$ dans
$\mathfrak S$, pour tout $d$ dans $\mathcal T_1$, on a
$\hat{p}^{sd}=s\left(\hat{p}^{d}\right)$.

\begin{Lem} \label{plongetriangle1}
Soit $d$ dans $\mathcal T_1$. L'application $p\mapsto \hat{p}^{d}$
induit un homéomor\-phisme de $Y$ dans $Y^{d}$. Réciproquement, un
point $p$ de $Y-\{Tu,Su\}$ appartient à $\hat{Y}$ si et seulement
si, pour tout $q$ dans $Y_p$, $q$ appartient à un triangle contenu
dans $Y_p$. Il existe alors un unique $d$ dans $\mathcal T_1$ et un
unique point $r$ de $Y$ tel que $p=\hat{r}^{d}$ et le triangle
contenant $p$ est
$\{\hat{r}^{a},\hat{r}^{b},\hat{r}^{c}\}$.\end{Lem}

La démonstration utilise le

\begin{Lem}\label{brechetvoisin}
Soit $p$ un point de $Y-\{Tu,Su\}$ tel que chaque point de $Y_p$
soit contenu dans un triangle. Alors, le triangle contenant $p$ est
soit $\{p,T^{-1}p,S^{-1}p\}$, soit $\{p,Tp,TS^{-1}p\}$, soit
$\{p,Sp,T^{-1}Sp\}$. S'il est de la forme $\{p,T^{-1}p,S^{-1}p\}$,
le troisième voisin $q$ de $p$ est soit $Tp$, soit $Sp$. Enfin, si
$q=Tp$, le triangle contenant $q$ est $\{q,Tq,TS^{-1}q\}$ et si
$q=Sp$, le triangle contenant $q$ est $\{q,Sq,T^{-1}Sq\}$.
\end{Lem}

\begin{demo} Soit $p$ comme dans l'énoncé. Quitte à faire agir
$\mathfrak S$, on peut supposer que le triangle contenant $p$
contient le point $T^{-1}p$. Alors, par dé\-finition, les seuls
voisins communs possibles pour $p$ et $T^{-1}p$ sont $T^{-1}Sp$ et
$S^{-1}p$. Or, comme $T^{-1}p$ appartient à $Y$, on a $p_{-1,0}=1$,
donc $p_{-1,1}=p_{0,0}+p_{-1,0}=0$ et $T^{-1}Sp\notin Y$. Par
conséquent, le triangle contenant $p$ est $\{p,T^{-1}p,S^{-1}p\}$.
Les autres cas s'en déduisent en faisant agir $\mathfrak S$.

Dans le cas où le triangle contenant $p$ est
$\{p,T^{-1}p,S^{-1}p\}$, le troisième voisin $q$ de $p$ est, par
construction, nécessairement dans $\{Tp,Sp\}$. Supposons, quitte à
encore faire agir $\mathfrak S$, qu'on a $q=Tp$. Alors, on a
$p_{0,1}=0=p_{1,-1}$ et, donc, $T^{-1}Sq$ et $S^{-1}q$
n'appartiennent pas à $Y$. Le triangle contenant $q$ est donc
$\{q,Tq,TS^{-1}q\}$. L'autre cas s'en déduit par symétrie.
\end{demo}

\begin{demo}[Démonstration du lemme \ref{plongetriangle1}]
On vérifie aisément que, pour $d$ dans $\mathcal T_1$, le point
$\hat{p}^{d}$ appartient à $Y^{d}$ et que l'application ainsi
définie induit un homéo\-morphisme de $Y$ dans $Y^{d}$.

Réciproquement, soit $p$ un point de $Y-\{Tu,Su\}$, tel que tout
élément de $Y_p$ soit contenu dans un triangle de $Y_p$. Alors, par
définition et d'après le lemme \ref{grapheinduit}, tout point de
$\tilde{Y}_p$ est contenu dans un triangle de $\tilde{Y_p}$. Soit
$(k,l)$, un point de $\tilde{Y}_p$. D'après le lemme
\ref{brechetvoisin}, le triangle contenant $(k,l)$ est de la forme
$\{(k,l),(k-1,l),(k,l-1)\}$, $\{(k,l),(k+1,l),(k+1,l-1)\}$ ou
$\{(k,l),(k,l+1),(k-1,l+1)\}$. On note
$(\varepsilon(k,l),\eta(k,l))$ l'unique élément de $(\mathbb
Z/2\mathbb Z)^2$ qui n'est pas congru aux éléments de ce triangle.

Montrons que, pour tous $(i,j)$ et $(k,l)$ dans $\tilde{Y}_p$ avec
$(i,j)\sim (k,l)$, on a
$(\varepsilon(i,j),\eta(i,j))=(\varepsilon(k,l),\eta(k,l))$. Si
$(i,j)$ et $(k,l)$ appartiennent au même triangle, c'est immédiat.
Sinon, quitte à faire agir $\mathfrak S$, d'après le lemme
\ref{brechetvoisin}, on peut supposer que le triangle contenant
$(k,l)$ est $\{(k,l),(k-1,l),(k,l-1)\}$, que $(i,j)=(k+1,l)$ et,
donc, que le triangle contenant $(i,j)$ est
$\{(k+1,l),(k+2,l),(k+2,l-1)\}$. Alors, par définition, on a
$(\varepsilon(i,j),\eta(i,j))=(\varepsilon(k,l),\eta(k,l))$.

Comme, d'après le lemme \ref{grapheinduit}, le graphe $\tilde{Y}_p$
est connexe, la fonction $(\varepsilon,\eta)$ y est constante. Par
définition, pour tous entiers $k$ et $l$, on a
$p_{2k+\varepsilon,2l+\eta}=0$, donc $p$ appartient à
$Y^{(\varepsilon,\eta)}$. La propriété sur les triangles découle
immédiatement de la définition des objets.
\end{demo}

Soit $p$ un point de $\hat{Y}$. On note $\bar{\Pi}p$ et
$\theta_1(p)$ les uniques éléments de $Y$ et de $\mathcal T_1$ pour
lesquels on a $\widehat{\bar{\Pi}p}^{\theta_1(p)}=p$. Par
construction, on a $\theta_1(p)=a$ (resp. $b$, resp. $c$) si et
seulement si le triangle contenant $p$ est $\{p,T^{-1}p,S^{-1}p\}$
(resp. $\{p,Tp,TS^{-1}p\}$, resp. $\{p,Sp,T^{-1}Sp\}$). Les
applications $\bar{\Pi}$ et $\theta_1$ sont $\mathfrak
S$-équivariantes. L'application $\bar{\Pi}$ est continue et
$\theta_1$ est localement constante. Pour tout $p$ dans $\hat{Y}$,
$\theta_1$ induit une bijection du $1$-triangle contenant $p$ dans
$\mathcal T_1$.

\begin{Lem} \label{plongetriangle2}
Soit $p$ dans $Y-\{Tu,Su\}$. Il existe un unique isomorphisme de
graphes $\sigma:\hat{Y_p}\rightarrow\bar{\Pi}^{-1}(Y_p)$ tel que
$\bar{\Pi}\sigma=\sigma\Pi$.\end{Lem}

\begin{demo} D'après le lemme \ref{automorphismes}, un tel isomorphisme est nécessairement
unique. Montrons son existence.
Soit $q$ le voisin de $p$ appartenant à $\{Tp,Sp\}$, $r$ son voisin
dans $\{T^{-1}p,T^{-1}Sp\}$ et $s$ son voisin dans
$\{S^{-1}p,TS^{-1}p\}$. On pose $\sigma(p,q)=\hat{p}^{a}$,
$\sigma(p,r)=\hat{p}^{b}$ et $\sigma(p,s)=\hat{p}^{c}$. Alors, les
trois points $\sigma(p,q)$, $\sigma(p,r)$ et $\sigma(p,s)$ sont
voisins dans $Y_p$. Vérifions, par exemple, que $\sigma(p,q)$ est
voisin de $\sigma(q,p)$. Quitte à faire agir $\mathfrak S$, on peut
supposer qu'on a $q=Tp$. Alors, on a $p=T^{-1} q$ et, donc,
$\sigma(q,p)=\hat{q}^{b}$. Par construction, on a alors
$T\hat{p}^a=\hat{q}^b$, donc $\sigma(q,p)=T\sigma(p,q)$, ce qu'il
fallait démontrer.
\end{demo}

Dorénavant, pour tout $p$ dans $\hat{Y}$, on identifie les graphes
$\hat{Y_p}$ et $\bar{\Pi}^{-1}(Y_p)$.

Nous allons à présent construire un élément $p$ de $Y$ pour lequel
le graphe $Y_p$ est isomorphe au graphe de Pascal. Posons, pour tous
$k,l\geq 0$, $p_{-k,-l}=p_{k+l+1,-l}={\rm C}^{k}_{l+k}$ dans
$\mathbb Z/2\mathbb Z$ et, pour tous $k,l$ dans $\mathbb Z$ avec
soit $l>0$, soit $k\geq 1$ et $k+l\geq 0$, $p_{k,l}=0$. On vérifie
aisément que $p$ appartient à $X$, et donc à $Y$ puisque
$p_{0,0}=1$. Nous avons la

\begin{Prop} \label{plongePascal}
Le point $p$ appartient à $\hat{Y}$ et on a $\bar{\Pi}p=p$ et
$\theta_1(p)=a$. Il existe un isomorphisme du graphe de Pascal
$\Gamma$ dans $Y_p$ envoyant $p_0$ sur $p$ et  $p_0^\vee$ sur
$Tp$.\end{Prop}

Cette représentation plane du graphe de pascal apparait dans la
figure \ref{graphep}. La démonstration utilise le

\begin{Lem} \label{binome}
Soient $0\leq k\leq n$ des entiers. Alors, les entiers ${\rm C}^{k}_{n}$,
${\rm C}^{2k}_{2n}$,  ${\rm C}^{2k}_{2n+1}$  et ${\rm C}^{2k+1}_{2n+1}$ sont congrus modulo $2$.\end{Lem}

\begin{demo} Soient $A$ et $B$ des indéterminées.
Dans l'anneau $\mathbb Z/2\mathbb Z[A,B]$ de caractéristique $2$, on
a $(A+B)^n=\sum_{k=0}^n{\rm C}^{k}_{n}A^kB^{n-k}$ et, donc,
$(A+B)^{2n}=\sum_{k=0}^n{\rm C}^{k}_{n}A^{2k}B^{2n-2k}$. Par
conséquent, par unicité, pour tout $0\leq k\leq n$, on a, dans
$\mathbb Z/2\mathbb Z$, ${\rm C}^{2k}_{2n}={\rm C}^{k}_{n}$ et ${\rm
C}^{2k-1}_{2n}={\rm C}^{2k+1}_{2n}=0$. Par l'identité classique, on
a alors ${\rm C}^{2k}_{2n+1}={\rm C}^{2k-1}_{2n}+{\rm
C}^{2k}_{2n}={\rm C}^{2k}_{2n}$ et ${\rm C}^{2k+1}_{2n+1}={\rm
C}^{2k}_{2n}+{\rm C}^{2k+1}_{2n}={\rm C}^{2k}_{2n}$.\end{demo}

\begin{demo}[Démonstration de la proposition \ref{plongePascal}]
En utilisant le lemme \ref{binome}, on vérifie qu'on a
$p=\hat{p}^{a}$. Par conséquent, $p$ appartient à $\hat{Y}$,
$\bar{\Pi}p=p$ et $\theta_1(p)=a$. Par récurence, en utilisant le
lemme \ref{plongetriangle2}, on en déduit que, pour tout entier $n$,
$p$ est le sommet d'un $n$-triangle contenu dans $Y_p$. Donc $p$ est
le sommet d'un triangle infini contenu dans $Y_p$. De même, on a
$\bar{\Pi}(Tp)=Tp$, $\theta_1(Tp)=b$ et $Tp$ est le sommet d'un
triangle infini contenu dans $Y_p$. D'après le lemme
\ref{grapheinduit}, le graphe $Y_p$ est connexe et, donc, il est la
réunion de ces deux triangles infinis. L'existence de l'isomorphisme
demandé en découle.\end{demo}

Dorénavant, on identifie $p$ à $p_0$, $Tp$ à $p_0^\vee$ et $\Gamma$
à $Y_p$. On note $\bar{\Gamma}$ l'adhérence de $\Gamma$ dans $Y$ et,
pour tout $p$ dans $\bar{\Gamma}$, on pose $\Gamma_p=Y_p$. On a
$\bar{\Pi}\bar{\Gamma}=\bar{\Gamma}$.

Nous allons décrire plus en détail l'ensemble $\bar{\Gamma}$. Pour
cela, introduisons une partition de $\hat{Y}$ en six sous-ensembles,
qui raffine la partition $\hat{Y}=Y^a\cup Y^b\cup Y^c$. Soit $p$ un
point de $\hat{Y}$. Alors, d'après le lemme \ref{brechetvoisin},
l'ensemble des voisins de $p$ est soit $\{Tp,T^{-1}p,S^{-1}p\}$,
soit $\{Sp,T^{-1}p,S^{-1}p\}$, soit $\{T^{-1}p,Tp,TS^{-1}p\}$, soit
$\{T^{-1}Sp,Tp,TS^{-1}p\}$, soit $\{S^{-1}p,Sp,T^{-1}Sp\}$ ou
$\{TS^{-1}p,Sp,T^{-1}Sp\}$. Appelons bréchet de $p$ l'ensemble des
$q$ dans $\hat{Y}$ pour lesquels on a $\{(k,l)\in \mathbb Z^2|T^kS^l
p\sim p\}=\{(k,l)\in \mathbb Z^2|T^kS^l q\sim q\}$. Les bréchets
constituent six sous-ensembles fermés de $\hat{Y}$ sur lesquels le
groupe $\mathfrak S$ agit simplement transitivement. On notera $B_0$
le bréchet de $p_0$, $i$ l'élément
$\begin{pmatrix}-1&-1\\0&1\end{pmatrix}$ de $\mathfrak S$ et $r$
l'élément $\begin{pmatrix}0&1\\-1&-1\end{pmatrix}$. L'élément $i$
s'identifie à la transposition $(ab)$ de $\{a,b,c\}$ et $r$
s'identifie au cycle $(cba)$.

Pour tout entier $n$, posons $\hat{Y}^{(n)}=\bar{\Pi}^{-n}Y$. Alors,
par une récurrence immédiate, d'après les lemmes
\ref{plongetriangle1} et \ref{plongetriangle2}, pour tout entier
$n$, $\hat{Y}^{(n)}$ est l'ensemble des éléments $p$ de
$Y-\{Tu,Su\}$ pour lesquels tout point de $Y_p$ appartient à un
$n$-triangle dans $Y_p$. On a donc
$\Gamma\subset\bigcap_{n\in\mathbb N}\hat{Y}^{(n)}$.

\begin{Lem} \label{brechet2}
Soient $n$ un entier et $p$ et $q$ dans $\hat{Y}^{(n+1)}$ tels que,
pour tout $0\leq m\leq n$, $\bar{\Pi}^mp$ et $\bar{\Pi}^mq$
appartiennent au même bréchet. Alors, pour tous $k$ et $l$ dans
$\mathbb Z$ avec $k\geq -2^n$,  $l\geq -2^n$ et $k+l\leq 2^n$, on a
$p_{k,l}=q_{k,l}$.
\end{Lem}

\begin{demo} Montrons ce résultat par récurrence sur $n$. Pour
$n=0$, supposons, quitte à faire agir $\mathfrak S$, qu'on a $p,q\in
B_0$. Alors, on a $p_{-1,0}=p_{0,0}=p_{1,0}=p_{0,-1}=1$, si bien que
$p_{-1,-1}=p_{-1,0}+p_{0,-1}=0$ et, de même,
$p_{1,-1}=p_{-1,1}=p_{0,1}=p_{-1,2}=0$ et $p_{2,-1}=1$. Comme ceci
est aussi vrai pour $q$, le lemme est vrai pour $n=0$.

Supposons donc $n\geq 1$ et le lemme démontré pour $n-1$.
Donnons-nous $p$ et $q$ comme dans l'énoncé. Alors, comme $p$ et $q$
sont dans le même bréchet, on a $\theta_1(p)=\theta_1(q)$. Quitte à
faire agir $\mathfrak S$, on peut supposer qu'on a $\theta_1(p)=a$,
si bien que $p=\widehat{\bar{\Pi}p}^a$ et que
$q=\widehat{\bar{\Pi}q}^a$. Le résultat découle alors de la
récurrence et de la définition de l'application
$r\mapsto\hat{r}^a$.\end{demo}

\begin{Lem} \label{brechet1} Soit $p$ dans $\hat{Y}$ tel que
le bréchet de $p$ soit $B_0$. Alors, le bréchet de $\hat{p}^a$ est
$B_0$, celui de $\hat{p}^b$ est $iB_0$ et celui de $\hat{p}^c$ est
$rB_0$ et $\hat{p}^a$ est le sommet d'un $2$-triangle dans
$\bar{\Gamma}$. Si $q$ et $r$ sont deux points de $\hat{Y}^{(2)}$
tels que $\bar{\Pi}q$ et $\bar{\Pi}r$ appartiennent au même bréchet,
il existe $r'$ dans le $1$-triangle contenant $r$ dans $Y_r$ tel que
$q$ et $r'$ appartiennent au même bréchet.\end{Lem}

\begin{demo} Le premier point découle immédiatement de la
construction des objets en question. Donnons-nous donc $q$ et $r$
comme dans l'énoncé. Quitte à faire agir $\mathfrak S$, on peut
supposer que le bréchet de $\bar{\Pi}q$ et de $\bar{\Pi}r$ est
$B_0$. La première partie du lemme implique alors clairement
l'énoncé.
\end{demo}

Soit $\Sigma$ l'ensemble des suites $(s_n)_{n\in\mathbb N}$
d'éléments de $\mathfrak S$ telles que, pour tout entier $n$, on ait
$s_n\in\{s_{n+1},s_{n+1}i,s_{n+1}r\}$. On munit $\Sigma$ de la
topologie induite par la topologie produit et on note
$\sigma:\Sigma\rightarrow\Sigma$ l'application de décalage. On fait
agir $\mathfrak S$ sur $\Sigma$ par multiplication à gauche sur tous
les facteurs. L'ensemble $\bar{\Gamma}$ est décrit par la

\begin{Prop} \label{codage}
On a $\bar{\Gamma}=\bigcap_{n\in\mathbb N}\hat{Y}^{(n)}$. Pour tout
$p$ dans $\bar{\Gamma}$, pour tout entier $n$, soit $s_n(p)$
l'unique élément de $\mathfrak S$ tel que le bréchet de
$\bar{\Pi}^np$ soit $s_n(p)B_0$. L'application $s$ ainsi définie
induit un homéomorphisme $\mathfrak S$-équi\-va\-riant de
$\bar{\Gamma}$ dans $\Sigma$ et on a $\sigma s=s\bar{\Pi}$. L'image
du point $p_0$ par $s$ est la suite constante de valeur $e$ et les
points fixes de $\bar{\Pi}$ dans $\bar{\Gamma}$ sont exactement les
six images de $p_0$ par l'action du groupe $\mathfrak S$. Enfin,
pour tout $p$ dans $\bar{\Gamma}$, l'ensemble $\Gamma_p$ est dense
dans $\bar{\Gamma}$.
\end{Prop}

\begin{demo} Comme l'ensemble $\Gamma$ est inclus dans $\bigcap_{n\in\mathbb N}\hat{Y}^{(n)}$,
l'ensem\-ble $\bar{\Gamma}$ l'est aussi.
Réciproquement, remarquons que chacun des six points $T^{-1}p_0$,
$T^{-2}p_0$, $T^{-2}S^{-1}p_0$, $T^{-1}S^{-2}p_0$, $S^{-2}p_0$ et
$S^{-1}p_0$ appartient à un bréchet différent. Par conséquent, si
$p$ est un point de $\bigcap_{n\in\mathbb N}\hat{Y}^{(n)}$, il
existe un point $q$ de $\Gamma$ tel que $p$ et $q$ appartiennent au
même bréchet. Par récurrence, en utilisant le lemme \ref{brechet1},
on en déduit que, pour tout entier $n$, il existe un point $q_n$ de
$\Gamma$ tel que, pour tout  $0\leq m\leq n$, $\bar{\Pi}^mp$ et
$\bar{\Pi}^mq_n$ appartiennent au même bréchet. D'après le lemme
\ref{brechet2}, on a alors $q_n\td{n}{\infty}p$ et $p$ appartient à
$\bar{\Gamma}$.

Soit $p$ dans $\bar{\Gamma}$. Comme on vient de le voir, le point
$p$ est complètement déterminé par la suite
$s(p)=(s_n(p))_{n\in\mathbb N}$. L'application $s$ est clairement
continue et $\mathfrak S$-équivariante et, par définition, on a
$s\bar{\Pi}=\sigma s$. Par ailleurs, d'après le lemme
\ref{plongetriangle1}, si $p$ est un point de $\bar{\Gamma}$, il
possède exactement trois an\-técédents par $\bar{\Pi}$ et, d'après
le lemme \ref{brechet1}, les bréchets de ces trois antécédents sont
$s_0(p) B_0$, $s_0(p)iB_0$ et $s_0(p) rB_0$. On en déduit que
l'application $s$ prend ses valeurs dans $\Sigma$ et qu'elle induit
un homéomorphisme de $\bar{\Gamma}$ dans $\Sigma$.

Par construction, on a $s_0(p_0)=e$ et, comme, d'après la
proposition \ref{plongePascal}, $\bar{\Pi}p_0=p_0$, pour tout entier
naturel $n$, $s_n(p_0)=e$. En particulier, les autres points fixes
de $\bar{\Pi}$ sont les images de $p_0$ par l'action de $\mathfrak
S$.

Enfin, si $p$ est un point de $\bar{\Gamma}$, d'après les lemmes
\ref{structuregrostriangles} et \ref{plongetriangle2}, pour tout
entier naturel $n$, le $n$-triangle contenant $p$ dans $\Gamma_p$
est l'ensemble $\bar{\Pi}^{-n}\left(\bar{\Pi}^np\right)$. Comme $r$
et $i$ engendrent le groupe $\mathfrak S$, on vérifie aisément que
le sous-décalage de type fini $(\Sigma,\sigma)$ est transitif, si
bien que, pour tout $t$ dans $\Sigma$, l'ensemble
$\bigcup_{n\in\mathbb N}\sigma^{-n}(\sigma^nt)$ est dense dans
$\Sigma$. Par conséquent, pour tout $p$ dans $\bar{\Gamma}$,
$\Gamma_p$ est dense dans $\bar{\Gamma}$.
\end{demo}

\section{Fonctions triangulaires et intégration sur $\bar{\Gamma}$}
\label{fonctiontriangulaire}

Dans cette section, nous étudions une classe particulière de
fonctions localement constantes sur $\bar{\Gamma}$. Nous utilisons
ces fonctions pour déterminer les propriétés d'une mesure de Radon
remarquable sur $\bar{\Gamma}$. Pour $p$ dans $\bar{\Gamma}$, on
note toujours, comme à la section \ref{compactifie},
$(s_n(p)B_0)_{n\in\mathbb N}$ la suite de bréchets associée.

Soient $n$ un entier $\geq 1$ et $\mathcal T$ un $n$-triangle dans
$\bar{\Gamma}$. Alors, d'après le lemme \ref{plongetriangle2},
l'ensemble $\bar{\Pi}^{n-1}\mathcal T$ est un $1$-triangle de
$\bar{\Gamma}$. Par conséquent, l'application
$\theta_1\circ\bar{\Pi}^{n-1}$ induit une bijection de l'ensemble
des sommets de $\mathcal T$ dans $\mathcal T_1=\{a,b,c\}$. On note
$a_n$ (resp. $b_n$, resp. $c_n$) l'ensemble des sommets $p$ de
$n$-triangles de $\bar{\Gamma}$ tels que
$\theta_1\left(\bar{\Pi}^{n-1}p\right)=a$ (resp. $b$, resp. $c$) et
$\theta_n$ l'application qui, à un sommet $p$ d'un $n$-triangle de
$\bar{\Gamma}$, associe l'élément de $\{a_n,b_n,c_n\}$ auquel il
appartient. Soit $\mathcal T_n$ le $n$-triangle $\mathcal
T_n(a_n,b_n,c_n)$. D'après le lemme \ref{grostriangles},
l'application $\theta_n$ se prolonge de manière unique en une
application $\bar{\Gamma}\rightarrow\mathcal T_n$, encore notée
$\theta_n$, qui, pour tout $n$-triangle $\mathcal T$ de
$\bar{\Gamma}$, induit un isomorphisme de graphes de $\mathcal T$
dans $\mathcal T_n$. Cette application est localement constante.
Pour $n=1$, cette définition est cohérente avec les notations de la
section \ref{compactifie}, en identifiant $a$ et $a_1$, $b$ et $b_1$
et $c$ et $c_1$. Par abus de langage, on considérera parfois que
$\mathcal T_0$ est un ensemble contenant un seul élément et que
$\theta_0$ est l'application constante
$\bar{\Gamma}\rightarrow\mathcal T_0$.

Pour tout $n\geq 1$, le groupe $\mathfrak S$ agit sur $\mathcal T_n$
en s'identifiant à $\mathfrak S(a_n,b_n,c_n)$ et l'application
$\theta_n$ est $\mathfrak S$-équivariante. On identifiera $\mathcal
T_{n+1}$ et $\hat{\mathcal T}_n$ à travers la bijection $\mathfrak
S$-équivariante de $\{a_n,b_n,c_n\}$ dans
$\{a_{n+1},b_{n+1},c_{n+1}\}$ qui envoie $a_n$ sur $a_{n+1}$, $b_n$
sur $b_{n+1}$ et $c_n$ sur $c_{n+1}$. En particulier, on notera
$\Pi:\mathcal T_{n+1}\rightarrow\mathcal T_n$ l'application de
contraction des triangles provenant de cette identification et
$\Pi^*$ et $\Pi$ les opérateurs associés $\ell^2(\mathcal
T_{n})\rightarrow \ell^2(\mathcal T_{n+1})$ et $\ell^2(\mathcal
T_{n+1})\rightarrow \ell^2(\mathcal T_{n})$.

Notons, pour tout $n\geq 2$, $a_nb_n$, $a_nc_n$, $b_na_n$, $b_nc_n$,
$c_na_n$ et $c_nb_n$ les points de $\mathcal T_n$ définis par le
corollaire \ref{decoupetriangles}. Les principales propriétés des
applications $\theta_n$, $n\geq 1$, que nous utiliserons par la
suite sont décrites par le

\begin{Lem} \label{facteurtriangle}
Soit $n$ un entier $\geq 1$. On a
$\Pi\theta_{n+1}=\theta_n\bar{\Pi}$. Si $p$ et $q$ sont des points
de $\bar{\Gamma}$ tels que $\theta_{n+1}(p)=\theta_{n+1}(q)$, on a
$\theta_n(p)=\theta_n(q)$. En particulier, on a $\theta_{n}(p)=a_n$
si et seulement si $\theta_{n+1}(p)$ est $a_{n+1}$, $b_{n+1}a_{n+1}$
ou $c_{n+1}a_{n+1}$. Soient $p$ et $q$ dans $\bar{\Gamma}$, tels que
$\theta_n(p)=\theta_n(q)$. Pour tout $0\leq m\leq n$, si
$\theta_n(p)$ n'est pas contenu dans le $m$-triangle issu d'un des
sommets de $\mathcal T_n$, on a $s_m(p)=s_m(q)$.
\end{Lem}

\begin{demo} Soit $\mathcal T$ un $(n+1)$-triangle de $\bar{\Gamma}$.
L'application $\theta_{n+1}$ induit un isomorphisme de $\mathcal T$
dans $\mathcal T_{n+1}$ et l'application $\theta_n$ induit un
isomorphisme du $n$-triangle $\bar{\Pi}\mathcal T$ dans $\mathcal
T_n$. Comme, par définition, les applications $\theta_n\bar{\Pi}$ et
$\Pi\theta_{n+1}$ coïncident sur l'ensemble $\partial\mathcal T_n$,
on a, d'après le lemme \ref{automorphismes},
$\Pi\theta_{n+1}=\theta_n\bar{\Pi}$.

Soient $p$, $q$ et $r$ les sommets de $\mathcal T$, de sorte que
$\theta_{n+1}(p)=a_{n+1}$, $\theta_{n+1}(q)=b_{n+1}$ et
$\theta_{n+1}(r)=c_{n+1}$. Par définition, on a $\theta_n(p)=a_n$.
Montrons qu'on a  $\theta_n(qp)=\theta_n(rp)=a_n$. Cela revient à
montrer qu'on a $\theta_1\left(\bar{\Pi}^{n-1}qp\right)=
\theta_1\left(\bar{\Pi}^{n-1}rp\right)=a_1$. Or, par le raisonnement
ci-dessus, on a
$\theta_2\left(\bar{\Pi}^{n-1}qp\right)=\Pi^{n-1}\theta_{n+1}(qp)=b_2a_2$
et
$\theta_2\left(\bar{\Pi}^{n-1}rp\right)=\Pi^{n-1}\theta_{n+1}(rp)=c_2a_2$,
si bien qu'on est ramené au cas où $n=1$. Alors, avec les notations
de la section \ref{compactifie}, si $s=\bar{\Pi}^2p$, on vérifie
qu'on a $qp=\widehat{\hat{s}^b}^a$ et $rp=\widehat{\hat{s}^c}^a$,
d'où le résultat.

En particulier, si $\mathcal S$ est un $n$-triangle de
$\bar{\Gamma}$, la restriction de $\theta_n$ à $\partial\mathcal S$
est complètement déterminée par la restriction de $\theta_{n+1}$ à
$\partial S$. Par définition, les valeurs de $\theta_{n+1}$ sont
donc déterminées par celles de $\theta_n$.

Soient enfin $p$ et $q$ tels que $\theta_n(p)=\theta_n(q)$ et
montrons par récurrence sur $n\geq 1$ la propriété du lemme. Si
$n=1$, cette propriété est vide. Supposons $n\geq 2$ et la propriété
établie pour $n-1$. Alors, on a
$\theta_{n-1}\left(\bar{\Pi}p\right)=\Pi\theta_n(p)=\Pi\theta_n(q)=\theta_{n-1}\left(\bar{\Pi}q\right)$
et, pour tout entier $m$ avec $1\leq m\leq n$, si $\theta_n(p)$
n'appartient pas au $m$-triangle issu d'un des sommets de $\mathcal
T_n$, $\theta_{n-1}\left(\bar{\Pi}p\right)$ n'appartient pas au
$(m-1)$-triangle issu d'un des sommets de $\mathcal T_{n-1}$ et,
donc, par récurrence,
$s_m(p)=s_{m-1}\left(\bar{\Pi}p\right)=s_{m-1}\left(\bar{\Pi}q\right)=s_m(q)$.
Il nous reste à traiter le cas où $m=0$. Supposons donc que
$\theta_n(p)$ n'est pas un sommet de $\mathcal T_n$ et montrons que
$s_0(p)=s_0(q)$. Remarquons, que, par la première partie de la
preuve, on a $\theta_1(p)=\theta_1(q)$. Quitte à faire agir
$\mathfrak S$, supposons que $\theta_1(p)=a_1$. Alors, soient
$\mathcal T$ le $n$-triangle de $\bar{\Gamma}$ contenant $p$ et $p'$
le voisin de $p$ qui n'appartient pas au $1$-triangle contenant $p$.
Comme $p$ n'est pas un sommet de $\mathcal T$, $p'$ appartient à
$\mathcal T$ et, d'après le lemme \ref{brechetvoisin}, $s_0(p)$ est
$B_0$ ou $riB_0$, suivant que $\theta_1(p')$ est $b_1$ ou $c_1$. Or,
comme $\theta_n$ induit un isomorphisme de $\mathcal T$ dans
$\mathcal T_n$, $\theta_n(p')$ ne dépend que de $\theta_n(p)$ et,
donc, toujours par la première partie du lemme, la valeur de
$\theta_1$ en $p'$ est complètement déterminée par la valeur de
$\theta_n$ en $p$. Par conséquent, la valeur de $s_0$ en $p$ est
déterminée par la valeur de $\theta_n$ en $p$, ce qu'il fallait
démontrer.
\end{demo}

Pour tout entier $n\geq 1$, d'après la proposition
\ref{plongePascal}, on a $\theta_n(p_0)=a_n=\theta_n(rip_0)$, si
bien que le codage de $\bar{\Gamma}$ par les applications
$\theta_n$, $n\geq 1$, possède des ambiguïtés. Elles sont décrites
par le

\begin{Cor}\label{trianglesepare}
Soient $p$ et $q$ dans $\bar{\Gamma}$ tels que, pour tout entier
$n$, on ait $\theta_n(p)=\theta_n(q)$. Alors, si $p\neq q$, il
existe $s$ dans $\mathfrak S$ tel que $p$ appartienne au triangle
infini issu de $sp_0$ dans $s\Gamma$ et que $q$ appartienne au
triangle infini issu de $srip_0$ dans $sri\Gamma$.
\end{Cor}

\begin{demo} Supposons qu'on a $p\neq q$. Alors, d'après la
proposition \ref{codage}, il existe un entier $m$ tel que
$s_m(p)\neq s_m(q)$. D'après le lemme \ref{facteurtriangle}, pour
tout $n\geq m$, le point $\theta_n(p)=\theta_n(q)$ appartient au
$m$-triangle issu de l'origine de $\mathcal T_n$. Posons
$p'=\bar{\Pi}^mp$ et $q'=\bar{\Pi}^mq$. D'après le lemme
\ref{facteurtriangle}, pour tout entier $n$, on a
$\theta_n(p')=\Pi^m\theta^{n+m}(p)=\Pi^m\theta^{n+m}(p)=\theta_n(q')$
et ce point est un des sommets de $\mathcal T_n$. Comme, pour tout
entier $n\geq 1$, on a $\theta_n(a_{n+1})=a_n$,
$\theta_n(b_{n+1})=b_n$ et $\theta_n(c_{n+1})=c_n$, on peut
supposer, quitte à faire agir $\mathfrak S$, qu'on a, pour tout
entier $n\geq 1$, $\theta_n(p')=\theta_n(q')=a_n$. Comme
$\theta_1(p')=a_1$, le bréchet de $p'$ est $B_0$ ou $riB_0$. Quitte
à faire agir $\mathfrak S$, supposons que ce bréchet soit $B_0$.
Alors, d'après le lemme \ref{brechet1}, le bréchet de $\bar{\Pi}p'$
est $B_0$, $iB_0$ ou $r^{-1} B_0$. Comme
$\theta_1\left(\bar{\Pi}p'\right)=\Pi\theta_1(p)=a_1$, le bréchet de
$\bar{\Pi}p'$ est $B_0$ et, par récurrence, pour tout entier $n$, on
a $s_n(p')=e$, si bien que, d'après la proposition \ref{codage},
$p'=p_0$ et que $p$ appartient au triangle infini $\mathcal
T_\infty(p_0)$ issu de $p_0$ dans $\Gamma$. De même, on a $q'=p_0$
ou $q'=ri p_0$ et $q$ appartient au triangle infini issu de $p_0$
dans $\Gamma$ ou au triangle infini issu de $rip_0$ dans $ri\Gamma$.
Pour tout entier $n$, l'application $\theta_n$ induit une bijection
du $n$-triangle de $\Gamma$ issu de $p_0$ dans $\mathcal T_n$. Par
conséquent, si $p''$ est un point de $\mathcal T_\infty(p_0)$ tel
que, pour tout entier $n$, on ait $\theta_n(p'')=\theta_n(p)$, on a
$p''=p$. Comme on a supposé $p\neq q$, $q$ appartient au triangle
infini issu de $rip_0$ dans $ri\Gamma$, ce qu'il fallait démontrer.
\end{demo}

Soit $n$ un entier. Nous dirons qu'une fonction
$\varphi:\bar{\Gamma}\rightarrow\mathbb C$ est $n$-triangulaire si
elle s'écrit $\psi\circ\theta_n$, pour une certaine fonction $\psi$
sur $\mathcal T_n$. Quand il n'y aura pas d'ambiguïté, pour alléger
les notations, on identifiera $\varphi$ et $\psi$. D'après le lemme
\ref{facteurtriangle}, une fonction $n$-triangulaire est
$(n+1)$-triangulaire. En particulier, les fonctions triangulaires
constituent une sous-algèbre de l'algèbre des fonctions localement
constantes sur $\bar{\Gamma}$. Comme, pour toute fonction
triangulaire $\varphi$, on a $\varphi(p_0)=\varphi(rip_0)$, cette
algèbre n'est pas dense dans $\mathcal C^0\left(\bar{\Gamma}\right)$
pour la topologie de la convergence uniforme.

Dorénavant, on notera $\mu$ la probabilité borélienne sur
$\bar{\Gamma}$ dont l'image par l'application de codage de la
proposition \ref{codage} est la mesure d'entropie maximale pour
$\sigma$ sur $\Sigma$. En d'autres termes, $\mu$ est l'unique mesure
telle que, pour toute suite $t_0,\ldots,t_n$ d'éléments de
$\mathfrak S$, si, pour tout $0\leq m\leq n-1$, on a
$t_m\in\{t_{m+1},t_{m+1}i,t_{m+1}r\}$, alors $\mu\left(t_0
B_0\cap\bar{\Pi}^{-1}t_1 B_0\cap \ldots\cap\bar{\Pi}^{-n}t_n
B_0\right)=\frac{1}{6.3^{n}}$. Par définition, la mesure $\mu$ est
$\bar{\Pi}$-invariante et $\mathfrak S$-invariante.

Pour toute fonction borélienne $\varphi$ sur $\bar{\Gamma}$, on note
$\bar{\Pi}^*\varphi=\varphi\circ\bar{\Pi}$. Pour tout $1\leq
p\leq\infty$, l'opérateur $\bar{\Pi}^*$ préserve la norme dans ${\rm
L}^p\left(\bar{\Gamma},\mu\right)$. On note $\bar{\Pi}$ son adjoint,
c'est-à-dire que, pour toute fonction borélienne $\varphi$ sur
$\bar{\Gamma}$, pour tout $p$ dans $\bar{\Gamma}$, on a
$\bar{\Pi}\varphi(p)=\frac{1}{3}\sum_{\bar{\Pi}(q)=p}\varphi(q)$. On
a $\bar{\Pi}1=1$ et, pour tout $1\leq p\leq\infty$, l'opérateur
$\bar{\Pi}$ est de norme $1$ dans ${\rm
L}^p\left(\bar{\Gamma},\mu\right)$. Enfin, on a
$\bar{\Pi}\bar{\Pi}^*=1$.

L'intégrale des fonctions triangulaires pour la mesure $\mu$ se
calcule naturellement :

\begin{Lem} \label{fonctiontriangulaire2}
Soient $n$ un entier naturel et $\varphi$ une fonction
$n$-tri\-angu\-laire. On a
$\int_{\bar{\Gamma}}\varphi\de\mu=\frac{1}{3^n}\sum_{p\in\mathcal
T_n}\varphi(p)$.
\end{Lem}

En d'autres termes, la mesure image de $\mu$ par $\theta_n$ est la
mesure de comptage normalisée de $\mathcal T_n$.

\begin{demo}
Démontrons le résultat par récurrence sur $n$. Si $n=0$, $\varphi$
est constante et le lemme est évident. Si $n\geq 1$, comme, d'après
le lemme \ref{facteurtriangle}, on a
$\Pi\theta_{n}=\theta_{n-1}\bar{\Pi}$, la fonction
$\bar{\Pi}\varphi$ est $(n-1)$-triangulaire et on a, par récurrence,
$\int_{\bar{\Gamma}}\varphi\de\mu=\int_{\bar{\Gamma}}\bar{\Pi}\varphi\de\mu
=\frac{1}{3^{n-1}}\sum_{p\in\mathcal
T_{n-1}}\frac{1}{3}\sum_{\Pi(q)=p}\varphi(q)
=\frac{1}{3^n}\sum_{p\in\mathcal T_n}\varphi(p)$, d'où le
résultat.\end{demo}

Dorénavant, pour tout entier $n$, on identifiera $\theta_n$ et la
partition associée de l'espace mesuré
$\left(\bar{\Gamma},\mu\right)$. D'après le lemme
\ref{facteurtriangle}, cette suite de partitions est croissante.
Comme $\mu$ n'a pas d'atomes, on a $\mu\left(\bigcup_{s\in\mathfrak
S}s\Gamma\right)=0$ et, d'après le lemme \ref{trianglesepare}, pour
tous $p$ et $q$ dans l'ensemble de mesure totale
$\bar{\Gamma}-\bigcup_{s\in\mathfrak S}s\Gamma$, si, pour tout
entier $n$, $\theta_n(p)=\theta_n(q)$, on a $p=q$. Pour tout
$\varphi$ dans ${\rm L}^1\left(\bar{\Gamma},\mu\right)$, pour tout
entier $n$, on note $\mathbb E(\varphi|\theta_n)$ l'espérance
conditionnelle de $\varphi$ sachant $\theta_n$, c'est-à-dire que,
pour tout $p$ dans $\mathcal T_n$, on a $\mathbb
E(\varphi|\theta_n)(p)=
\frac{1}{\mu(\theta_n^{-1}(p))}\int_{\theta_n^{-1}(p)}\varphi\de\mu$.

\begin{Lem} \label{fonctiontriangulaire3}
Pour tout $1\leq p< \infty$, pour tout $\varphi$ dans ${\rm
L}^p\left(\bar{\Gamma},\mu\right)$, on a $\mathbb
E(\varphi|\theta_n)\td{n}{\infty}\varphi$ dans ${\rm
L}^p\left(\bar{\Gamma},\mu\right)$. En particulier, l'espace des
fonctions triangulaires est dense dans ${\rm
L}^p\left(\bar{\Gamma},\mu\right)$. De même, l'espace des fonctions
triangulaires qui sont nulles aux sommets de leur triangle de
définition est dense dans ${\rm L}^p\left(\bar{\Gamma},\mu\right)$.
Enfin, pour tout entier $n$, pour tout $\varphi$ dans ${\rm
L}^1\left(\bar{\Gamma},\mu\right)$, on a $\mathbb
E\left(\bar{\Pi}^*\varphi|\theta_{n+1}\right)=\Pi^*\mathbb
E(\varphi|\theta_n)$, $\mathbb
E\left(\bar{\Pi}\varphi|\theta_{n}\right)=\frac{1}{3}\Pi\mathbb
E(\varphi|\theta_{n+1})$ et, pour $p$ dans $\mathcal T_n$,
$$\mathbb
E(\varphi|\theta_n)(p)=\frac{1}{3}\sum_{\substack{q\in\mathcal
T_{n+1}\\ \theta_n(q)=p}}\mathbb E(\varphi|\theta_{n+1})(q).$$
\end{Lem}

\begin{demo} La convergence dans ${\rm
L}^p\left(\bar{\Gamma},\mu\right)$ résulte de la discussion
précé\-dente et de propriétés générales des partitions des espaces
de probabilité. La densité des fonctions triangulaires nulles aux
sommets de leur triangle de définition en découle, vu que, d'après
le lemme \ref{fonctiontriangulaire3}, pour tout $n\geq 1$, la mesure
de l'ensemble des éléments de $\bar{\Gamma}$ qui sont sommet d'un
$n$-triangle est $3^{n-1}$. Enfin, les formules liant les espérances
conditionnelles sachant $\theta_{n+1}$ et $\theta_n$ découlent du
lemme \ref{facteurtriangle} et du fait que, d'après le lemme
\ref{fonctiontriangulaire2}, la mesure image de $\mu$ par $\theta_n$
est la mesure de comptage normalisée de $\mathcal T_n$.
\end{demo}

Décrivons enfin un homéomorphisme
$\bar{\Gamma}\rightarrow\bar{\Gamma}$ qui nous sera utile pour la
suite. Pour tout $p$ dans $\bar{\Gamma}$, on note $\alpha(p)$
l'unique voisin de $p$ dans $\bar{\Gamma}$ qui n'appartient pas au
triangle contenant $p$. L'application
$\alpha:\bar{\Gamma}\rightarrow\bar{\Gamma}$ est une involution sans
points fixes. D'après le corollaire \ref{sommettriangles}, pour tout
entier $n\geq 1$, $\alpha$ laisse stable l'ensemble des points de
$\bar{\Gamma}$ qui sont sommet d'un $n$-triangle. Pour $p$ dans
$\mathcal T_n-\partial\mathcal T_n$, on note encore $\alpha_n(p)$
l'unique voisin de $p$ dans $\mathcal T_n$ qui n'appartient pas au
triangle contenant $p$.

\begin{Lem}\label{involution}
Pour tout entier $n$, pour tout $p$ dans $\bar{\Gamma}$, si $p$
n'est pas le sommet d'un $n$-triangle de $\bar{\Gamma}$, on a
$\theta_n(\alpha(p))=\alpha_n(\theta_n(p))$. L'application $\alpha$
préserve la mesure $\mu$ et, pour tous $n\geq 1$,  $\varphi$ dans
${\rm L}^1\left(\bar{\Gamma},\mu\right)$ et $p$ dans $\mathcal
T_n-\partial\mathcal T_n$, on a $\mathbb
E\left(\varphi\circ\alpha|\theta_{n}\right)(p)=\mathbb
E(\varphi|\theta_n)(\alpha_n(p))$.
\end{Lem}

Comme les fonctions triangulaires ne sont pas denses dans $\mathcal
C^0\left(\bar{\Gamma}\right)$, pour vérifier que $\alpha$ préserve
la mesure $\mu$, nous aurons recours au

\begin{Lem}\label{compactatrou}
Soit $X$ un espace métrique compact et soit $A$ une sous-algèbre de
$\mathcal C^0(X)$ stable par conjugaison complexe et fermée pour la
topologie de la convergence uniforme. Soit $Y$ l'ensemble des
éléments $x$ de $X$ pour lesquels il existe $y\neq x$ dans $X$ tel
que, pour tout $\varphi$ dans $A$, $\varphi(y)=\varphi(x)$.
L'ensemble $Y$ est borélien et, si $\lambda$ est une mesure complexe
borélienne sur $X$ telle que $\lambda_{|Y}=0$ et que, pour tout
$\varphi$ dans $A$, $\int_X\varphi\de\lambda=0$, on a $\lambda=0$.
\end{Lem}

\begin{demo} Soient $S$ le spectre de la ${\rm C}^*$-algèbre commutative
$A$ et $\pi:X\rightarrow S$ la surjection continue duale de
l'injection naturelle de $A$ dans $\mathcal C^0(X)$. Par hypothèse,
la mesure complexe $\pi_*\lambda$ est nulle sur $S$. Notons $p$ la
projection sur la première composante $X\times X\rightarrow X$ et
posons $D=\{(x,x)|x\in X\}\subset X\times X$ et $E=\{(x,y)\in
X\times X|\pi(x)=\pi(y)\}$. On a $Y=p(E-D)$. Comme $E$ et $D$ sont
des fermés de l'espace compact métrisable $X$, l'ensemble $E-D$ est
une réunion dénombrable de compacts, donc $Y$ et $\pi(Y)$ le sont
aussi. En particulier, ces ensembles sont boréliens et $\pi$ induit
un isomorphisme borélien de $X-Y$ dans $S-\pi(Y)$. Par conséquent,
la restriction de $\lambda$ à $X-Y$ est nulle. Comme sa restriction
à $Y$ est nulle, on a $\lambda=0$.
\end{demo}

\begin{demo}[Démonstration du lemme \ref{involution}]
La première partie du lemme résulte de la définition de $\alpha$ et
du fait que $\theta_n$ induit un morphisme de graphes du
$n$-triangle contenant $p$ dans $\mathcal T_n$.

Soit $n\geq 1$. Comme, $\alpha$ permute les points de $\bar{\Gamma}$
qui sont sommet d'un $n$-triangle, si $\varphi$ est une fonction
$n$-triangulaire nulle aux sommets de $\mathcal T_n$, on a, d'après
le lemme \ref{fonctiontriangulaire2},
$\int_{\bar{\Gamma}}\varphi\circ\alpha\de\mu=\int_{\bar{\Gamma}}\varphi\de\mu$.
Comme, à nouveau d'après le lemme \ref{fonctiontriangulaire2}, pour
tout entier $n\geq 1$, la mesure de l'ensemble des points de
$\bar{\Gamma}$ qui sont sommet d'un $n$-triangle est $3^{n-1}$, on
en déduit que, pour toute fonction triangulaire $\varphi$, on a
$\int_{\bar{\Gamma}}\varphi\circ\alpha\de\mu=\int_{\bar{\Gamma}}\varphi\de\mu$.

Soit $A\subset\mathcal C^0\left(\bar{\Gamma}\right)$ l'adhérence de
l'algèbre des fonctions triangulaires pour la topologie de la
convergence uniforme. D'après le corollaire \ref{trianglesepare},
l'ensemble des éléments $p$ de $\bar{\Gamma}$ pour lesquels il
existe $q\neq p$ tel que, pour tout $\varphi$ dans $A$, on ait
$\varphi(p)=\varphi(q)$ est $\bigcup_{s\in\mathfrak S}s\Gamma$.
Comme $\mu$ n'a pas d'atomes, cet ensemble est de mesure nulle pour
$\mu$ et pour $\alpha_*\mu$. Pour tout $\varphi$ dans $A$, on a
$\int_{\bar{\Gamma}}\varphi\circ\alpha\de\mu=\int_{\bar{\Gamma}}\varphi\de\mu$.
D'après le lemme \ref{compactatrou}, on a donc $\alpha_*\mu=\mu$.
\end{demo}

Notons enfin $\bar{\Theta}$ le quotient de $\bar{\Gamma}$ par
l'application $\alpha$ et munissons-le de la mesure $\lambda$, image
de $\mu$ par la projection naturelle. L'espace $\bar{\Theta}$ peut
se voir naturellement comme l'espace des arêtes de $\bar{\Gamma}$ et
il peut-être muni d'une structure de graphe régulier de valence $4$.
L'image de $\Gamma$ dans $\Theta$ s'identifie alors de manière
naturelle au graphe de Sierpi\'nski $\Theta$ et forme une partie
dense de $\bar{\Theta}$. On appelle fonctions triangulaires sur
$\bar{\Theta}$ les fonctions qui proviennent de fonctions
triangulaires nulles aux sommets de leurs triangles de définition et
$\alpha$-invariantes sur $\bar{\Gamma}$. Alors, les résultats de
cette section se transportent en des résultats analogues sur
$\bar{\Theta}$.

\section{L'opérateur $\bar{\Delta}$ et ses mesures harmoniques}
\label{deltacomp}

Nous allons à présent étudier un opérateur sur $\bar{\Gamma}$ qui
est un analogue de l'opérateur $\Delta$ sur $\Gamma$.

Soit $\varphi$ une fonction borélienne sur $\bar{\Gamma}$. Pour tout
$p$ dans $\bar{\Gamma}$, on pose $\bar{\Delta}\varphi(p)=\sum_{q\sim
p}\varphi(q)$. Pour étudier les propriétés de cet opérateur sur les
fonctions triangulaires, notons encore, pour tout entier naturel
$n$, pour toute fonction $\varphi$ sur $\mathcal T_n$, pour tout $p$
dans $\mathcal T_n-\partial\mathcal T_n$,
$\Delta\varphi(p)=\sum_{q\sim p}\varphi(q)$ et, pour tout $p$ dans
$\partial\mathcal T_n$, $\Delta\varphi(p)=\varphi(p)+\sum_{q\sim
p}\varphi(q)$.

\begin{Lem} \label{fonctiontriangulaire4}
Pour tout entier $n\geq 1$, l'opérateur $\Delta$ est auto-adjoint
dans $\ell^2(\mathcal T_n)$. Si $\varphi$ est une fonction
$n$-triangulaire constante sur $\partial\mathcal T_n$, on a
$\Delta\varphi=\bar{\Delta}\varphi$.
\end{Lem}

\begin{demo} Le caractère auto-adjoint de l'opérateur $\Delta$ dans
$\mathcal T_n$ se vérifie aisément. Par ailleurs, comme $\theta_n$
induit des isomorphismes de graphes des $n$-triangles de
$\bar{\Gamma}$ dans $\mathcal T_n$, pour tout $p$ dans $\mathcal
T_n-\partial\mathcal T_n$, on a
$$\bar{\Delta}1_{\theta_n^{-1}(p)}=\sum_{q\sim
p}1_{\theta_n^{-1}(q)}=\Delta 1_{\theta_n^{-1}(p)}.$$ De même, comme
pour tout point $p$ de $\bar{\Gamma}$ qui est le sommet d'un
$n$-triangle, l'unique voisin de $p$ n'appartenant pas à ce
$n$-triangle est lui-même le sommet d'un $n$-triangle, on a
$\bar{\Delta}\sum_{p\in\partial\mathcal
T_n}1_{\theta_n^{-1}(p)}=\Delta\sum_{p\in\partial\mathcal
T_n}1_{\theta_n^{-1}(p)}$ et, donc, pour toute fonction
$n$-triangulaire $\varphi$ constante sur $\partial\mathcal T_n$,
$\bar{\Delta}\varphi=\Delta\varphi$.
\end{demo}

Nous pouvons alors énoncer les propriétés principales de
$\bar{\Delta}$ dans la

\begin{Prop}\label{operateurcomp}
L'opérateur $\bar{\Delta}$ commute à l'action de $\mathfrak S$. Il
est continu de norme $3$ dans l'espace des fonctions continues sur
$\bar{\Gamma}$ et on a $\bar{\Delta}^*\mu=3\mu$. Pour tout $1\leq
p\leq\infty$, l'opérateur $\bar{\Delta}$ est continu de norme $3$
dans ${\rm L}^p\left({\bar{\Gamma}},\mu\right)$ et, pour
$\frac{1}{p}+\frac{1}{q}=1$, pour tous $\varphi$ dans ${\rm
L}^p\left({\bar{\Gamma}},\mu\right)$ et $\psi$ dans ${\rm
L}^q\left({\bar{\Gamma}},\mu\right)$, on a
$\langle\bar{\Delta}\varphi,\psi\rangle=\langle\varphi,\bar{\Delta}\psi\rangle$.
\end{Prop}

\begin{demo} La première assertion est évidente. Comme $\bar{\Delta}$ est positif et
que $\bar{\Delta} 1=3$, $\bar{\Delta}$ est de norme $3$ dans l'espace des fonctions continues.

Rappelons qu'on a noté $\alpha$ l'application
$\bar{\Gamma}\rightarrow \bar{\Gamma}$ qui, à un point $p$, associe
son unique voisin qui n'appartient pas au triangle contenant $p$.
Pour $\varphi$ dans $\mathcal C^0\left(\bar{\Gamma}\right)$, on a
$\bar{\Delta}\varphi=3\bar{\Pi}^*\bar{\Pi}\varphi+\varphi\circ\alpha-\varphi$.
Comme les opérateurs $\bar{\Pi}$ et $\bar{\Pi}^*$ préservent $\mu$
et que, d'après le lemme \ref{involution}, l'homéomorphisme $\alpha$
préserve $\mu$, on a $\bar{\Delta}^*\mu=3\mu$.

Pour tout $1\leq p\leq\infty$, l'opérateur positif $\bar{\Delta}$
agit donc bien dans ${\rm L}^p\left({\bar{\Gamma}},\mu\right)$ et il
y est borné et de norme $3$. Soient $1<p,q<\infty$ avec
$\frac{1}{p}+\frac{1}{q}=1$. D'après le lemme
\ref{fonctiontriangulaire3}, les fonctions triangulaires nulles aux
sommets de leurs triangles de définition sont denses dans ${\rm
L}^p\left({\bar{\Gamma}},\mu\right)$ et dans ${\rm
L}^q\left({\bar{\Gamma}},\mu\right)$. D'après le lemme
\ref{fonctiontriangulaire4}, on a donc, pour tous $\varphi$ dans
${\rm L}^p\left({\bar{\Gamma}},\mu\right)$ et $\psi$ dans ${\rm
L}^q\left({\bar{\Gamma}},\mu\right)$,
$\langle\bar{\Delta}\varphi,\psi\rangle=\langle\varphi,\bar{\Delta}\psi\rangle$.
Comme les opérateurs apparaissant dans cette identité sont continus
dans ${\rm L}^1\left({\bar{\Gamma}},\mu\right)$ et dans ${\rm
L}^\infty\left({\bar{\Gamma}},\mu\right)$, elle est encore vraie
pour $p=1$ et $q=\infty$.
\end{demo}

Nous allons à présent montrer que la mesure $\mu$ est, à
multiplication par un scalaire près, la seule mesure complexe
borélienne $\lambda$ sur $\bar{\Gamma}$ telle que
$\bar{\Delta}^*\lambda=3\lambda$. Commen\c cons par nous intéresser
aux mesures de ce type qui sont $\mathfrak S$-invariantes :

\begin{Lem}\label{espacep31comp}
Soit $\lambda$ une mesure complexe borélienne $\mathfrak
S$-invariante sur $\bar{\Gamma}$ telle que
$\bar{\Delta}^*\lambda=3\lambda$. On a
$\lambda=\lambda\left(\bar{\Gamma}\right)\mu$.\end{Lem}

\begin{demo} Soit
$\varphi:p\mapsto\lambda(p),\Gamma\rightarrow\mathbb C$. Alors
$\varphi$ appartient à $\ell^1(\Gamma)$ et on a
$\Delta\varphi=3\varphi$. D'après le principe du maximum, on a donc
$\varphi=0$. Par conséquent, la restriction de $\lambda$ à
$\bigcup_{s\in\mathfrak S}s\Gamma$ est nulle. D'après le corollaire
\ref{trianglesepare} et le lemme \ref{compactatrou}, il suffit donc
de vérifier que $\lambda$ est proportionnelle à $\mu$ sur l'espace
des fonctions triangulaires. Soit $n$ un entier $\geq 1$ et, pour
tout $p$ dans $\mathcal T_n$, soit
$\varphi_n(p)=\lambda(\theta_n^{-1}(p))
=\int_{\bar{\Gamma}}1_{\theta_n^{-1}(p)}\de\lambda$. D'après le
lemme \ref{fonctiontriangulaire4}, si $p$ n'est pas un sommet de
$\mathcal T_n$, on a $\bar{\Delta}1_{\theta_n^{-1}(p)}=\Delta
1_{\theta_n^{-1}(p)}$ et, donc, comme
$\bar{\Delta}^*\lambda=3\lambda$,
$\Delta\varphi_n(p)=3\varphi_n(p)$. De plus, comme $\lambda$ est
$\mathfrak S$-invariante, $\varphi_n$ est constante sur
$\partial\mathcal T_n$. Par le principe du maximum, $\varphi_n$ est
constante. Comme $\sum_{p\in\mathcal
T_n}\varphi_n(p)=\lambda\left(\bar{\Gamma}\right)$, on a, pour tout
$p$ dans $\mathcal T_n$,
$\varphi_n(p)=\frac{1}{3^n}\lambda\left(\bar{\Gamma}\right)$, d'où
le résultat, d'après le lemme \ref{fonctiontriangulaire2}.
\end{demo}

\'Etudions maintenant l'espace propre de valeur propre $3$ pour
l'action de $\bar{\Delta}$ dans ${\rm
L}^1\left({\bar{\Gamma}},\mu\right)$. Nous aurons à utiliser le

\begin{Lem} \label{fonctiontriangulaire5}
Soient $n$ un entier $\geq 1$, $\varphi$ dans ${\rm
L}^1\left({\bar{\Gamma}},\mu\right)$ et $p$ dans $\mathcal
T_n-\partial\mathcal T_n$. On a $\Delta\mathbb
E(\varphi|\theta_n)(p)=\mathbb
E\left(\bar{\Delta}\varphi|\theta_n\right)(p)$.
\end{Lem}

\begin{demo} Notons toujours $\alpha$ et $\alpha_n$ comme dans le
lemme \ref{involution}. On a
$\bar{\Delta}\varphi=3\bar{\Pi}^*\bar{\Pi}\varphi+\varphi\circ\alpha-\varphi$
et, donc, d'après les lemmes \ref{fonctiontriangulaire3} et
\ref{involution}, \begin{multline*}\mathbb
E\left(\bar{\Delta}\varphi|\theta_n\right)(p)=\Pi^*\Pi\mathbb
E\left(\varphi|\theta_n\right)(p)+\mathbb
E\left(\varphi|\theta_n\right)(\alpha_n(p))-\mathbb
E\left(\varphi|\theta_n\right)(p)\\ =\Delta\mathbb
E\left(\varphi|\theta_n\right)(p),\end{multline*} ce qu'il fallait
démontrer.
\end{demo}

Pour tout $n\geq 1$, on note $H_n$ l'espace des fonctions $\varphi$
sur $\mathcal T_n$ telles que, pour tout $p$ dans $\mathcal T_n$ qui
ne soit pas un sommet, on ait $\Delta\varphi(p)=3\varphi(p)$.
D'après le lemme \ref{fonctiontriangulaire5}, pour tout $\varphi$
dans ${\rm L}^1\left({\bar{\Gamma}},\mu\right)$, si
$\bar{\Delta}\varphi=3\varphi$, on a $\mathbb E(\varphi|\theta_n)\in
H_n$. On identifie $\mathbb C^3$ et l'espace des fonctions à valeurs
complexes sur $\mathcal T_1$ en considérant $(1,0,0)$ (resp.
$(0,1,0)$, resp. $(0,0,1)$) comme la fonction caractéristique du
singleton $\{a_1\}$ (resp. $\{b_1\}$, resp. $\{c_1\}$) et on note
$\eta_n$ l'application linéaire $\mathfrak S$-équivariante
$H_n\rightarrow\mathbb
C^3,\varphi\mapsto(\varphi(a_n),\varphi(b_n),\varphi(c_n))$. Par
ailleurs, on note $(s_n)_{n\geq 1}$ la suite de nombres réels telle
que $s_1=1$ et que, pour tout $n\geq 1$, on ait
$s_{n+1}=\frac{3s_n}{3s_n+5}$. On montre facilement qu'on a
$s_n\td{n}{\infty}0$. Soit $\mathbb C^3_0$ l'ensemble des éléments
de $\mathbb C^3$ dont la somme des coordonnées est nulle. Nous avons
le

\begin{Lem} \label{espacep32comp}
Soit $n\geq 1$. Pour tous $\varphi$ dans $H_n$ et $p$ dans $\mathcal
T_n$, on a
$\abs{\varphi(p)}\leq\max\{\abs{\varphi(a_n)},\abs{\varphi(b_n)},\abs{\varphi(c_n)}\}$.
En particulier, l'application $\eta_n$ est un isomorphisme.
Supposons $n\geq 2$. Soient $\varphi$ dans $H_n$ tel que
$\eta_n(\varphi)$ appartienne à $\mathbb C^3_0$ et $\psi=\mathbb
E(\varphi|\theta_{n-1})$. Alors $\psi$ appartient à $H_{n-1}$ et on
a $\eta_{n-1}(\psi)=\frac{2}{3s_{n-1}+5}\eta_n(\varphi)$.\end{Lem}

\begin{demo} La majoration découle du
principe du maximum appliqué à l'opérateur $\frac{1}{3}\Delta$. Elle
implique que, pour tout $n\geq 1$, l'opérateur $\eta_n$ est
injectif. Soit $\varphi$ une fonction sur $\mathcal T_n$ et, pour
tout $p$ dans $\mathcal T_n$, posons
$\delta_n\varphi(p)=\Delta\varphi(p)-3\varphi(p)$ si $p$ n'est pas
un sommet et $\delta_n\varphi(p)=\varphi(p)$ si $p$ est un sommet.
Comme $\eta_n$ est injectif, $\delta_n$ l'est aussi et est donc un
isomorphisme ; en particulier, $\eta_n$ est surjectif et, donc,
c'est un isomorphisme.

Pour tout $n\geq 2$, posons $t_n=\frac{3s_{n-1}+2}{3s_{n-1}+5}$ et
$u_n=\frac{1}{3s_{n-1}+5}$. Rappelons que, comme dans le corollaire
\ref{decoupetriangles}, si $\mathcal S$ est un $n$-triangle et si
$p$ et $q$ sont deux sommets de $\mathcal S$, on note $pq$ l'unique
point de $\mathcal S$ appartenant au $(n-1)$-triangle contenant $p$
et dont un voisin appartient au $(n-1)$-triangle contenant $q$.
Soient $d_n$ et $e_n$ les deux voisins de $a_n$ dans $\mathcal T_n$.
Montrons par récurrence sur $n\geq 2$ que, pour tout $\varphi$ dans
$H_n$, on a
$\varphi(d_n)+\varphi(e_n)=s_n(\varphi(b_n)+\varphi(c_n))+2(1-s_n)\varphi(a_n)$
et
$\varphi(a_nb_n)=t_n\varphi(a_n)+u_n(2\varphi(b_n)+\varphi(c_n))$.
Pour $n=2$, c'est un calcul immédiat. Si $n\geq 3$ et si la formule
est vraie pour $n-1$, donnons-nous une fonction $\varphi$ dans
$H_{n}$. Alors, comme $\Delta\varphi(a_nb_n)=3\varphi(a_nb_n)$, en
appliquant la récurrence à la restriction de $\varphi$ au
$(n-1)$-triangle contenant $a_n$, on a
$$s_{n-1}(\varphi(a_n)+\varphi(a_nc_n))+2(1-s_{n-1})\varphi(a_nb_n)+\varphi(b_na_n)
=3\varphi(a_nb_n).$$ Comme $\eta_n$ est un isomorphisme, il existe
un unique $(x,y,z)$ dans $\mathbb C^3$ tel que, pour tout $\varphi$
dans $H_n$, on ait
$\varphi(a_nb_n)=x\varphi(a_n)+y\varphi(b_n)+z\varphi(c_n)$. Comme
$\eta_n$ est $\mathfrak S$-équivariant, on a, pour tout $\varphi$
dans $H_n$,
$\varphi(b_na_n)=x\varphi(b_n)+y\varphi(a_n)+z\varphi(c_n)$ et
$\varphi(a_nc_n)=x\varphi(a_n)+y\varphi(c_n)+z\varphi(b_n)$. Il
vient
\begin{align*}
s_{n-1}(1+x)+2(1-s_{n-1})x+y&=3x\\
s_{n-1}z+2(1-s_{n-1})y+x&=3y\\
s_{n-1}y+2(1-s_{n-1})z+z&=3z.
\end{align*}
En résolvant ce système, on obtient $x=t_n$, $y=2u_n$ et $z=u_n$.
Enfin, par récurrence, on a
\begin{align*}
\varphi(d_n)+\varphi(e_n)&=s_{n-1}(\varphi(a_nb_n)+\varphi(a_nc_n))
+2(1-s_{n-1})\varphi(a_n)\\
&=3s_{n-1}u_n(\varphi(b_n)+\varphi(c_n))
+2(1-s_{n-1}+s_{n-1}t_n)\varphi(a_n),
\end{align*}
d'où le résultat puisque $3s_{n-1}u_n=s_n=s_{n-1}(1-t_n)$.

Alors, si $\psi=\mathbb E(\varphi|\theta_{n-1})$, on a, d'après le
lemme \ref{fonctiontriangulaire3}, pour tout $p$ dans $\mathcal
T_{n-1}$, $\psi(p)=\frac{1}{3}\sum_{\theta_{n-1}(q)=p}\varphi(q)$.
Comme $\theta_{n-1}$ induit un isomorphisme de graphes de chacun des
$(n-1)$-triangles de $\mathcal T_n$ dans $\mathcal T_{n-1}$, on en
déduit que $\psi$ appartient à $H_{n-1}$ et que, en particulier,
d'après le lemme \ref{facteurtriangle}, si $\eta_n(\varphi)$ est
dans $\mathbb C^3_0$, on a
\begin{align*}\psi(a_n)&=\frac{1}{3}(\varphi(a_n)+\varphi(b_na_n)+\varphi(c_na_n))\\
&=\frac{1}{3}\left((1+4u_n)\varphi(a_n)+\left(t_n+u_n\right)(\varphi(b_n)+\varphi(c_n))\right)\\
&=\frac{1+4u_n-t_n-u_n}{3}\varphi(a_n)=\frac{2}{3s_{n-1}+5}\varphi(a_n),\end{align*}
où, pour l'avant-dernière égalité, on a utilisé la relation
$\varphi(a_n)+\varphi(b_n)+\varphi(c_n)=0$. Par $\mathfrak
S$-équivariance, on a la formule analogue aux autres sommets de
$\mathcal T_n$ et, donc,
$\eta_{n-1}(\psi)=\frac{2}{3s_{n-1}+5}\eta_n(\varphi)$.
\end{demo}

\begin{Cor}\label{espacep33comp}
Soit $\varphi$ dans ${\rm L}^\infty\left({\bar{\Gamma}},\mu\right)$
tel que $\bar{\Delta}\varphi=3\varphi$ et que $\sum_{s\in\mathfrak
S}\varphi\circ s=0$. On a $\varphi=0$.\end{Cor}

\begin{demo} Pour tout entier $n\geq
1$, posons $\varphi_n=\mathbb E(\varphi|\theta_n)$. D'après le lemme
\ref{fonctiontriangulaire5}, on a $\varphi_n\in H_n$. Comme
$\sum_{s\in\mathfrak S}\varphi\circ s=0$, on a
$\eta_n(\varphi_n)\in\mathbb C^3_0$. Par conséquent, si $n\geq 2$,
d'après le lemme \ref{espacep32comp}, comme $\varphi_{n-1}=\mathbb
E(\varphi_n|\theta_{n-1})$, on a
$\eta_{n-1}(\varphi_{n-1})=\frac{2}{3s_{n-1}+5}\eta_n(\varphi_n)$.
Or, pour tout $n\geq 1$, on a $\N{\varphi_n}_\infty\leq
\N{\varphi}_\infty$ et, comme $s_n\td{n}{\infty}0$,
$\prod_{n=1}^\infty\frac{2}{3s_n+5}=0$. Par conséquent, on a
nécessairement, pour tout $n\geq 1$, $\eta_n(\varphi_n)=0$, donc,
d'après le lemme \ref{espacep32comp}, $\varphi_n=0$, et, d'après le
lemme \ref{fonctiontriangulaire3}, $\varphi=0$.
\end{demo}

Nous pouvons alors décrire les vecteurs propres de valeur propre $3$
dans ${\rm L}^1\left({\bar{\Gamma}},\mu\right)$ :

\begin{Lem}\label{espacep34comp}
Soit $\varphi$ dans ${\rm L}^1\left({\bar{\Gamma}},\mu\right)$ tel
que $\bar{\Delta}\varphi=3\varphi$. La fonction $\varphi$ est
constante $\mu$-presque partout.\end{Lem}

\begin{demo} On peut supposer que $\varphi$ est à valeurs réelles.
Ramenons-nous alors au cas où $\varphi$ est positive. Comme
$\bar{\Delta}$ est positif, on a
$\bar{\Delta}\abs{\varphi}\geq\abs{\bar{\Delta}\varphi}=3\abs{\varphi}$
et, donc, comme $\bar{\Delta}$ est de norme $\leq 3$,
$\bar{\Delta}\abs{\varphi}=3\abs{\varphi}$. Quitte à étudier les
fonctions $\abs{\varphi}-\varphi$ et $\abs{\varphi}+\varphi$, on
peut donc supposer qu'on a $\varphi\geq 0$. Posons
$\psi=\sum_{s\in\mathfrak S}\varphi\circ s$. La mesure
$\lambda=\psi\mu$ est $\mathfrak S$-invariante et on a
$\bar{\Delta}^*\lambda=3\lambda$. D'après le lemme
\ref{espacep31comp}, $\lambda$ est proportionnelle à $\mu$,
c'est-à-dire que $\psi$ est constante $\mu$-presque partout. En
particulier, $\psi$ est dans ${\rm
L}^\infty\left({\bar{\Gamma}},\mu\right)$. Comme on a
$0\leq\varphi\leq\psi$, $\varphi$ est dans ${\rm
L}^\infty\left({\bar{\Gamma}},\mu\right)$. Alors, d'après le
corollaire \ref{espacep33comp}, on a
$\varphi-\frac{1}{6}\psi=0$.\end{demo}

Pour étendre ce résultat à toutes les mesures complexes sur
$\bar{\Gamma}$, nous utiliserons un lemme général sans doute
classique. Soit $X$ un espace métrique compact. On munit l'espace
$\mathcal C^0(X)$ de la topologie de la convergence uniforme. Si
$\lambda$ est une mesure complexe borélienne sur $X$, rappelons que
la variation totale $\abs{\lambda}$ de $\lambda$ est la mesure
borélienne positive et finie sur $X$ telle que, pour toute fonction
continue positive $g$ sur $X$, on ait
$$\int_X g\de\abs{\lambda}=
\sup_{\substack{h\in \mathcal C^0(X)\\ \N{h}_\infty\leq
1}}\abs{\int_X gh\de\lambda}$$ (on pourra se référer à \cite[Chapter
6]{Rud}). En particulier, $\abs{\lambda}$ est la plus petite mesure
de Radon positive telle que, pour toute fonction continue positive
$g$ sur $X$, on ait $\abs{\int_X g\de\lambda}\leq\int_X
g\de\abs{\lambda}$.

\begin{Lem}\label{vartotinvariante}
Soient $X$ un espace métrique compact et $P$ un opérateur positif de
norme $\leq 1$ sur l'espace des fonctions continues sur $X$. Pour
toute mesure complexe borélienne $\lambda$ sur $X$, on a
$\abs{P^*\lambda}\leq P^*\abs{\lambda}$. En particulier, si
$P^*\lambda=\lambda$, on a
$P^*\abs{\lambda}=\abs{\lambda}$.\end{Lem}

\begin{demo} Pour toute fonction continue positive $g$ sur $X$, on a $Pg\geq 0$, donc
$\abs{\int_X Pg\de\lambda}\leq \int_X Pg\de\abs{\lambda}$. Comme la
mesure $P^*\abs{\lambda}$ est positive, il vient
$\abs{P^*\lambda}\leq P^*\abs{\lambda}$. Si $P^*\lambda=\lambda$, on
a $\abs{\lambda}\leq P^*\abs{\lambda}$, d'où l'égalité, puisque $P$
est de norme $1$.\end{demo}

Nous en déduisons enfin la

\begin{Prop}\label{espacep35comp}
Soit $\lambda$ une mesure complexe borélienne sur $\bar{\Gamma}$
telle que $\bar{\Delta}^*\lambda=3\lambda$. On a
$\lambda=\lambda\left(\bar{\Gamma}\right)\mu$.\end{Prop}

\begin{demo} On peut supposer que $\lambda$ est à valeurs réelles.
D'après le lemme \ref{vartotinvariante}, appliqué à l'opérateur
$\frac{1}{3}\bar{\Delta}$, on a
$\bar{\Delta}^*\abs{\lambda}=3\abs{\lambda}$ et, donc, quitte à
étudier les mesures $\abs{\lambda}-\lambda$ et
$\abs{\lambda}+\lambda$, on peut supposer que $\lambda$ est
positive. Alors, d'après le lemme \ref{espacep31comp}, la mesure
$\sum_{s\in\mathfrak S}s_*\lambda$ est proportionnelle à $\mu$.
Comme on a $0\leq\lambda\leq\sum_{s\in\mathfrak S}s_*\lambda$,
$\lambda$ est absolument continue par rapport à $\mu$. D'après le
lemme \ref{espacep34comp}, $\lambda$ est donc proportionnelle à
$\mu$.
\end{demo}

\section{Spectre et mesures spectrales de $\bar{\Gamma}$}
\label{secspectrecomp}

Nous allons à présent aborder l'étude spectrale de l'opérateur
$\bar{\Gamma}$. Commen\c cons par remarquer que, comme dans le lemme
\ref{relation1}, on a le

\begin{Lem}\label{relation5}
On a
$(\bar{\Delta}^2-\bar{\Delta}-3)\bar{\Pi}^*=\bar{\Pi}^*\bar{\Delta}$
et
$\bar{\Pi}(\bar{\Delta}^2-\bar{\Delta}-3)=\bar{\Delta}\bar{\Pi}$.\end{Lem}

Rappelons que nous avons noté $\alpha$ l'application qui, à un point
$p$ de $\bar{\Gamma}$ associe le voisin de $p$ qui n'appartient pas
au $1$-triangle contenant $p$. Comme dans la section
\ref{secspectre}, on déduit du lemme \ref{relation5} le

\begin{Cor}\label{spectrecomp}
Le spectre de $\bar{\Delta}$ est la réunion de $\Lambda$ et de
$\bigcup_{n\in\mathbb N}f^{-n}(0)$. L'espace propre associé à la
valeur propre $-2$ est l'espace des fonctions $\varphi$ dans ${\rm
L}^2\left(\bar{\Gamma},\mu\right)$ telles que $\bar{\Pi}\varphi=0$
et que $\varphi\circ\alpha=-\varphi$. L'espace propre associé à la
valeur propre $0$ est l'espace des fonctions $\varphi$ dans ${\rm
L}^2\left(\bar{\Gamma},\mu\right)$ telles que $\bar{\Pi}\varphi=0$
et que $\varphi\circ\alpha=\varphi$.
\end{Cor}

\begin{demo} Soient, comme dans le corollaire
\ref{transformespectre1}, $K=\bar{\Pi}^*{\rm
L}^2\left(\bar{\Gamma},\mu\right)$ et $H$ le sous-espace fermé de
${\rm L}^2\left(\bar{\Gamma},\mu\right)$ engendré par $K$ et par
$\bar{\Delta}K$. D'après le lemme \ref{relation5}, on a
$f\left(\bar{\Delta}\right)K\subset K$ et, comme $\bar{\Pi}^*$ est
une isométrie de ${\rm L}^2\left(\bar{\Gamma},\mu\right)$ dans $K$,
le spectre de $f\left(\bar{\Delta}\right)$ dans $K$ est égal au
spectre de $\bar{\Delta}$ dans ${\rm
L}^2\left(\bar{\Gamma},\mu\right)$. Nous allons chercher à appliquer
le lemme \ref{spectral2} à l'opérateur $\bar{\Delta}$ dans $H$. Pour
cela, montrons que $\bar{\Delta}^{-1}K\cap K$ est réduit aux
fonctions constantes. Soient $\varphi$ et $\psi$ dans ${\rm
L}^2\left(\bar{\Gamma},\mu\right)$ telles que
$\bar{\Delta}\bar{\Pi}^*\varphi=\bar{\Pi}^*\psi$. Pour tout entier
$n\geq 1$, posons $\varphi_n=\mathbb E(\varphi|\theta_n)$ et
$\psi_n=\mathbb E(\psi|\theta_n)$. D'après les lemmes
\ref{fonctiontriangulaire3} et \ref{fonctiontriangulaire5}, pour
tout $p$ dans $\mathcal T_{n+1}-\partial\mathcal T_{n+1}$, on a
$\Delta\Pi^*\varphi_n(p)=\Pi^*\psi_n(p)$. En raisonnant comme dans
la démonstration du corollaire \ref{transformespectre1}, on en
déduit que, pour tout $q$ dans $\mathcal T_n$, $\varphi_n$ est
constante sur les voisins de $q$. Comme tout point de $\mathcal T_n$
est contenu dans un triangle et que $\mathcal T_n$ est connexe,
$\varphi_n$ est constante. Comme, d'après le lemme
\ref{fonctiontriangulaire3}, on a $\varphi_n\td{n}{\infty}\varphi$
dans ${\rm L}^2\left(\bar{\Gamma},\mu\right)$, $\varphi$ est
constante. L'espace $\bar{\Delta}^{-1}K\cap K$ est donc égal à la
droite des fonctions constantes. D'après les lemmes \ref{spectral2}
et \ref{relation5}, le spectre de $\bar{\Delta}$ dans $H$ est donc
égal à la réunion de $\{3\}$ et de l'image inverse par $f$ du
spectre de $\bar{\Delta}$ dans l'espace des fonctions d'intégrale
nulle dans ${\rm L}^2\left(\bar{\Gamma},\mu\right)$.

Par ailleurs, en raisonnant comme dans le lemme
\ref{spectreresiduel}, on voit que l'orthogonal $L$ de $H$ dans
${\rm L}^2\left(\bar{\Gamma},\mu\right)$ est la somme directe de
l'espace $L_{-2}$ des éléments $\varphi$ de ${\rm
L}^2\left(\bar{\Gamma},\mu\right)$ tels que $\bar{\Pi}\varphi=0$ et
que $\varphi\circ\alpha=-\varphi$ et de l'espace $L_{0}$ des
éléments $\varphi$ de ${\rm L}^2\left(\bar{\Gamma},\mu\right)$ tels
que $\bar{\Pi}\varphi=0$ et que $\varphi\circ\alpha=\varphi$. On a
$\bar{\Delta}=-2$ sur $L_{-2}$ et $\bar{\Delta}=0$ sur $L_0$. En
raisonnant comme dans le lemme \ref{existevecteurpropre} et en
utilisant le lemme \ref{fonctiontriangulaire4}, on voit que ces deux
espaces ne sont pas réduits à $\{0\}$, puisqu'ils contiennent des
fonctions triangulaires. Comme dans la démonstration du corollaire
\ref{spectre}, on en déduit que le spectre de $\bar{\Delta}$ dans
${\rm L}^2\left(\bar{\Gamma},\mu\right)$ est égal à la réunion de
$\Lambda$ et de $\bigcup_{n\in\mathbb N}f^{-n}(0)$.

Enfin, comme dans la démonstration du lemme \ref{spectreresiduel},
il nous reste à montrer que $L_{-2}$ et $L_0$ sont exactement les
espaces propres de $\bar{\Delta}$ associés aux valeurs propres $-2$
et $0$, c'est-à-dire que $\bar{\Delta}$ ne possède pas de vecteur
propre de valeur propre $-2$ ou $0$ dans $H$. Soit $\varphi$ dans
$H$ tel que $\bar{\Delta}\varphi=-2\varphi$. D'après le lemme
\ref{relation5}, on a
$\bar{\Delta}\bar{\Pi}\varphi=3\bar{\Pi}\varphi$ et, donc, d'après
le lemme \ref{espacep34comp}, $\bar{\Pi}\varphi$ est constante.
Comme $\varphi$ est orthogonale aux fonctions constantes, on a
$\bar{\Pi}\varphi=0$ et
$\bar{\Pi}\bar{\Delta}\varphi=-2\bar{\Pi}\varphi=0$. Comme $\varphi$
est dans $H$, il vient $\varphi=0$. De même, si $\varphi$ est dans
$H$ et si $\bar{\Delta}\varphi=0$, on a
$\bar{\Delta}\bar{\Pi}\varphi=-3\bar{\Pi}\varphi$. Or, par un calcul
immédiat, $-3$ n'appartient pas au spectre de $\bar{\Delta}$. Il
vient $\bar{\Pi}\varphi=0$ et, donc, $\varphi=0$, ce qu'il fallait
démontrer.
\end{demo}

Nous avons aussi un analogue du lemme \ref{relation2} :

\begin{Lem} \label{relation6} On a
$\bar{\Pi}\bar{\Delta}\bar{\Pi}^*=2+\frac{1}{3}\bar{\Delta}$ et,
donc, pour tous $\varphi$ et $\psi$ dans ${\rm
L}^2\left(\bar{\Gamma},\mu\right)$,
\begin{multline*}\left\langle\bar{\Delta}\bar{\Pi}^*\varphi,\bar{\Pi}^*\psi\right\rangle=
2\left\langle\varphi,\psi\right\rangle+\frac{1}{3}\left\langle\bar{\Delta}\varphi,\psi\right\rangle
\\=2\left\langle\bar{\Pi}^*\varphi,\bar{\Pi}^*\psi\right\rangle+
\frac{1}{3}\left\langle\left(\bar{\Delta}^2-\bar{\Delta}-3\right)\bar{\Pi}^*\varphi,\bar{\Pi}^*\psi\right\rangle.\end{multline*}
\end{Lem}

Comme dans la section \ref{mesurespectrale}, on en déduit le

\begin{Cor}\label{mesurespectralecomp}
Soient $\varphi$ dans  ${\rm L}^2\left(\bar{\Gamma},\mu\right)$, $\mu$ la mesure spectrale de
$\varphi$ pour $\bar{\Delta}$ dans ${\rm L}^2\left(\bar{\Gamma},\mu\right)$  et $\nu$ la mesure
spectrale de $\bar{\Pi}^*\varphi$ pour $\bar{\Delta}$ dans ${\rm L}^2\left(\bar{\Gamma},\mu\right)$.
Alors, on a $\nu(\frac{1}{2})=0$ et, si, pour tout
$x\neq\frac{1}{2}$, on pose $\tau(x)=\frac{x(x+2)}{3(2x-1)}$, on a
$\nu=L_{f,\tau}^*\mu$.
\end{Cor}

\section{Fonctions propres dans ${\rm L}^2\left(\bar{\Gamma},\mu\right)$}
\label{valeurpcomp}

Dans cette section, nous allons reprendre le plan de la section
\ref{valeurp}, pour décrire les espaces propres de $\bar{\Delta}$
dans ${\rm L}^2\left(\bar{\Gamma},\mu\right)$. Comme dans la section
\ref{valeurp}, en utilisant les lemmes \ref{relation5} et
\ref{relation6}, on montre l'analogue suivant du lemme
\ref{transformevalp} :

\begin{Lem}\label{transformevalpcomp}
Soit $H$ le sous-espace fermé de ${\rm
L}^2\left(\bar{\Gamma},\mu\right)$ engendré par l'image de
$\bar{\Pi}^*$ et par celle de $\bar{\Delta}\bar{\Pi}^*$. Alors, pour
tout $x$ dans $\mathbb R-\{0,-2\}$, $x$ est valeur propre de
$\bar{\Delta}$ dans $H$ si et seulement si $y=f(x)$ est valeur
propre de $\bar{\Delta}$ dans ${\rm
L}^2\left(\bar{\Gamma},\mu\right)$. Dans ce cas, l'application
$\bar{R}_x$ qui, à une fonction propre $\varphi$ de valeur propre
$y$ dans ${\rm L}^2\left(\bar{\Gamma},\mu\right)$, associe
$(x-1)\bar{\Pi}^*\varphi+\bar{\Delta}\bar{\Pi}^*\varphi$ induit un
isomorphisme entre l'espace propre de valeur propre $y$ dans ${\rm
L}^2\left(\bar{\Gamma},\mu\right)$ et l'espace propre de valeur
propre $x$ dans $H$ et, pour tout $\varphi$, on a
$\N{\bar{R}_x\varphi}_2^2=\frac{1}{3}x(x+2)(2x-1)\N{\varphi}_2^2$.
\end{Lem}

Pour décrire les fonctions propres de valeur propre dans
$\bigcup_{n\in\mathbb N}f^{-n}(0)$, nous allons raisonner comme dans
la section \ref{valeurp}. Pour cela, remarquons à nouveau que, pour
tout entier $n\geq 1$, l'espace des arêtes extérieures aux
$n$-triangles de $\bar{\Gamma}$ s'identifie de manière naturelle à
$\bar{\Theta}$. Si $\varphi$ est une fonction sur $\bar{\Gamma}$ qui
est constante sur les arêtes extérieures aux $n$-triangles, on note
$\bar{P}_n(\varphi)$ la fonction qui, en tout point de
$\bar{\Theta}$, a pour valeur la valeur de $\varphi$ sur l'arête
extérieure aux $n$-triangles de $\bar{\Gamma}$ associée. Par
ailleurs, on note encore $\bar{\Delta}$ l'opérateur qui, à une
fonction $\psi$ sur $\bar{\Theta}$, associe la fonction dont la
valeur en un point $p$ de $\bar{\Theta}$ est $\sum_{q\sim
p}\psi(q)$. Cet opérateur vérifie $\bar{\Delta}^*\lambda=4\lambda$
et il est auto-adjoint de norme $4$ dans ${\rm
L^2}\left(\bar{\Theta},\lambda\right)$, où $\lambda$ est la mesure
sur $\bar{\Theta}$ introduite à la fin de la section
\ref{fonctiontriangulaire}.

\begin{Lem} \label{espacep01comp}
L'application $\bar{P}_2$ induit un isomorphisme d'espaces de Banach
de l'espace propre de ${\rm L}^2\left(\bar{\Gamma},\mu\right)$
associé à la valeur propre $0$ dans ${\rm
L^2}\left(\bar{\Theta},\lambda\right)$. On note $\bar{Q}_{0}$ sa
réciproque. Pour tout $\psi$ dans ${\rm
L^2}\left(\bar{\Theta},\lambda\right)$, on a
$\N{\bar{Q}_{0}\psi}^2_{{\rm L}^2\left(\bar{\Gamma},\mu\right)}
=\frac{1}{2}\N{\psi}^2_{{\rm L^2}\left(\bar{\Theta},\lambda\right)}
-\frac{1}{12}\langle\Delta\psi,\psi\rangle_{{\rm
L^2}\left(\bar{\Theta},\lambda\right)}$.
\end{Lem}

\begin{demo} On raisonne comme dans le lemme \ref{espacep01}
en utilisant la caractérisation des fonctions propres de valeur
propre $0$ donnée dans le corollaire \ref{spectrecomp}. La formule
se vérifie aisément sur les fonctions triangulaires nulles aux
sommets de leur triangle de définition et le cas général s'en déduit
par densité.
\end{demo}

Rappelons que, pour $x$ dans $\bigcup_{n\in\mathbb N}f^{-n}(0)$, on
a noté $n(x)$ l'entier $n$ tel que $f^n(x)=0$ et
$$\kappa(x)=\prod_{k=0}^{n(x)-1}\frac{f^k(x)(2f^k(x)-1)}{f^k(x)+2}.$$
Des lemmes \ref{transformevalpcomp} et \ref{espacep01comp}, on
déduit l'analogue suivant de la proposition \ref{espacep02} :

\begin{Prop} \label{espacep02comp}
Soit $x$ dans $\bigcup_{n\in\mathbb N}f^{-n}(0)$. Les fonctions
propres de valeur propre $x$ dans ${\rm
L}^2\left(\bar{\Gamma},\mu\right)$ sont constantes sur les arêtes
extérieures aux $(n(x)+2)$-triangles dans $\bar{\Gamma}$.
L'application $\bar{P}_{n(x)+2}$ induit un isomorphisme d'espaces de
Banach de l'espace propre de ${\rm
L}^2\left(\bar{\Gamma},\mu\right)$ associé à la valeur propre $x$
dans ${\rm L}^2\left(\bar{\Theta},\lambda\right)$. On note
$\bar{Q}_{x}$ sa réciproque. Alors, pour tout $\psi$ dans ${\rm
L}^2\left(\bar{\Theta},\lambda\right)$, on a
$$\N{\bar{Q}_{x}\psi}^2_{{\rm L}^2\left(\bar{\Gamma},\mu\right)}
=\frac{\kappa(x)}{3^{n(x)}}\left(\frac{1}{2}\N{\psi}^2_{{\rm
L}^2\left(\bar{\Theta},\lambda\right)}
-\frac{1}{12}\langle\bar{\Delta}\psi,\psi\rangle_{{\rm
L}^2\left(\bar{\Theta},\lambda\right)}\right).$$
\end{Prop}

\begin{Cor} \label{espacep03comp}
Pour tout $x$ dans $\bigcup_{n\in\mathbb N}f^{-n}(0)$, l'espace
propre associé à $x$ dans ${\rm L}^2\left(\bar{\Gamma},\mu\right)$
est de dimension infinie et engendré par des fonctions triangulaires
nulles aux sommets de leur triangle de définition.
\end{Cor}

Comme dans la section \ref{valeurp}, la description des espaces
propres associés aux éléments de $\bigcup_{n\in\mathbb N}f^{-n}(-2)$
est moins précise.

Commen\c cons par le cas de la valeur propre $-2$. Nous allons avoir
besoin de renseignements supplémentaires sur les fonctions
triangulaires qui sont vecteurs propres de valeur propre $-2$. Pour
cela, donnons-nous un entier $n\geq 1$ et un $n$-triangle $\mathcal
S$ et notons, pour tout entier $n\geq 1$, $E_{\mathcal S}$ l'espace
des fonctions $\varphi$ sur $\mathcal S$ telles que $\Pi\varphi=0$
et que, pour tout point $p$ de $\mathcal S$ qui ne soit pas un
sommet, si $q$ est le voisin de $p$ n'appartenant pas au triangle
contenant $p$, on a $\varphi(q)=-\varphi(p)$. Si $\mathcal S$ est
$\mathcal T_n$, on note $E_n$ pour $E_{\mathcal T_n}$. D'après les
lemmes \ref{fonctiontriangulaire3} et \ref{involution} et le
corollaire \ref{spectrecomp}, si $\varphi$ est un élément de ${\rm
L}^2\left(\bar{\Gamma},\mu\right)$ tel que
$\bar{\Delta}\varphi=-2\varphi$, pour tout entier $n\geq 1$, on a
$\mathbb E(\varphi|\theta_n)\in E_n$. En raisonnant comme dans le
lemme \ref{espacep-21}, on montre le

\begin{Lem}\label{sommenulle}
Soient $n\geq 1$, $\mathcal S$ un $n$-triangle de sommets $p$, $q$ et $r$, et
$\varphi$ dans $E_{\mathcal S}$. On a
$\varphi(p)+\varphi(q)+\varphi(r)=0$.\end{Lem}

L'espace $\mathbb C^3_0=\{(s,t,u)\in\mathbb R^3|s+t+u=0\}$ est
stable par l'action de $\mathfrak S$ sur $\mathbb C^3$. On le munit
de la norme hermitienne $\mathfrak S$-invariante $\N{.}_0$ telle
que, pour tout $(s,t,u)$ dans $\mathbb C^3_0$, on ait
$\N{(s,t,u)}_0^2=\frac{1}{3}\left(\abs{s}^2+\abs{t}^2+\abs{u}^2\right)$.
Pour tout $n\geq 1$, on note $\rho_n$ l'application linéaire
$\mathfrak S$-équivariante $E_n\rightarrow\mathbb
C^3_0,\varphi\mapsto(\varphi(a_n),\varphi(b_n),\varphi(c_n))$, $F_n$
le noyau de $\rho_n$ et $G_n$ l'orthogonal de $F_n$ dans $E_n$ pour
la norme de ${\rm L}^2\left(\bar{\Gamma},\mu\right)$. D'après le
lemme \ref{fonctiontriangulaire4}, les éléments de $F_n$ sont des
vecteurs propres de valeur propre $-2$ de $\bar{\Delta}$.

\begin{Lem}\label{espacep-20comp}
Soit $n\geq 1$. On a $\dim F_n=\frac{1}{2}(3^{n-1}-1)$.
L'application $\rho_n$ est surjective et, pour tout $\varphi$ dans
$G_n$, on a $\N{\varphi}^2_{{\rm L}^2\left(\bar{\Gamma},\mu\right)}
=\left(\frac{5}{9}\right)^{n-1}\N{\rho_n(\varphi)}^2_0$. Enfin, si
$n\geq 2$ et si $\psi=\mathbb E(\varphi|\theta_{n-1})$, $\psi$
appartient à $G_{n-1}$ et
$\rho_{n-1}(\psi)=\frac{2}{3}\rho_n(\varphi)$.\end{Lem}

\begin{demo} Soient $n\geq 1$, $\mathcal S$ un $n$-triangle et $p$
et $q$ des sommets distincts de $\mathcal S$. Définissons une
fonction $\varphi_{\mathcal S}^{p,q}$ sur $\mathcal S$ de la fa\c
con suivante. Si $n=1$, on pose $\varphi^{p,q}_{\mathcal S}(p)=1$ et
$\varphi^{p,q}_{\mathcal S}(q)=-1$ et on dit que $\varphi_{\mathcal
S}^{p,q}$ est nulle au troisième point de $\mathcal S$. Si $n\geq
2$, notons toujours $pq$ et $qp$ les point décrits dans le
corollaire \ref{decoupetriangles} : le point $pq$ appartient au
$(n-1)$-triangle $\mathcal P$ de sommet $p$ dans $\mathcal S$, le
point $qp$ appartient au $(n-1)$-triangle $\mathcal Q$ de sommet $q$
dans $\mathcal S$ et les points $pq$ et $qp$ sont voisins. On
définit alors $\varphi^{p,q}_{\mathcal S}$ comme la fonction dont la
restriction à $\mathcal P$ est $\varphi^{p,pq}_{\mathcal P}$, dont
la restriction à $\mathcal Q$ est $\varphi^{qp,q}_{\mathcal Q}$ et
dont la restriction au troisième $(n-1)$-triangle de $\mathcal S$
est nulle. On vérifie aisément par récurrence que
$\varphi^{p,q}_{\mathcal S}$ appartient à $E_\mathcal S$. Si
$\mathcal S=\mathcal T_n$, on note $\varphi^{p,q}_{n}$ pour
$\varphi^{p,q}_{\mathcal T_n}$. Comme on a
$\rho_n(\varphi^{a_n,b_n}_{n})=(1,-1,0)$ et
$\rho_n(\varphi^{a_n,c_n}_{n})=(1,0,-1)$, l'application $\rho_n$ est
surjective.

Pour $n\geq 2$, notons $\psi_n$ la fonction sur $\mathcal T_n$ dont
la restriction au $(n-1)$-triangle $\mathcal A_n$ (resp. $\mathcal
B_n$, resp. $\mathcal C_n$) de sommet $a_n$ (resp. $b_n$, resp.
$c_n$) est égale à $\varphi^{a_nb_n,a_nc_n}_{\mathcal A_n}$ (resp.
$\varphi^{b_nc_n,b_na_n}_{\mathcal B_n}$, resp.
$\varphi^{c_na_n,c_nb_n}_{\mathcal C_n}$). Alors, on vérifie
aisément que $\psi_n$ appartient à $F_n$.

Ces fonctions sont représentées à la figure \ref{fonctionphi}.

\begin{figure}\begin{center}\input{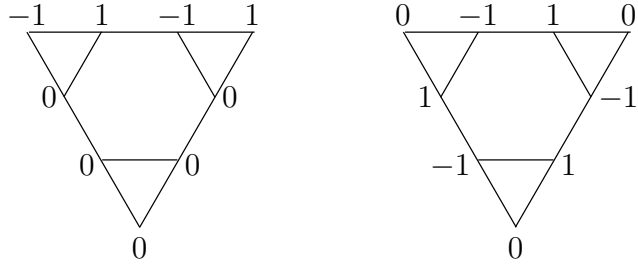}\caption{Les fonctions $\varphi_2^{a_2,b_2}$ et $\psi_2$}
\label{fonctionphi}\end{center}\end{figure}

\'Etablissons par récurrence sur $n\geq 1$ les formules du lemme sur
la dimension de $F_n$ et la norme des éléments de $G_n$. Pour $n=1$,
on a $F_1=\{0\}$ et l'application $\rho_1$ est un isomorphisme, si
bien que la formule sur les normes découle du lemme
\ref{fonctiontriangulaire2}. Supposons donc $n\geq 2$ et les
formules démontrées pour $n-1$. Nous allons construire explicitement
l'inverse de $\rho_n$ en fonction de celui de $\rho_{n-1}$. Pour
tout triangle $\mathcal S$, on désigne par $F_{\mathcal S}$
l'ensemble des éléments de $E_{\mathcal S}$ nuls aux sommets de
$\mathcal S$ et par $G_{\mathcal S}$ l'orthogonal de $F_{\mathcal
S}$ dans $E_{\mathcal S}$ pour le produit scalaire de
$\ell^2(\mathcal S)$. Pour tout $(s,t,u)$ dans $\mathbb C^3_0$,
notons $\tau(s,t,u)$ l'unique fonction sur $\mathcal T_n$ qui vaut
$s$ en $a_n$, $t$ en $b_n$, $u$ en $c_n$, $\frac{t-s}{3}$ en
$a_nb_n$, $\frac{u-s}{3}$ en $a_nc_n$, $\frac{s-t}{3}$ en $b_na_n$,
$\frac{u-t}{3}$ en $b_nc_n$, $\frac{s-u}{3}$ en $c_na_n$ et
$\frac{t-u}{3}$ en $c_nb_n$ et dont la restriction à $\mathcal A_n$
(resp. $\mathcal B_n$, resp. $\mathcal C_n$) est dans $G_{\mathcal
A_n}$ (resp. $G_{\mathcal B_n}$, resp. $G_{\mathcal C_n}$). Alors,
clairement, $\tau(s,t,u)$ appartient à $E_n$ et
$\rho_n(\tau(s,t,u))=(s,t,u)$. Par ailleurs, on a, d'après le lemme
\ref{fonctiontriangulaire2} et par récurrence,
\begin{align*}\left\langle\tau(s,t,u),\psi_n\right\rangle_{{\rm L}^2\left(\bar{\Gamma},\mu\right)}
&=\frac{1}{3^n}\left(\left\langle\tau(s,t,u),\varphi^{a_nb_n,a_nc_n}_{\mathcal
A_n}\right\rangle_{\ell^2(\mathcal
A_n)}\right.\\&\qquad\left.+\left\langle\tau(s,t,u),\varphi^{b_nc_n,b_na_n}_{\mathcal
B_n}\right\rangle_{\ell^2(\mathcal
B_n)}\right.\\&\qquad\left.+\left\langle\tau(s,t,u),\varphi^{c_na_n,c_nb_n}_{\mathcal
C_n}\right\rangle_{\ell^2(\mathcal C_n)}\right)\\
&=\frac{1}{3}\frac{5^{n-2}}{9^{n-1}}\left((t-s)-(u-s)+(u-t)-(s-t)\right.\\
&\qquad\left.+(s-u)-(t-u)\right)=0\end{align*} Réciproquement, on
vérifie aisément, grâce à un calcul de produit scalaire analogue
que, si $\varphi$ est un élément de $E_n$ orthogonal à $\psi_n$ dont
la restriction à $\mathcal A_n$ (resp. $\mathcal B_n$, resp.
$\mathcal C_n$) est dans $G_{\mathcal A_n}$ (resp. $G_{\mathcal
B_n}$, resp. $G_{\mathcal C_n}$), alors $\varphi$ appartient à
l'image de $\tau$. En particulier, l'espace $G_n$ est contenu dans
l'image de $\tau$. Comme ces deux espaces sont de dimension $2$, ils
coïncident et $\tau$ est l'inverse de $\rho_{n}$. En particulier,
$F_n$ est engendré par $\psi_n$ et des éléments nuls aux sommets des
$(n-1)$-triangles, si bien que $\dim F_n=3\dim F_{n-1}+1$, d'où le
calcul de la dimension, par récurrence.  Par ailleurs, à nouveau
d'après le lemme \ref{fonctiontriangulaire2} et par récurrence, pour
tout $\varphi$ dans $G_n$, si $\rho_n(\varphi)=(s,t,u)$, on a, vue
la définition de $\tau$,
\begin{align*}\N{\varphi}_{{\rm L}^2\left(\bar{\Gamma},\mu\right)}^2&=
\frac{5^{n-2}}{9^{n-1}}\left(\abs{s}^2+\abs{t}^2+\abs{u}^2+2\abs{\frac{s-t}{3}}^2
+2\abs{\frac{t-u}{3}}^2+2\abs{\frac{s-u}{3}}^2\right)\\
&=\frac{1}{3}\frac{5^{n-1}}{9^{n-1}}\left(\abs{s}^2+\abs{t}^2+\abs{u}^2\right)\end{align*}
(sans oublier, pour la dernière égalité, que $s+t+u=0$). La formule
concernant les normes en découle par récurrence.

Enfin, pour $n\geq 2$, donnons-nous $\varphi$ dans $G_n$ et posons
$\psi=\mathbb E(\varphi|\theta_{n-1})$. Comme dans la démonstration
du lemme \ref{espacep32comp}, on déduit du lemme
\ref{fonctiontriangulaire3} et du fait que $\theta_{n-1}$ induit des
isomorphismes entre les $(n-1)$-triangles de $\mathcal T_n$ et
$\mathcal T_{n-1}$ que, comme $\varphi$ appartient à $E_n$, $\psi$
appartient à $E_{n-1}$. Comme les éléments de $F_{n-1}$ sont nuls
aux sommets des $(n-1)$-triangles, ils appartiennent aussi à $F_n$,
donc ils sont orthogonaux à $\varphi$, si bien que $\psi$ appartient
à $G_{n-1}$. D'après les lemmes \ref{facteurtriangle} et
\ref{fonctiontriangulaire3}, on a
$\psi(a_n)=\frac{1}{3}(\varphi(a_n)+\varphi(b_na_n)+\varphi(c_na_n))$
et, donc, d'après les formules ci-dessus, si
$\rho_n(\varphi)=(s,t,u)$, on a
$\psi(a_n)=\frac{1}{3}(s+\frac{s-t}{3}+\frac{s-u}{3})=\frac{2}{3}s$
et $\rho_{n-1}(\psi)=\frac{2}{3}\rho_n(\varphi)$.
\end{demo}

Nous pouvons alors décrire l'espace propre de valeur propre $-2$ de
$\bar{\Delta}$ :

\begin{Lem} \label{espacep-21comp}
L'espace propre associé à
la valeur propre $-2$ est de dimension infinie et engendré par des
fonctions triangulaires nulles aux sommets de leur triangle de
définition.\end{Lem}

\begin{demo}
Comme, d'après le lemme \ref{espacep-21comp}, pour tout $n\geq 1$,
l'espace $F_n$ est de dimension $\frac{1}{2}(3^{n-1}-1)$ et que,
d'après le lemme \ref{fonctiontriangulaire4}, ses éléments sont des
fonctions propres de valeur propre $-2$, l'espace propre associé à
la valeur propre $-2$ est de dimension infinie.

Soit $\varphi$ une fonction propre de valeur propre $-2$ dans ${\rm
L}^2\left(\bar{\Gamma},\mu\right)$ qui soit orthogonale à toutes les
fonctions propres triangulaires nulles aux sommets de leur triangle
de définition. Montrons que $\varphi$ est nulle. Pour tout entier
$n\geq 1$, soit $\varphi_n=\mathbb E(\varphi|\theta_n)$. D'après le
corollaire \ref{spectrecomp}, on a $\bar{\Pi}\varphi=0$ et
$\varphi\circ\alpha=-\varphi$ et, donc, d'après les lemmes
\ref{fonctiontriangulaire3} et \ref{involution}, pour tout $n\geq
1$, $\varphi_n$ appartient à $E_n$. Comme $\varphi$ est orthogonale
aux éléments de $F_n$, $\varphi_n$ appartient à $G_n$. Si $n\geq 2$,
comme $\varphi_{n-1}=\mathbb E(\varphi_n|\theta_{n-1})$, d'après le
lemme \ref{espacep-21comp}, on a
$\rho_{n-1}(\varphi_{n-1})=\frac{2}{3}\rho_n(\varphi_n)$. Il existe
donc $v$ dans $\mathbb C^3_0$ tel que, pour tout $n\geq 1$, on ait
$\rho_n(\varphi_n)=\left(\frac{3}{2}\right)^{n-1}v$, si bien que,
toujours d'après le lemme \ref{espacep-21comp},
$$\N{\varphi_n}^2_{{\rm L}^2\left(\bar{\Gamma},\mu\right)}
=\left(\frac{5}{9}\right)^{n-1}\N{\rho_n(\varphi_n)}_0^2
=\left(\frac{5}{4}\right)^{n-1}\N{v}_0^2.$$ Comme, d'après le lemme
\ref{fonctiontriangulaire3}, on a $\varphi_n\td{n}{\infty}\varphi$
dans ${\rm L}^2\left(\bar{\Gamma},\mu\right)$, il vient
nécessairement $v=0$, donc, pour tout $n\geq 1$, $\varphi_n=0$ et
$\varphi=0$, ce qu'il fallait démontrer.
\end{demo}

Des lemmes \ref{transformevalpcomp} et \ref{espacep-21comp}, on
déduit par récurrence le

\begin{Cor}\label{espacep-22comp}
Pour tout $x$ dans $\bigcup_{n\in\mathbb N}f^{-n}(-2)$, l'espace
propre associé à la valeur propre $x$ est de dimension infinie et
engendré par des fonctions triangulaires nulles aux sommets de leur
triangle de définition.\end{Cor}

\section{Décomposition spectrale de ${\rm L}^2\left(\bar{\Gamma},\mu\right)$}
\label{decompocomp}

Dans ce paragraphe, nous allons montrer que ${\rm
L}^2\left(\bar{\Gamma},\mu\right)$ est la somme directe orthogonale
de l'espace des fonctions constantes, des espaces propres associés
aux éléments de l'ensemble $\bigcup_{n\in\mathbb
N}f^{-n}(-2)\cup\bigcup_{n\in\mathbb N}f^{-n}(0)$ et des espaces
cycliques engendrés par les fonctions $1$-triangulaires $\varphi$
telles que $\bar{\Pi}\varphi=0$. Comme à la section
\ref{valeurpcomp}, on note $E_1$ l'espace de ces fonctions. Commen\c
cons par décrire leurs espaces cycliques.

\begin{Lem} \label{relation7}
Soit $\varphi$ dans $E_1$. On a
$\left(\bar{\Delta}+2\right)\varphi=\left(\bar{\Delta}-1\right)\bar{\Pi}^*\varphi$
et
$\bar{\Pi}\bar{\Delta}\varphi=\left(1+\frac{1}{3}\bar{\Delta}\right)\varphi$.
\end{Lem}

\begin{demo} Soit $(s,t,u)=(\varphi(a_1),\varphi(b_1),\varphi(c_1))$. On a,
par définition, $s+t+u=0$. Soit $p$ dans $\bar{\Gamma}$. Quitte à
faire agir le groupe $\mathfrak S$, on peut supposer que
$\bar{\Pi}p$ appartient au bréchet $B_0$ de la section
\ref{compactifie}. Alors, les valeurs de $\varphi$ et de
$\bar{\Delta}\varphi$ sur le $1$-triangle contenant $p$ et sur ses
voisins sont décrites par la figure \ref{valeurs1}.
\begin{figure}\begin{center}\input{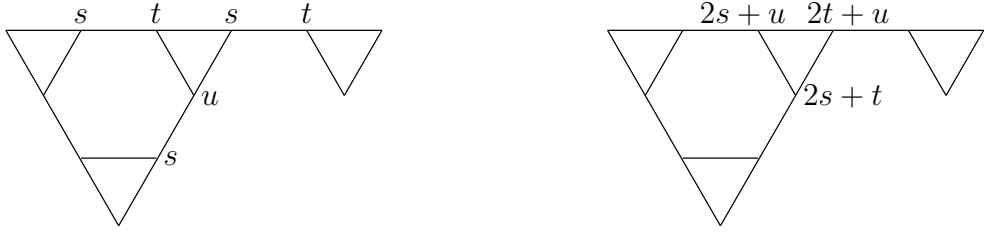}
\caption{Valeurs de $\varphi$ et de
$\bar{\Delta}\varphi$}\label{valeurs1}\end{center}\end{figure} De
même, les valeurs de $\bar{\Pi}^*\varphi$ et
$\bar{\Delta}\bar{\Pi}^*\varphi$ sur le $1$-triangle contenant $p$
et sur ses voisins sont représentées par la figure \ref{valeurs2}.
\begin{figure}\begin{center}\input{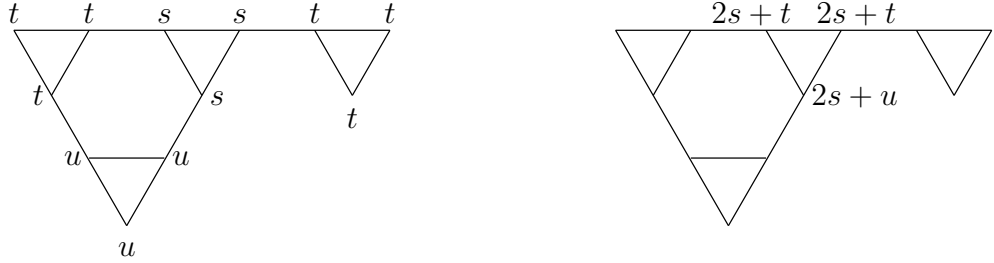}
\caption{Valeurs de $\bar{\Pi}^*\varphi$ et de
$\bar{\Delta}\bar{\Pi}^*\varphi$}\label{valeurs2}\end{center}\end{figure}
Si $\theta_1(p)=a_1$ ou $\theta_1(p)=b_1$, on a donc
$\left(\bar{\Delta}+2\right)\varphi(p)=2s+2t+u=s+t
=\left(\bar{\Delta}-1\right)\bar{\Pi}^*\varphi(p)$ ; si
$\theta_1(p)=c_1$, on a
$\left(\bar{\Delta}+2\right)\varphi(p)=2s+t+2u=s+u
=\left(\bar{\Delta}-1\right)\bar{\Pi}^*\varphi(p)$. Il vient bien
$\left(\bar{\Delta}+2\right)\varphi=\left(\bar{\Delta}-1\right)\bar{\Pi}^*\varphi$.

Par définition, on a $\bar{\Pi}\varphi=0$, si bien que, en
appliquant $\bar{\Pi}$ à l'identité précédente, il vient
$\bar{\Pi}\bar{\Delta}\varphi=\bar{\Pi}\left(\bar{\Delta}+2\right)\varphi
=\bar{\Pi}\bar{\Delta}\bar{\Pi}^*\varphi-\varphi=
\left(1+\frac{1}{3}\bar{\Delta}\right)\varphi$, où, pour la dernière
égalité, on a utilisé le lemme \ref{relation6}.
\end{demo}

Grâce au lemme \ref{relation7}, nous allons procéder comme dans la section \ref{decompo}
pour déterminer les mesures spectrales des éléments de $E_1$.
Commen\c cons par montrer que ces mesures ne chargent pas les points $-2$ et $0$.

\begin{Lem} \label{produitnul1comp}
Soient $\varphi$ dans $E_1$ et $\psi$ une fonction propre de valeur $-2$ ou $0$ dans
${\rm L}^2\left(\bar{\Gamma},\mu\right)$. On a $\langle\varphi,\psi\rangle=0$.
\end{Lem}

\begin{demo}
Supposons que $\psi$ est un vecteur propre de valeur propre $0$.
D'après le corollaire \ref{spectrecomp}, on a $\bar{\Pi}\psi=0$ et
$\psi\circ\alpha=\psi$ et, d'après le corollaire
\ref{espacep03comp}, on peut supposer que, pour un certain entier
$n\geq 2$, $\psi$ est $n$-triangulaire, de valeur nulle aux sommets
de $\mathcal T_n$. Soient $p$, $q$ et $r$ les sommets d'un
$2$-triangle $\mathcal S$ de $\mathcal T_n$ et soient $pq$, $qp$,
$pr$, $rp$, $qr$ et $rq$ les autres points de $\mathcal S$, avec la
convention du corollaire \ref{decoupetriangles}. Alors, on a
$\psi(qp)=\psi(pq)$ et $\psi(rp)=\psi(pr)$, donc
$\psi(p)+\psi(qp)+\psi(rp)=0$ et, en utilisant les identités
analogues sur les autres $1$-triangles de $\mathcal S$, d'après le
lemme \ref{facteurtriangle}, comme $\varphi$ est $1$-triangulaire,
on a
\begin{multline*}\sum_{s\in\mathcal S}\overline{\varphi(s)}\psi(s)=
\overline{\varphi(p)}\left(\psi(p)+\psi(qp)+\psi(rp)\right)\\
+\overline{\varphi(q)}\left(\psi(q)+\psi(pq)+\psi(rq)\right)
+\overline{\varphi(r)}\left(\psi(r)+\psi(pr)+\psi(qr)\right)=0\end{multline*}
et, donc, d'après le lemme \ref{fonctiontriangulaire2},
$\langle\varphi,\psi\rangle=0$.

Intéressons-nous à présent au cas de la valeur propre $-2$. Pour
tout $n\geq 1$, notons $E_n$ et $F_n$ comme à la section
\ref{valeurpcomp}. Soit
$(s,t,u)=(\varphi(a_1),\varphi(b_1),\varphi(c_1))$. Montrons par
récurrence sur $n$ que, si $\psi$ appartient à $E_n$, on a
$$\sum_{p\in\mathcal T_n}\varphi(p)\psi(p)=2^{n-1}(s\psi(a_n)+t\psi(b_n)+u\psi(c_n)).$$
Pour $n=1$, le résultat est trivial. Supposons $n\geq 2$ et le
résultat établi pour $n$. Alors, en appliquant la récurrence à la
restriction de $\psi$ au $(n-1)$-triangles de $\mathcal T_n$, on
obtient, comme $\varphi$ est $1$-triangulaire, d'après le lemme
\ref{facteurtriangle},
\begin{multline*}\sum_{p\in\mathcal T_n}\varphi(p)\psi(p)
=2^{n-2}(\psi(a_n)s+\psi(a_nb_n)t+\psi(a_nc_n)u\\
+\psi(b_n)t+\psi(b_na_n)s+\psi(b_nc_n)u+\psi(c_n)u+\psi(c_na_n)s+\psi(c_nb_n)t).\end{multline*}
Or, on a $\psi(a_nb_n)+\psi(b_na_n)=0$ et, d'après le lemme
\ref{sommenulle}, $\psi(a_n)+\psi(a_nb_n)+\psi(a_nc_n)=0$, si bien
que $\psi(b_na_n)+\psi(c_na_n)=\psi(a_n)$. En utilisant cette
identité et les formules analogues aux autres sommets de $\mathcal
T_n$, il vient
$$\sum_{p\in\mathcal T_n}\varphi(p)\psi(p)
=2^{n-1}(s\psi(a_n)+t\psi(b_n)+u\psi(c_n)),$$ ce qu'il fallait
démontrer. En particulier, pour $\psi$ dans $F_n$, on a, d'après le
lemme \ref{fonctiontriangulaire2}, $\langle\varphi,\psi\rangle=0$
et, donc, d'après le lemme \ref{espacep-21comp}, ceci est encore
vrai pour tout vecteur propre $\psi$ de valeur propre $-2$.
\end{demo}

\begin{Cor}  \label{produitnul2comp} Soient $\varphi$ dans $E_1$ et $\psi$
une fonction propre associée à une valeur propre dans
$\bigcup_{n\in\mathbb N}f^{-n}(-2)\cup\bigcup_{n\in\mathbb
N}f^{-n}(0)$. On a $\langle\varphi,\psi\rangle=0$.\end{Cor}

\begin{demo} D'après le corollaire \ref{spectrecomp} et les lemmes
\ref{transformevalpcomp} et \ref{produitnul1comp}, il suffit de
montrer que, pour $x$ dans $\mathbb R$, si $\psi$ est un vecteur
propre de valeur propre $x$ et si $\langle \varphi,\psi\rangle=0$,
on a
$\left\langle\varphi,\bar{\Pi}^*\psi\right\rangle=\left\langle\varphi,\bar{\Delta}\bar{\Pi}^*\psi\right\rangle=0$.
Or, d'une part, par définition, on a $\bar{\Pi}\varphi=0$, donc
$\left\langle\varphi,\bar{\Pi}^*\psi\right\rangle=0$. D'autre part,
d'après le lemme \ref{relation7}, on a
$$\left\langle\varphi,\bar{\Delta}\bar{\Pi}^*\psi\right\rangle
=\left\langle\bar{\Pi}\bar{\Delta}\varphi,\psi\right\rangle
=\left\langle\varphi,\left(1+\frac{1}{3}\bar{\Delta}\right)\psi\right\rangle
=\left(1+\frac{x}{3}\right)\langle\varphi,\psi\rangle=0,$$ ce qu'il
fallait démontrer.
\end{demo}

Posons, pour tout $x\neq -3$, $j(x)=\frac{1}{3}\frac{3-x}{x+3}$ et,
pour $x\neq\frac{1}{2}$,
$\zeta(x)=\frac{1}{3}\frac{(x+3)(x-1)}{2x-1}$. Comme pour le
corollaire \ref{spectrediffus1}, nous déduisons du lemme
\ref{relation7} et du corollaire \ref{mesurespectralecomp} le

\begin{Cor} \label{spectrediffus1comp}
Soit $\nu_\zeta$ l'unique probabilité borélienne sur $\Lambda$ telle
qu'on ait $L_\zeta^*\nu_\zeta=\nu_\zeta$. Pour tout $\varphi$ dans
$E_1$, la mesure spectrale de $\varphi$ est
$\N{\varphi}_2^2j\nu_\zeta$.
\end{Cor}

\begin{demo} Comme la démonstration de ce résultat est
 analogue à celle du corollaire \ref{spectrediffus1},
nous en reprenons seulement les grandes lignes. Soit $\lambda$ la
mesure spectrale de $\varphi$. D'après le lemme
\ref{produitnul1comp}, on a $\lambda(-2)=0$. Posons, pour
$x\notin\left\{-2,\frac{1}{2}\right\}$,
$\theta(x)=\frac{x(x-1)^2}{3(x+2)(2x-1)}$. On a
$\theta=\frac{j}{j\circ f}\zeta$ et, d'après le corollaire
\ref{mesurespectralecomp} et le lemme \ref{relation7},
$\lambda=L_\theta^*\lambda$. D'après le lemme
\ref{fonctiontriangulaire2}, $\varphi$ est orthogonale aux fonctions
constantes. Par conséquent, d'après le lemme \ref{espacep34comp}, on
a $\lambda(3)=0$. De plus, d'après les corollaires \ref{spectrecomp}
et \ref{produitnul2comp}, la mesure $\lambda$ est concentrée sur
$\Lambda$.

La fonction $\zeta$ est strictement positive sur $\Lambda$ et
$L_\zeta(1)=1$. D'après le lemme \ref{transfert}, il existe une
unique probabilité borélienne $\nu_\zeta$ sur $\Lambda$ telle que
$L_\zeta^*\nu_\zeta=\nu_\zeta$. En raisonnant comme dans la
démonstration du corollaire \ref{spectrediffus1}, on montre que les
mesures $\lambda$ et $j\nu_\zeta$ sont proportionnelles. Comme on a
$L_{\zeta}j=1$, il vient
$\lambda=\N{\varphi}_2^2j\nu_\zeta$.\end{demo}

Notons toujours $l$ la fonction $x\mapsto x$ sur $\Lambda$ et
posons, pour $x\neq 1$, $m(x)=\frac{x+2}{x-1}$ et, pour $x\neq
\frac{1}{2}$, $\xi(x)=\frac{1}{3}\frac{x(x-1)}{2x-1}$. Notons
$\bar{\Phi}$ le sous-espace fermé de ${\rm
L}^2\left(\bar{\Gamma},\mu\right)$ engendré par les éléments de
$E_1$ et par leurs images par les puissances de $\bar{\Delta}$ et,
comme à la section \ref{valeurpcomp}, désignons par $\rho_1$
l'isomorphisme $\mathfrak S$-équivariant de $E_1$ dans $\mathbb
C^3_0$. Munissons toujours $\mathbb C^3_0$ de la norme hermitienne
$\N{.}_0$ égale à un tiers de la norme canonique et notons
$\langle.,.\rangle_0$ le produit scalaire associé. D'après le lemme
\ref{fonctiontriangulaire2}, l'application $\rho_1$ est une
isométrie. Identifions les espaces de Hilbert ${\rm
L}^2\left(j\nu_\zeta,\mathbb C^3_0\right)$ et ${\rm
L}^2\left(j\nu_\zeta\right)\otimes\mathbb C^3_0$ et, pour tout
polynôme $p$ dans $\mathbb C[X]$ et pour tout $v$ dans $\mathbb
C^3_0$, posons $\widehat{p\otimes
v}=p\left(\bar{\Delta}\right)\rho_1^{-1}(v)$. Nous avons un analogue
de la proposition \ref{spectrediffus2} :

\begin{Prop}\label{spectrediffus2comp}
L'application $g\mapsto\hat{g}$ induit une isométrie $\mathfrak
S$-équivariante de ${\rm L}^2\left(j\nu_\zeta,\mathbb C^3_0\right)$
dans $\bar{\Phi}$. Le sous-espace $\bar{\Phi}$ est stable par les
opérateurs $\bar{\Delta}$, $\bar{\Pi}$ et $\bar{\Pi}^*$. Pour tout
$g$ dans ${\rm L}^2\left(j\nu_\zeta,\mathbb C^3_0\right)$, on a
\begin{align*}\bar{\Delta}\hat{g}&=\widehat{lg}\\
\bar{\Pi}\hat{g}&=\widehat{L_\xi g}\\
\bar{\Pi}^*\hat{g}&=\widehat{m(g\circ f)}.\end{align*}
\end{Prop}

\begin{demo}
Soit $p$ dans $\mathbb C[X]$. L'application
\begin{align*}\mathbb C^3_0\times\mathbb C^3_0&\rightarrow\mathbb C\\
(v,w)&\mapsto\left\langle p\left(\bar{\Delta}\right)\rho_1^{-1}(v),
\rho_1^{-1}(w)\right\rangle_{{\rm
L}^2\left(\bar{\Gamma},\mu\right)}\end{align*} est une forme
sesquilinéaire $\mathfrak S$-invariante. Comme la représentation de
$\mathfrak S$ dans $\mathbb C^3_0$ est irréductible, cette forme
sesquilinéaire est proportionnelle au produit scalaire $\langle
.,.\rangle_0$. D'après le lemme \ref{fonctiontriangulaire2} et le
corollaire \ref{spectrediffus1comp}, pour $v$ dans $\mathbb C^3_0$,
on a
$$\left\langle p\left(\bar{\Delta}\right)\rho_1^{-1}(v),
\rho_1^{-1}(v)\right\rangle_{{\rm
L}^2\left(\bar{\Gamma},\mu\right)}=\N{v}^2_0\int_{\bar{\Gamma}}pj\de\nu_\zeta,$$
par conséquent, pour tous $p$ et $q$ dans $\mathbb C[X]$, pour tous
$v$ et $w$ dans $\mathbb C^3_0$, on a $$\left\langle
\widehat{p\otimes v},\widehat{q\otimes w}\right\rangle_{{\rm
L}^2\left(\bar{\Gamma},\mu\right)}=\langle v,w\rangle_0 \langle
p,q\rangle_{{\rm L}^2\left(j\nu\zeta\right)}$$ et, donc,
l'application $g\mapsto\hat{g}$ induit une isométrie de ${\rm
L}^2\left(j\nu_\zeta,\mathbb C^3_0\right)$ dans un sous-espace fermé
de ${\rm L}^2\left(\bar{\Gamma},\mu\right)$. Comme ce sous-espace
est engendré par les éléments de $E_1$ et leurs images par les
puissances de $\bar{\Delta}$, il est, par définition, égal à
$\bar{\Phi}$.

La suite de la démonstration est analogue à celle de la proposition
\ref{spectrediffus2}.

La stabilité de $\bar{\Phi}$ par $\bar{\Delta}$ et la formule pour
$\bar{\Delta}$ résultent de la définition même des objets concernés.

Un calcul direct montre que $L_\xi(1)=0$ et que $L_\xi
(l)=1+\frac{1}{3}l$, si bien que, pour tout entier naturel $n$, on a
$L_\xi(f^n)=0$ et $L_\xi(f^nl)=l^n\left(1+\frac{1}{3}l\right)$. Or,
d'après les lemmes \ref{relation5} et \ref{relation7}, pour tout
$\varphi$ dans $E_1$, on a
$\bar{\Pi}\left(f\left(\bar{\Delta}\right)^n\varphi\right)=0$ et
$\bar{\Pi}\left(f\left(\bar{\Delta}\right)^n\bar{\Delta}\varphi\right)
=\bar{\Delta}^n\left(1+\frac{1}{3}\bar{\Delta}\right)\varphi$.
L'espace $\bar{\Phi}$ est donc stable par $\bar{\Pi}$ et, pour tous
$p$ dans $\mathbb C[X]$ et $v$ dans $\mathbb C^3_0$, on a
$\bar{\Pi}\widehat{p\otimes v}=\widehat{L_\xi(p)\otimes v}$. Comme
$\zeta$ est partout $>0$ sur $\Lambda$, il existe un réel $c>0$ tel
que, pour tout $x$ de $\Lambda$, on ait $\abs{\xi(x)}\leq
c\zeta(x)$, si bien que, pour toute fonction borélienne $g$ sur
$\Lambda$, on a $\abs{L_\xi(g)}\leq cL_\zeta\left(\abs{g}\right)$.
En raisonnant comme dans la démonstration de la proposition
\ref{spectrediffus2}, on montre que $L_\zeta$ est borné dans ${\rm
L}^2(j\nu_\zeta)$. On en déduit que $L_\xi$ est borné et l'identité
concernant $\bar{\Pi}$ en découle, par densité.

Enfin, d'après les lemmes \ref{relation5} et \ref{relation7}, pour
tous $p$ dans $\mathbb C[X]$ et $\varphi$ dans $E_1$, on a
$\left(\bar{\Delta}-1\right)\bar{\Pi}^*\left(p\left(\bar{\Delta}\right)\varphi\right)
=p\left(f\left(\bar{\Delta}\right)\right)\left(\bar{\Delta}+2\right)\varphi$.
D'après le corollaire \ref{spectrecomp}, $1$ n'appartient pas au
spectre de $\bar{\Delta}$, si bien que, par densité, pour toute
fonction rationnelle $p$ dont les pôles n'appartiennent pas au
spectre de $\bar{\Delta}$, on a
$\bar{\Pi}^*\left(p\left(\bar{\Delta}\right)\varphi\right)
=(m(p\circ f))\left(\bar{\Delta}\right)\varphi$ et, donc, l'espace
$\bar{\Phi}$ est stable par $\bar{\Pi}^*$. De plus, comme, pour tout
$x$ dans $\Lambda$, on a
$m(x)^2\frac{j(x)}{j(f(x))}=\frac{x(x+2)}{(x-1)(x+3)}$, il vient,
par un calcul élémentaire, $L_\zeta\left(m^2\frac{j}{j\circ
f}\right)=1$ et, pour tout $g$ dans ${\rm L}^2(j\nu_\zeta)$,
$$\int_\Lambda
\abs{m (g\circ f)}^2j\de\nu_\zeta =\int_\Lambda
\left(m^2\frac{j}{j\circ f}\right)\abs{g\circ f}^2(j\circ
f)\de\nu_\zeta=\int_\Lambda\abs{g}^2j\de\nu_\zeta.$$ La formule pour
$\bar{\Pi}^*$ en découle, par densité.
\end{demo}

Intéressons-nous à présent aux autres composantes $\mathfrak
S$-isotypiques de l'espace ${\rm L}^2\left(\bar{\Gamma},\mu\right)$.
Notons $\varepsilon:\mathfrak S\rightarrow\{-1,1\}$ le morphisme de
signature. Nous dirons qu'une fonction $\varphi$ sur $\bar{\Gamma}$
est $(\mathfrak S,\varepsilon)$-semi-invariante si, pour tout $s$
dans $\mathfrak S$, on a $\varphi\circ s=\varepsilon(s)\varphi$.
Notons toujours $k$ et $l$ les fonctions $x\mapsto x+2$ et $x\mapsto
x$.

\begin{Prop} \label{spectrediscret1comp}
Pour tout entier $n\geq 1$, l'espace des fonctions
$n$-trian\-gu\-laires $\mathfrak S$-invariantes sur $\bar{\Gamma}$
est stable par $\bar{\Delta}$ et le polynôme caractéristique de
$\bar{\Delta}$ y est
$$(X-3)\prod_{p=0}^{n-2}(l\circ f^p(X))^{\frac{3^{n-2-p}+2n-2p-1}{4}}(k\circ f^p(X))^{\frac{3^{n-2-p}-2n+2p+3}{4}}.$$
Pour tout entier $n\geq 2$, l'espace des fonctions $n$-triangulaires
$(\mathfrak S,\varepsilon)$-semi-invariantes sur $\bar{\Gamma}$ est
stable par $\bar{\Delta}$ et le polynôme caractéristique de
$\bar{\Delta}$ y est
$$\prod_{p=0}^{n-2}(l\circ f^p(X))^{\frac{3^{n-2-p}-2n+2p+3}{4}}(k\circ f^p(X))^{\frac{3^{n-2-p}+2n-2p-1}{4}}.$$
\end{Prop}

\begin{demo} Ces espaces sont stables d'après le lemme \ref{fonctiontriangulaire4}.
Le calcul des polynômes caractéristiques s'obtient en raisonnant
comme dans la démons\-tration de la proposition
\ref{spectrePascalfinifinal}.\end{demo}

De cette proposition on déduit, grâce au lemme
\ref{fonctiontriangulaire3}, le

\begin{Cor}\label{spectrediscret2comp}
Le spectre de $\bar{\Delta}$ dans l'espace des éléments $\mathfrak
S$-invariants de ${\rm L}^2\left(\bar{\Gamma},\mu\right)$ est
discret. Les valeurs propres de $\bar{\Delta}$ y sont $3$, qui est
simple, et les éléments de $\bigcup_{n\in\mathbb
N}f^{-n}(-2)\cup\bigcup_{n\in\mathbb N}f^{-n}(0)$. Le spectre de
$\bar{\Delta}$ dans l'espace des éléments $(\mathfrak
S,\varepsilon)$-semi-invariants de ${\rm
L}^2\left(\bar{\Gamma},\mu\right)$ est discret. Les valeurs propres
de $\bar{\Delta}$ y sont les éléments de $\bigcup_{n\in\mathbb
N}f^{-n}(-2)\cup\bigcup_{n\in\mathbb N}f^{-n}(0)$.
\end{Cor}

La démonstration du théorème \ref{spectrePascalcomp} s'achève avec
la

\begin{Prop} \label{spectrePascalfinalcomp}
Soit $\bar{\Phi}^\perp$ l'orthogonal de $\bar{\Phi}$ dans ${\rm
L}^2\left(\bar{\Gamma},\mu\right)$. Le spectre de $\bar{\Delta}$
dans  $\bar{\Phi}^\perp$ est discret. Ses valeurs propres y sont
$3$, qui est simple, et les éléments de $\bigcup_{n\in\mathbb
N}f^{-n}(-2)\cup\bigcup_{n\in\mathbb N}f^{-n}(0)$.
\end{Prop}

La démonstration de cette proposition est analogue à celle de la
proposition \ref{spectrePascalfinal}. Elle exige que nous
introduisions des objets qui joueront un rôle semblable à celui des
espaces $L_n$, $n\in\mathbb N$, de cette démonstration.

Reprenons les notations de la section \ref{compactifie} et rappelons
que, par construction, si $p$ est un point de $\bar{\Gamma}$ tel que
$\theta_1(p)=a_1$, le bréchet de $p$ est $B_0$ ou $riB_0$. Pour tout
entier $n\geq 1$, notons $\mathcal B_n$ l'ensemble constitué de la
réunion de $\mathcal T_n-\partial\mathcal T_n$ et de l'ensemble des
six couples de la forme $(d,B)$ où $d$ est un élément de
$\partial\mathcal T_n$ et $B$ un des deux bréchets pour lesquels il
existe des points $p$ de $B$ tels que $\theta_n(p)=d$. Notons
$\tau_n$ l'application localement constante
$\bar{\Gamma}\rightarrow\mathcal B_n$ telle que, pour tout $p$ dans
$\bar{\Gamma}$, si $p$ n'est pas le sommet d'un $n$-triangle, on a
$\tau_n(p)=\theta_n(p)$ et, si $p$ est le sommet d'un $n$-triangle,
$\tau_n(p)$ est le couple formé de $\theta_n(p)$ et du bréchet qui
contient $p$. Enfin, disons qu'une fonction $\varphi$ sur
$\bar{\Gamma}$ est $\tau_n$-mesurable si on a
$\varphi=\psi\circ\tau_n$, où $\psi$ est une fonction définie sur
$\mathcal B_n$. L'intérêt de cette définition provient du

\begin{Lem} \label{recurrencecomp}
Soient $n$ un entier $\geq 1$ et $\varphi$ une fonction $\tau_{n+1}$-mesurable sur $\bar{\Gamma}$.
Alors les fonctions $\bar{\Pi}\varphi$ et $\bar{\Pi}\bar{\Delta}\varphi$ sont $\tau_n$-mesurables.\end{Lem}

\begin{demo} Soit $p$ un point de $\bar{\Gamma}$.
Si $p$ n'est pas le sommet d'un $n$-triangle, le triangle
$\bar{\Pi}^{-1}p$ ne contient pas de sommet d'un $(n+1)$-triangle.
De même, aucun des voisins des points de $\bar{\Pi}^{-1}p$ n'est le
sommet d'un $(n+1)$-triangle. Pour tous les points $q$ en question
dans le calcul de $\bar{\Pi}\varphi(p)$ et de
$\bar{\Pi}\bar{\Delta}\varphi(p)$, on a donc
$\tau_n(q)=\theta_n(q)$. Par conséquent, par définition de
$\theta_n$ et d'après le lemme \ref{facteurtriangle},
$\bar{\Pi}\varphi(p)$ et $\bar{\Pi}\bar{\Delta}\varphi(p)$ ne
dépendent que de $\theta_n(p)$.

Si, à présent, $p$ est le sommet d'un $n$-triangle, le voisin $q$ de
$p$ qui n'appartient pas à ce $n$-triangle est lui-même le sommet
d'un $n$-triangle, et, d'après le lemme \ref{brechetvoisin}, le
bréchet de $q$ est déterminé par le bréchet de $p$. En particulier,
$\tau_n(q)$ est déterminé par $\tau_n(p)$. Un seul des trois
antécédents du point $p$ par l'application $\bar{\Pi}$ est le sommet
d'un $(n+1)$-triangle. D'après le lemme \ref{brechet1}, il s'agit de
celui dont le bréchet est égal à celui de $p$. En particulier,
l'image par $\tau_{n+1}$ de ce point $r$ est déterminée par
$\tau_n(p)$. De même, l'image par $\tau_{n+1}$ de l'unique
antécédent $s$ de $q$ qui est le sommet d'un $(n+1)$-triangle ne
dépend que de $\tau_n(q)$, et donc de $\tau_n(p)$. Le point $s$ est
le voisin de $r$ qui n'appartient pas à $\bar{\Pi}^{-1}p$. Enfin,
les deux autres points de $\bar{\Pi}^{-1}p$ et leurs voisins qui
n'appartiennent pas à $\bar{\Pi}^{-1}p$ ne sont pas des sommets d'un
$(n+1)$-triangle et, donc, leur image par $\tau_{n+1}$ est leur
image par $\theta_{n+1}$ qui ne dépend que de $\theta_n(p)$. \`A
nouveau, $\bar{\Pi}\varphi(p)$ et $\bar{\Pi}\bar{\Delta}\varphi(p)$
ne dépendent que de $\tau_n(p)$.
\end{demo}

\begin{demo}[Démonstration de la proposition \ref{spectrePascalfinalcomp}]
D'après le lemme \ref{espacep34comp}, la valeur propre $3$ de
$\bar{\Delta}$ est simple. D'après les corollaires
\ref{espacep03comp} et \ref{espacep-22comp}, les espaces propres
associés aux éléments de $\bigcup_{n\in\mathbb
N}f^{-n}(-2)\cup\bigcup_{n\in\mathbb N}f^{-n}(0)$ sont non-triviaux.
Notons $\bar{P}$ le projecteur orthogonal de ${\rm
L}^2\left(\bar{\Gamma},\mu\right)$ dans $\bar{\Phi}^\perp$ et, pour
tous $\varphi$ et $\psi$ dans
 ${\rm L}^2\left(\bar{\Gamma},\mu\right)$, notons $\lambda_{\varphi,\psi}$
l'unique mesure complexe borélienne sur $\mathbb R$ telle que, pour
tout polynôme $p$ dans $\mathbb C[X]$, on ait $\int_{\mathbb R}p
\de\lambda_{\varphi,\psi}=\left\langle
p\left(\bar{\Delta}\right)\varphi,\psi\right\rangle$. D'après la
proposition \ref{spectrediffus2comp}, l'opérateur $\bar{P}$ commute
à $\bar{\Delta}$, $\bar{\Pi}$ et $\bar{\Pi}^*$. D'après le lemme
\ref{fonctiontriangulaire3}, pour démontrer la proposition, il
suffit d'établir que, pour tout entier $n\geq 1$, pour toute
fonction $\tau_n$-mesurable $\varphi$, pour tout $\psi$ dans ${\rm
L}^2\left(\bar{\Gamma},\mu\right)$, la mesure
$\lambda_{\bar{P}\varphi,\psi}$ est atomique et concentrée sur
l'ensemble $\bigcup_{n\in\mathbb N}f^{-n}(3)\cup\bigcup_{n\in\mathbb
N}f^{-n}(0)$. Montrons ce résultat par récurrence sur $n$.

Pour $n=1$, les fonctions $\tau_1$-mesurables sont les fonctions qui
ne dé\-pendent que du bréchet. On vérifie aisément que cet espace
est engendré par les fonctions constantes, une droite de fonctions
$(\mathfrak S,\varepsilon)$-semi-invariantes, les éléments de $E_1$
et leurs images par $\bar{\Delta}$. Dans ce cas, la description des
mesures spectrales découle immédiatement des corollaires
\ref{spectrediffus1comp} et \ref{spectrediscret2comp}.

Si le résultat est vrai pour un entier $n\geq 1$, donnons-nous une
fonction $\tau_{n+1}$-mesurable $\varphi$. Alors, d'après le lemme
\ref{recurrencecomp}, les fonctions $\bar{\Pi}\varphi$ et
$\bar{\Pi}\bar{\Delta}\varphi$ sont $\tau_n$-mesurables et, par
récurrence, pour tout $\psi$ dans  ${\rm
L}^2\left(\bar{\Gamma},\mu\right)$, les mesures
$\lambda_{\bar{\Pi}\bar{P}\varphi,\psi}
=\lambda_{\bar{P}\bar{\Pi}\varphi,\psi}$ et
$\lambda_{\bar{\Pi}\bar{\Delta}\bar{P}\varphi,\psi}
=\lambda_{\bar{P}\bar{\Pi}\bar{\Delta}\varphi,\psi}$ sont atomiques
et concentrées sur l'ensemble $\bigcup_{n\in\mathbb
N}f^{-n}(3)\cup\bigcup_{n\in\mathbb N}f^{-n}(0)$. En raisonnant
comme dans le lemme \ref{poussemesure}, on en déduit que les mesures
$\lambda_{\bar{P}\varphi,\bar{\Pi}^*\psi}$ et
$\lambda_{\bar{\Delta}\bar{P}\varphi,\bar{\Pi}^*\psi}
=\lambda_{\bar{P}\varphi,\bar{\Delta}\bar{\Pi}^*\psi}$ sont
atomiques et concentrées sur l'ensemble $\bigcup_{n\in\mathbb
N}f^{-n}(3)\cup\bigcup_{n\geq 1}f^{-n}(0)$. Or, d'après le
corollaire \ref{spectrecomp}, le spectre de $\bar{\Delta}$ dans
l'orthogonal du sous-espace de
 ${\rm L}^2\left(\bar{\Gamma},\mu\right)$ engendré par l'image de $\bar{\Pi}^*$ et par celle de $\bar{\Delta}\bar{\Pi}^*$ est égal à $\{-2,0\}$.
Par conséquent, pour tout $\psi$ dans ${\rm
L}^2\left(\bar{\Gamma},\mu\right)$, la mesure
$\lambda_{\bar{P}\varphi,\psi}$ est atomique et concentrée sur
l'ensemble $\bigcup_{n\in\mathbb N}f^{-n}(3)\cup\bigcup_{n\in\mathbb
N}f^{-n}(0)$. Le résultat en découle.
\end{demo}

\section{Le graphe de Sierpi{\'n}ski}
\label{PascalSierp}

Dans cette section, nous expliquons rapidement comment les résultats
obtenus dans cet article pour le triangle de Pascal $\Gamma$ se
transportent au graphe de  Sierpi{\'n}ski $\Theta$ représenté à la
figure \ref{graphes}. Comme on l'a vu à la section
\ref{decritgraphes}, le graphe $\Theta$ s'identifie au graphe des
arêtes de $\Gamma$. Si $\varphi$ est une fonction sur $\Gamma$, on
note $\Xi^*\varphi$ la fonction sur $\Theta$ telle que, pour tous
points voisins $p$ et $q$ de $\Gamma$, la valeur de $\Xi^*\varphi$
sur l'arête associée à $p$ et $q$ soit $\varphi(p)+\varphi(q)$. On
note $\Xi$ l'adjoint de $\Xi^*$ et on vérifie immédiatement le

\begin{Lem}\label{relationSierp} On a
$(\Delta-1)\Xi^*=\Xi^*\Delta$ et $\Xi\Xi^*=3+\Delta$. La restriction
de $\Delta$ à l'orthogonal de l'image de $\Xi^*$ dans
$\ell^2(\Theta)$ est une homothétie de rapport $-2$.
\end{Lem}

\`A travers ce lemme, l'ensemble des résultats de cet article se
transportent du graphe de Pascal au graphe de Sierpi{\'n}ski. Ils
pourraient d'ailleurs s'obtenir directement dans le graphe de
Sierpi{\'n}ski, en considérant les opérateurs adéquats dans
$\ell^2(\Theta)$. Nous nous contenterons ici de décrire le spectre
continu de $\Theta$ et de traduire le théorème \ref{spectrePascal} :
ceci répond à la question posée par Teplyaev dans \cite[§ 6.6]{Tep}.

Pour $x$ dans $\mathbb R$, posons $k(x)=x+2$ et $t(x)=x+1$. Du lemme
\ref{relationSierp}, on déduit le

\begin{Lem}\label{mesurePascalSierp}
Soient $\varphi$ dans $\ell^2(\Gamma)$, $\mu$ la mesure spectrale de
$\varphi$ pour $\Delta$ dans $\ell^2(\Gamma)$ et $\lambda$ la mesure
spectrale de $\Xi^*\varphi$ pour $\Delta$ dans $\ell^2(\Theta)$.
Alors, on a $\lambda=k(t_*\mu)$.\end{Lem}

Pour tout $x$ dans $\mathbb R$, on pose $g(x)=x^2-3x=f(x-1)+1$. On
note $\Sigma=t(\Lambda)$ l'ensemble de Julia de $g$. Pour tout $x$
dans $\mathbb R$, soit $c(x)=(x+2)(4-x)=k(x)h(x-1)$ et, pour
$x\neq\frac{3}{2}$, $\gamma(x)=\frac{x-1}{2x-3}=\rho(x-1)$. On note
$\nu_\gamma=t_*\nu_\rho$ l'unique mesure de probabilité
$L_{g,\gamma}$-invariante sur $\Sigma$.

Notons toujours $\varphi_0$ la fonction sur $\Gamma$ apparaissant à
la section \ref{decompo} et posons $\theta_0=\Xi^*\varphi_0$ (c'est
la fonction notée $1_{\partial\partial V}$ dans \cite[§ 6]{Tep}). Du
théorème \ref{spectrePascal} et des lemmes \ref{relationSierp} et
\ref{mesurePascalSierp}, on déduit le théorème suivant, qui complète
la description du spectre de $\Theta$ effectuée par Teplyaev dans
\cite{Tep} :

\begin{Thm} Le spectre de $\Delta$ dans $\ell^2(\Theta)$ est
constitué de la réunion de $\Sigma$ et de l'ensemble
$\bigcup_{n\in\mathbb N}g^{-n}(-2)$. La mesure spectrale de
$\theta_0$ pour $\Delta$ dans $\ell^2(\Theta)$ est $c\nu_\gamma$,
les valeurs propres de $\Delta$ dans $\ell^2(\Theta)$ sont les
éléments de $\bigcup_{n\in\mathbb
N}g^{-n}(-2)\cup\bigcup_{n\in\mathbb N}g^{-n}(-1)$ et les
sous-espaces propres associés sont engendrés par des fonctions à
support fini. Enfin, l'orthogonal de la somme des sous-espaces
propres de $\Delta$ dans $\ell^2(\Theta)$ est le sous-espace
cyclique engendré par $\theta_0$.\end{Thm}

\noindent Jean-François Quint\\
LAGA\\
Université Paris 13\\
99, avenue Jean-Baptiste Clément\\
93430 Villetaneuse\\
France\\
quint@math.univ-paris13.fr

\end{document}